\documentclass[3p,11pt,review,preprint]{elsarticle}

\usepackage{hyperref}
\usepackage{graphicx,subcaption}
\usepackage[labelfont=bf,figurename=Fig.,labelsep=period]{caption}
\usepackage[usenames, dvipsnames]{color}
\usepackage{bm}
\usepackage{amsmath,amssymb,amsfonts}
\usepackage{mathrsfs}
\usepackage{multirow}
\usepackage{algpseudocode}
\usepackage[inline]{enumitem}
\usepackage{soul}
\usepackage[linesnumbered,ruled,vlined]{algorithm2e}

\begin{document}
\begin{frontmatter}



\title{A non-intrusive reduced-order modelling for uncertainty propagation of time-dependent problems using a B-splines B\'{e}zier elements-based method and Proper Orthogonal Decomposition: application to dam-break flows}
\date{}

\author[1]{Azzedine Abdedou}
\author[1]{Azzeddine Soula\"imani\corref{cor1}}
\cortext[cor1]{Corresponding author. Tel.: +1 514 396 8977; fax: +1 514 396 8530.}
\ead{azzeddine.soulaimani@etsmtl.ca}
\address[1]{Department of Mechanical Engineering, Ecole de Technologie Sup\'{e}rieure, 1100 Notre-Dame W., Montr\'{e}al (QC),
Canada H3C 1K3}

\begin{abstract}
A proper orthogonal decomposition-based B-splines B\'{e}zier elements method (POD-BSBEM) is proposed as a non-intrusive reduced-order model for uncertainty propagation analysis for stochastic time-dependent problems. The method uses a two-step proper orthogonal decomposition (POD) technique to extract the reduced basis from a collection of high-fidelity solutions called snapshots. A third POD level is then applied on the data of the projection coefficients associated with the reduced basis to separate the time-dependent modes from the stochastic parametrized coefficients. These are approximated in the stochastic parameter space using B-splines basis functions defined in the corresponding B\'{e}zier element. The accuracy and the efficiency of the proposed method are assessed using benchmark steady-state and time-dependent problems and compared to the reduced-order model-based artificial neural network (POD-ANN) and to the full-order model-based polynomial chaos expansion (Full-PCE). The POD-BSBEM is then applied to analyze the uncertainty propagation through a flood wave flow stemming from a hypothetical dam-break in a river with a complex bathymetry. The results confirm the ability of the POD-BSBEM to accurately predict the statistical moments of the output quantities of interest with a substantial speed-up for both offline and online stages compared to other techniques.\\  
\end{abstract}
\begin{keyword}
Uncertainty propagation\sep Proper orthogonal decomposition\sep B-Splines B\'{e}zier elements method\sep Dam-break flows



\end{keyword}
\end{frontmatter}


\section{Introduction}\label{intro}
The recent developments in computational power have led to great improvements in numerical modeling, leading to the solutions of more and more complex physical problems in several engineering fields. However, when these problems are related to large scale configurations and involve mathematical models expressed as time-dependent and parametrized governing equations whose discretization may be significantly demanding in terms of the computational efforts required \citep{berzicnvs2020standardized,sanderse2020non}. This situation becomes all the more demanding when the input parameters are subject to uncertainties that propagate through numerical solvers \citep{kalinina2020metamodeling}, and which must be considered, further increasing the requirements both in terms of CPU time and memory demands. To address this issue, considerable efforts have been dedicated over the last decades to develop reduced-order modeling methods (ROM) to approximate the original full-order model based on one with a  significantly reduced dimensionality \citep{zokagoa2018pod,georgaka2020hybrid}.\\

One of the widely-used methods to construct such a reduced-order model is proper orthogonal decomposition (POD) \citep{sirovich1987turbulence,chatterjee2000introduction}, which adopts the singular value decomposition technique (SVD) to extract the reduced basis from a collection of high-fidelity solutions, called a snapshot matrix. Two broad categories can be considered to compute the associated coefficients of a reduced basis, intrusive and non-intrusive techniques. The intrusive methods combine the POD with Galerkin's projection and require a modification of the source code that solves the resulting discretized governing equations \citep{zokagoa2012pod,fang2013non,couplet2005calibrated}. Despite the improvements that have been introduced \citep{willcox2006unsteady,ballarin2015supremizer}, intrusive methods may require cumbersome modifications in complex codes, and the inaccessibility of these codes in most commercial software poses another difficulty. As a valuable alternative to overcome these issues, non-intrusive techniques have become increasingly attractive due to their flexibility, given that they treat the source code describing the physical model as a black box. A growing number of studies based on regression interpolations have been proposed to compute the associated coefficients of the reduced basis, including Gaussian process regression (GPR) \citep{guo2018reduced}, Radial basis functions interpolation (RBF)\citep{walton2013reduced,xiao2016non}, Kriging interpolation \citep{12081940}, and Smolyak sparse grid interpolation \citep{XIAO2015522}. Recently, a panoply of configurations dealing with regression techniques based on artificial neural networks (ANN) has been introduced in the ROM construction process to accurately approximate unknown coefficients \citep{hesthaven2018non,zhiwei2020non}. Other non-intrusive approaches have been developed to deal with parametrized time-dependent problems that require specific treatments \citep{guo2019data,swischuk2019projection}. Some of them use the two-level POD technique to extract the reduced basis from a large snapshot matrix \citep{wang2019non}, more particularly when the associated parametric values are generated by sampling in the random parametric domain \citep{JACQUIER2021109854}.\\

The effect of the variability stemming from the input parameters has to be quantified in the output quantities of interest through the uncertainty propagation analysis framework. The sampling approaches, i.e. Monte Carlo (MC) and Latin Hypercube Sampling (LHS) \citep{mckay1979comparison} are the most popular and are known to require high sample sizes to reach an accurate estimation of the statistics. Being able to combine them with the above-mentioned reduced-order models renders them a useful tool for uncertainty propagation analysis, and with lower cost. In addition to the sampling methods mentioned above, polynomial chaos expansions (PCEs) \citep{ghanem1991stochastic,sudret2014polynomial} are known to present alternative techniques to deal with uncertainty propagation analysis with a straightforward method to approximate the statistics. They are mainly associated with POD to propose a non-intrusive stochastic reduced-order model for polynomial chaos representation \citep{hijazi2020non,raisee2015non,raisee2019non}. Recently, the so-called POD-PCE methods \citep{sun2019non,el2020stochastic}, which are based on an offline-online paradigm, have been introduced as a stochastic reduced-order model, in which the projection data of the associated coefficient of the reduced basis is expressed as a stochastic expansion with respect to the orthogonal polynomials of the random input parameters. The same paradigm has been extended by Sun and Choi \citep{sun2021non} to space-time-dependent parametrized problems where the projection time-parameter-dependent coefficients are approximated by POD-PCE, combining the generated temporal modes with the orthonormal polynomial basis functions in the random parametric domain. Despite the multiple advantages they offer, these approaches may present some limitations to accurately reproduce the stochastic outputs that present strong hyperbolic behavior or even discontinuity where spurious oscillations may result, as is the case for dam-break flows.\\
 
This work presents a non-intrusive reduced order model for stochastic time-dependent problems. The method, referred to as the POD-BSBEM, combines the advantages of the POD and the B-splines B\'{e}zier element method (BSBEM) \citep{abdedou2019non,abdedou2020uncertainty} to extract the reduced basis and to build a stochastic expansion over the random space, respectively. A two-level POD is applied to the snapshot matrix obtained by the collection of time-trajectory high-fidelity solutions for each given value of the random parametric vector. A third POD level is then performed for each associated coefficient of the reduced basis to separate the time-dependent modes from the stochastic parameterized weights. A series of appropriate transformations are then applied on the resulting projection data, considered as outputs, with which a local stochastic expansion is constructed based on local B-splines basis functions within each B\'{e}zier element constituting the parametric domain. Despite the similarity with the method proposed by Sun and Choi \citep{sun2021non} in the way decoupling the time-parameter dependency for each projection coefficient through an additional level of POD, the technique proposed in this work presents a notable difference in the approximation of the parameter-dependent coefficients by using a multi-elements approach based on local piecewise B-splines basis functions in each B\'{e}zier element in contrast to the polynomial chaos basis functions which cover the whole parametric domain. The use of such local interpolations with the multi-element aspect offers a promising tool to accurately estimate the statistics of the output quantities of interest, particularly  those with a strong hyperbolic or even discontinuous behavior.\\

The paper is organized as follows. Section 2 presents the fundamental concepts of the POD (subsection \ref{POD}) and the mathematical framework of the proposed POD-BSBEM (subsection \ref{BSBEM}). Section \ref{res} evaluates the performance and accuracy of the proposed method through two benchmark numerical test cases for the stochastic steady-state Ackley function (subsection \ref{Ack_sol}), the time-dependent Burgers' equation (subsection \ref{Burg_sol}), and for its application to a hypothetical dam-break with a real terrain database (subsection \ref{dam_break}). Finally, a summary and concluding remarks are presented in Section \ref{conc}.

\section{Methodology}\label{metho}
This section presents  the main aspects of the reduced-order model formulation based on the proper orthogonal decomposition concept combined with the non-intrusive B-splines B\'{e}zier element method (POD-BSBEM) for both steady-state and time-dependent problems with uncertain input parameters.\\

Consider a parametrized time-dependent relationship $u(\mathbf{x},\boldsymbol{\eta},t)=\mathcal{M}(\boldsymbol{\eta},t)$ that links an input of uncertain parameters to the output response, where $\mathcal{M}$ stands for the computational model that may be considered as a black box, $\mathbf{x}\in\Omega\subset\mathbb{R}^{q}$,  $q=1,\,2$ or $3$ is the space domain composed of $N_{e}$ computational elements, and $t\in\mathcal{T}=\left[0,\;T \right]$ represents the temporal domain which is decomposed into $N_{t}$ time steps. $\boldsymbol{\eta}\in\mathcal{P}\subset\mathbb{R}^{m}$ denotes the parametric domain that contains a random vector $\boldsymbol{\eta}=\{\eta_{1},\eta_{2},\ldots,\eta_{m}\}$, with $m$ independent input uncertain parameters described by their probability density functions $ \varrho_{i}(\eta_{i})$ on the probability space $\left( \Theta, \varSigma, P \right)$, where $\Theta$ is the event space, $ P $ the probability measure and $\varSigma$ is the $ \sigma$-algebra on $\Theta$.
\subsection{Proper Orthogonal Decomposition}\label{POD}
The main framework of reduced-order models is based on the approximation of the aforementioned high-fidelity solution of parameterized time-dependent problems by a combination of modes in the reduced-order basis stemming from the proper orthogonal decomposition (POD) as follows:
\begin{equation}\label{<eq_u_pod>}
 \mathbf{u}(\mathbf{x},\boldsymbol{\eta},t)\approx\widehat{\mathbf{u}}(\mathbf{x},\boldsymbol{\eta},t)=\sum_{\ell=1}^{L}\beta_{\ell}(\boldsymbol{\eta},t)\Phi_{\ell}(\mathbf{x})
\end{equation}
where $\Phi_{\ell}(\mathbf{x})\;(\ell=1,\ldots,L)$ represents the orthonormal modes that define the low dimensional basis in which $L$ is the total number of modes. $\beta_{\ell}(\boldsymbol{\eta},t)$ denotes the coefficients of the $\ell^{th}$ mode that are dependent upon time and the random parameter vector $\boldsymbol{\eta}$. As outlined above, the extraction of the low dimensional basis is performed with the proper orthogonal decomposition procedure (POD) \citep{liang2002proper,hesthaven2018non} from a snapshot matrix $\boldsymbol{\mathcal{U}}\in\mathbb{R}^{N_{e}\times N_{s}N_{t}}$, a collection of high-fidelity solutions vectors, called snapshots \citep{sirovich1987turbulence}, generated by the numerical solver $\left\lbrace \mathbf{u}(\boldsymbol{\eta}^{(s)},t_{j})\in\mathbb{R}^{N_{e}\times 1};\,j=1,\ldots,N_{t};\,s=1,\ldots,N_{s}\right\rbrace$, with $N_{s}$ denotes the sample size. In the present work, the aforementioned global snapshot matrix $\boldsymbol{\mathcal{U}}$ is organized by concatenating local snapshot matrices $\mathcal{U}_{s}\in\mathbb{R}^{N_{e}\times N_{t}}$ that collect the high fidelity solutions for a given value of the random parametric vector $\boldsymbol{\eta^{(s)}},\;s=1,\ldots,N_{s}$ over all of the time-steps $t=t_{1},\ldots,t_{N_{t}}$. It is expressed as \citep{guo2019data}:
\begin{equation}\label{<eq_loc_snap_mat>}
 \mathcal{U}_{s}=\left[\mathbf{u}(\boldsymbol{\eta}^{(s)},t_{1})\mid\ldots\mid\mathbf{u}(\boldsymbol{\eta}^{(s)},t_{N_{t}}) \right]\in\mathbb{R}^{N_{e}\times N_{t}} 
\end{equation}
The global snapshot matrix assembling all the snapshots is defined as:
\begin{equation}\label{<eq_glob_snap_mat>}
 \boldsymbol{\mathcal{U}}=\left[\;\mathcal{U}_{1}\mid\ldots\mid\mathcal{U}_{s}\mid\ldots\mid\mathcal{U}_{N_{s}} \right]\in\mathbb{R}^{N_{e}\times N_{s}N_{t}} 
\end{equation}
The reduced basis can be recovered with POD using the singular value decomposition (SVD) of the snapshot matrix $\boldsymbol{\mathcal{U}}$ \citep{wang2019non}:
\begin{equation}\label{<eq_svd>}
 \boldsymbol{\mathcal{U}}=\mathbb{W}\begin{bmatrix} \mathbb{D} & 0 \\ 0 & 0 \end{bmatrix}
 \mathbb{V}^{T} 
\end{equation}
where $\mathbb{W}=\left[\;\mathcal{W}_{1}\mid\ldots\mid\mathcal{W}_{N_{e}} \right]\in\mathbb{R}^{N_{e}\times N_{e}}$ and $\mathbb{V}=\left[\;\mathcal{V}_{1}\mid\ldots\mid\mathcal{V}_{N_{s}N_{t}} \right]\in\mathbb{R}^{N_{s}N_{t}\times N_{s}N_{t}}$ are orthogonal matrices, and $\mathbb{D}=diag\,(\sigma_{1},\ldots,\sigma_{N_{r}})$ is the diagonal matrix containing the $N_{r}$ singular values $\sigma_{1}\geq\sigma_{2}\geq\ldots\geq\sigma_{N_{r}}>0$, in which  $N_{r}\leq min(N_{e},N_{s}N_{t})$ represents the number of non-zero singular values. The low-rank POD basis can be constructed by selecting the $L$-most dominant singular values, i.e., those that capture the essential of the cumulative energy, by adopting the following criterion \citep{sun2019non}:
\begin{equation}\label{<eq_epsi_pod>}
\frac{\sum_{\ell=1}^{L}\sigma_{\ell}^2}{\sum_{\ell=1}^{N_{r}}\sigma_{\ell}^2}>1-\epsilon
\end{equation}
where  $\epsilon$ is a user-defined hyperparameter representing the relative error tolerance. Thus, the POD mode vector $\Phi_{\ell}\in\mathbb{R}^{N_{e}\times 1}$ can be obtained by combining the snapshot matrix $\boldsymbol{\mathcal{U}}$ and the $\ell^{th}$ column of  $\mathbb{V},\;\mathcal{V_{\ell}}\in\mathbb{R}^{N_{s}N_{t}\times 1}$ :
\begin{equation}\label{<eq_pod_mode>}
\Phi_{\ell}=\frac{\boldsymbol{\mathcal{U}}\mathcal{V_{\ell}}}{\sigma_{\ell}} 
\end{equation}
The orthonormal reduced POD  basis obtained above   $\boldsymbol{\Phi}=\left[\Phi_{1}\mid\ldots\mid\Phi_{L} \right]\in\mathbb{R}^{N_{e}\times L}=POD\,(\boldsymbol{\mathcal{U}},\epsilon)$ provides an optimal low-dimensional approximation of the snapshot data.\\

The POD construction procedure described above (Eqs.\eqref{<eq_svd>}, \eqref{<eq_epsi_pod>} and \eqref{<eq_pod_mode>}) for parameterized time-dependent problems with large computational domains and many time-steps will likely be cumbersome and thus incur significant computational costs and require too much memory to perform the SVD of the collected snapshot matrix. To overcome this problem, the so-called two-step POD algorithm \citep{wang2019non} is adopted by performing the first-level POD on the time-trajectory matrix $ \mathcal{U}_{s}\in\mathbb{R}^{N_{e}\times N_{t}}$ given by Eq.\eqref{<eq_loc_snap_mat>}, corresponding to the $s^{th}$ value of the random parametric vector sampling set $\boldsymbol{\eta^{(s)}}$, i.e., $\mathbb{T}_{s}=\left[T_{1}\mid\ldots\mid T_{L_{s}} \right]\in\mathbb{R}^{N_{e}\times L_{s}}=POD\,(\mathcal{U}_{s},\epsilon_{t})$. The compressed time-trajectory matrices thereby obtained for the whole parametric sampling set components ($s=1,\ldots,N_{s}$) are concatenated, and then a second POD level is applied on the resulting matrix to extract the final POD basis, expressed as   $\boldsymbol{\Phi}=\left[\Phi_{1}\mid\ldots\mid\Phi_{L} \right]\in\mathbb{R}^{N_{e}\times L}=POD\left( \left[\mathbb{T}_{1}\mid\ldots\mid\mathbb{T}_{N_{s}}\right]\in\mathbb{R}^{N_{e}\times \sum_{s=1}^{N_{s}}L_{s}},\epsilon_{s}\right)$, where $\epsilon_{t}$ and $\epsilon_{s}$ represent the relative error tolerance for time and random parametric spaces, respectively.  Detailed algorithms with illustrating schematics of the two-step POD procedure can be found in \citep{wang2019non,JACQUIER2021109854}.
\subsection{Regression-based non-intrusive B-Splines B\'{e}zier elements method}\label{BSBEM}
Once the construction of the reduced POD basis has been performed, the corresponding parametrized time-dependent POD coefficients involved in Eq.\eqref{<eq_u_pod>} must be computed. The B-splines B\'{e}zier element method is combined with the POD model reduction technique to propose a novel stochastic non-intrusive reduced order-based regression approach (POD-BSBEM) to accurately estimate the statistics of the output response of quantities of interest with a low computational cost. This approach  belongs to the multi-element techniques that use local piecewise bases  expressed as a function of input random parameters. It is worth noting that the basic framework of the BSBEM in its full-order version has been presented in detail in \citep{abdedou2019non,abdedou2020uncertainty}. Thus, only a summary description of its fundamental aspects is outlined in the following.\\

Let us consider $\boldsymbol{\xi}=\{\xi_{1},\xi_{2},\ldots,\xi_{m}\} \in \Gamma=\left[0,\,1\right]^{m}$ as an ensemble of $m$ independent parametric variables defined in the parametric domain and linked to their corresponding sets of input random parameters in the physical domain $\boldsymbol{\eta}=\{\eta_{1},\eta_{2},\ldots,\eta_{m}\}\in \mathbb{R}^{m}$ by the marginal cumulative density function (CDF) whose image domain corresponds to that of the B-splines basis functions, i.e., $\left[0,\,1\right]$. This allows $\xi_{i}(\eta_{i})=F_{i}(\eta_{i})=\int_{-\infty}^{\eta_{i}}\varrho_{i}(\eta')d\eta'$, which implies that $ d\xi_{i}=\varrho_{i}(\eta_{i})d\eta_{i} $. The parametric domain of each variable $\xi_{i}\in\left[0,\,1 \right] $ is decomposed into $nx_{i}$ non-overlapping segments. The tensor product of these univariate intervals forms a multidimensional element defined as a B\'{e}zier element $I_{e}\subset\Gamma=\left[0,\,1\right]^{m}$ whose total number in the entire parametric domain is given by $N_{elt}=\prod_{i=1}^{m}nx_{i}$. A set of $n_{b}=\prod_{i=1}^{m}(p_{i}+1)$ multivariate local piecewise basis functions is constructed within each B\'{e}zier element using a tensor product of the univariate basis, expressed as  
\begin{equation}\label{<eq_bsbem_multi_basis>}
 \Psi_{j}^{e}(\boldsymbol{\xi})=\psi_{1}^{e}(\xi_{1})\times\ldots\times\psi_{m}^{e}(\xi_{m}),\; j=1,\ldots,n_{b}
\end{equation}
where $e$ stands for the $e^{th}$ B\'{e}zier element $I_{e}$ and $p_{i}$ denotes the polynomial order of the univariate B-splines basis function $\psi_{i}(\xi_{i})$. Details of the construction of the B-splines basis and their extraction from the local Bernstein's functions can be found in \citep{abdedou2019non,abdedou2020uncertainty,borden2011isogeometric}. In each B\'{e}zier element, a set of $n_{s}^{e}\geq n_{b}$ collocation points are selected $\boldsymbol{\xi}^{e(r)}\in I_{e};\;r=1,\ldots,n_{s}^{e}$; $e=1,\ldots,N_{elt}$ and used to construct the snapshot matrix by running the numerical solver with their corresponding values in the physical domain $\left\lbrace \mathbf{u}(\boldsymbol{\eta}^{e(r)},t_{j})\in\mathbb{R}^{N_{e}\times 1};\,j=1,\ldots,N_{t};\,r=1,\ldots,n_{s}^{e}\right\rbrace$, obtained by an appropriate inverse CDF mapping:   
\begin{equation} \label{<eq_cdf_mapping>}
\eta_{i}^{e\,(r)}=F_{i}^{-1}(\xi_{i}^{e\,(r)}),\,r=1,\ldots,n_{s}^{e},\,i=1,\ldots,m
\end{equation}
where $F_{i}^{-1}$ denotes the cumulative density function that links the random parametric variable $\xi_{i}$ to its corresponding random physical variable $\eta_{i}$. The total number of sampling points over the entire parametric domain is $N_{s}=n_{s}^{e}N_{elt}$.\\
\subsubsection{A third level proper orthogonal decomposition}
The unknown parametric time-dependent coefficients associated with the orthonormal column bases of $\boldsymbol{\Phi}$, involved in Eq.\eqref{<eq_u_pod>}, can be expressed as: 
\begin{equation}\label{<eq_coeff_B_matrix>}
\boldsymbol{\mathcal{B}}(\boldsymbol{\eta},t)=\boldsymbol{\Phi}^{T}\boldsymbol{\mathcal{U}}=\begin{bmatrix} 
\beta_{11}^{1} & \cdots & \beta_{1N_{t}}^{1}\mid & \cdots &\mid  \beta_{11}^{N_{s}} & \cdots & \beta_{1N_{t}}^{N_{s}} \\ 
\vdots & \ddots & \vdots &  & \vdots & \ddots & \vdots & \\
\beta_{\ell1}^{1} & \cdots & \beta_{\ell N_{t}}^{1}\mid & \cdots &\mid \beta_{\ell1}^{N_{s}} & \cdots & \beta_{\ell N_{t}}^{N_{s}}\\
\vdots & \ddots & \vdots &  & \vdots & \ddots & \vdots \\
\beta_{L1}^{1} & \cdots & \beta_{L N_{t}}^{1}\mid & \cdots &\mid  \beta_{L1}^{N_{s}} & \cdots & \beta_{L N_{t}}^{N_{s}}
 \end{bmatrix}\in\mathbb{R}^{L\times N_{s}N_{t}}        
\end{equation}
It is worth noting that the procedure adopted to compute these coefficients is detailed in the following paragraphs  only for the $\ell^{th}$ projection coefficient, i.e., $\boldsymbol{\beta}_{\ell}(\boldsymbol{\eta},t)\in\mathbb{R}^{1\times N_{s}N_{t}}$, representing the $\ell^{th}$ row vector reported in Eq.\eqref{<eq_coeff_B_matrix>}, and will be generalized for the whole $L$ modes. Thus, this $\ell^{th}$ row vector can be written in a matrix form, $\boldsymbol{\beta}_{\ell}\in\mathbb{R}^{1\times N_{s}N_{t}}\equiv\boldsymbol{\beta}^{\ell}\in\mathbb{R}^{N_{t}\times N_{s}}$, with $N_{t}$ rows and $N_{s}$ columns, representing the time and the parameter locations, respectively, as \citep{guo2019data}:
\begin{equation}\label{<eq_trans_line_vec_to_matrix>}
\boldsymbol{\beta}^{\ell}(t,\boldsymbol{\eta})=\begin{bmatrix} 
\beta_{11}^{\ell} & \cdots & \beta_{1N_{s}}^{\ell}\\ 
\vdots & \ddots & \vdots \\
\beta_{N_{t}1}^{\ell} & \cdots & \beta_{N_{t}N_{s}}^{\ell}\\
 \end{bmatrix}\in\mathbb{R}^{N_{t}\times N_{s}}        
\end{equation}
A third POD level is then applied on $\boldsymbol{\beta}^{\ell}(t,\boldsymbol{\eta})$ to decompose the data into time-dependent modes with their associated parameter-dependent coefficients, expressed as:
\begin{equation}  \label{<eq_third_Pod>}
\beta^{\ell}(t,\boldsymbol{\eta})\approx\sum_{k=1}^{\mathcal{K}_{\ell}}\lambda_{k}^{\ell}(\boldsymbol{\eta})\chi_{k}^{\ell}(t)
\end{equation}
where $\boldsymbol{X}^{\ell}(t)=\left[\boldsymbol{\chi}^{\ell}_{1}\mid\ldots\mid\boldsymbol{\chi}^{\ell}_{\mathcal{K}_{\ell}}\right]\in\mathbb{R}^{N_{t}\times \mathcal{K}_{\ell}}=POD\left(\boldsymbol{\beta}^{\ell}(t,\boldsymbol{\eta})\in\mathbb{R}^{N_{t}\times N_{s}},\epsilon_{s}\right)$ and $\boldsymbol{\Lambda}^{\ell}(\boldsymbol{\eta})\in\mathbb{R}^{\mathcal{K}_{\ell}\times N_{s}}$ represent the time-dependent orthogonal basis modes  with a truncating rank of $\mathcal{K}_{\ell}$ and the associated parameter-dependent coefficients, respectively, of the $\ell^{th}$ mode. These associated coefficients can be expressed in a matrix form obtained by projection from Eq.\eqref{<eq_third_Pod>} by exploiting the orthonormality feature of the $\boldsymbol{X}^{\ell}$ basis:
\begin{equation} \label{<eq_lamda_coeff>}
\boldsymbol{\Lambda}^{\ell}=\boldsymbol{X}^{\ell^{T}}\boldsymbol{\beta}^{\ell}=\begin{bmatrix} 
\lambda_{11}^{\ell}&\cdots&\lambda_{1N_{s}}^{\ell}\\ 
\vdots & \ddots & \vdots \\
\lambda_{k1}^{\ell}&\cdots&\lambda_{kN_{s}}^{\ell}\\
\vdots & \ddots & \vdots\\
\lambda_{\mathcal{K}_{\ell}1}^{\ell}&\cdots&\lambda_{\mathcal{K}_{\ell}N_{s}}^{\ell}\\
 \end{bmatrix}\in\mathbb{R}^{\mathcal{K}_{\ell}\times N_{s}}        
\end{equation}
The procedure used to construct a regression as a function of stochastic basis functions is detailed below. As with the above coefficient calculation, the calculations are restricted, in this case, to the $k^{th}$ row vector, and must be generalized for the $\mathcal{K}_{\ell}$ associated coefficients given by Eq.\eqref{<eq_lamda_coeff>}. The collocation points are selected locally within each B\'{e}zier element $I_{e}$ with a minimum required number of $n_{s}^{e}$, as noted earlier. The set of sampling points $N_{s}=n_{s}^{e}N_{elt}$ is obtained by concatenating the local points $n_{s}^{e}$ over the complete set of $N_{elt}$ B\'{e}zier elements that contains entire parametric domain. Thus, the $k^{th}$ row vector of $\boldsymbol{\Lambda}^{\ell}$ is a concatenation of local vectors corresponding to each B\'{e}zier element, i.e., $\boldsymbol{\lambda}_{k}^{\ell,e}\in\mathbb{R}^{1\times n_{s}^{e}}\subset \boldsymbol{\lambda}_{k}^{\ell}\in\mathbb{R}^{1\times N_{s}}$, where $e=1,\ldots,N_{elt}$ denotes the $e^{th}$ B\'{e}zier element. These local data $\boldsymbol{\delta}_{k}^{\ell,e}\in\mathbb{R}^{n_{s}^{e}\times 1}=(\boldsymbol{\lambda}_{k}^{\ell,e})^{T}$ are considered as  outputs with which a regression is constructed in each B\'{e}zier element as the sum of the local piecewise basis functions:
\begin{equation}\label{<eq_bsbem>}
 \delta_{k}^{\ell,e}=\sum_{j=1}^{nb}\alpha_{j,k}^{\ell,e}\Psi_{j}^{e}(\boldsymbol{\xi}) 
\end{equation}
where $\alpha_{j,k}^{\ell,e}$ and $\Psi_{j}^{e}$ denote the local coefficients and the $n_{b}$ local multivariate piecewise B-splines basis functions, given by Eq.\eqref{<eq_bsbem_multi_basis>}, over the $e^{th}$ B\'{e}zier element $I_{e}$, respectively. These local basis functions are evaluated on each component $\boldsymbol{\xi}^{e(i)}$ of the collocation points set in the $e^{th}$ B\'{e}zier element, i.e. $\boldsymbol{\psi}_{n_{s}^{e}\times nb}^{e}=\left\lbrace \Psi_{i,j}^{e}=\Psi_{j}^{e}(\boldsymbol{\xi}^{e\,(i)}),i=1,\ldots,n_{s}^{e},\,j=1,\ldots,nb \right\rbrace$. Thus, the resulting local design matrix and the output response vector are expressed as:
\begin{equation} \label{<eq_psi_loc>}
\boldsymbol{\Psi}^{e}=\boldsymbol{\psi}^{e^{T}}\boldsymbol{\psi}^{e}\in\mathbb{R}^{nb\times nb}\quad \textrm{and}\quad \boldsymbol{\Delta}^{\ell,e}_{k}=\boldsymbol{\psi}^{e^{T}}\boldsymbol{\delta}_{k}^{\ell,e}\in\mathbb{R}^{nb\times 1} 
\end{equation}
The unknown coefficients are computed by solving a global system of equations instead of the constructed local systems in each B\'{e}zier element. This global system of equations is constructed from the assembly of local design matrices and local output vectors through the assembly operator \citep{hughes2012finite}. The global design matrix and output vector are thus given by:
\begin{equation} \label{<eq_Assembly>}
\boldsymbol{\Psi}=\overset{N_{elt}}{\underset{e=1}{\mathbf{A}}}(\boldsymbol{\Psi}^{e}) \in\mathbb{R}^{M\times M} \quad \textrm{and}\quad  \boldsymbol{\Delta}_{k}^{\ell}=\overset{N_{elt}}{\underset{e=1}{\mathbf{A}}}(\boldsymbol{\Delta}^{\ell,e}_{k}) \in\mathbb{R}^{M\times 1}
\end{equation}
where $M=\prod_{i=1}^{m}{(nx_{i}+p_{i})}$ denotes the total number of the global coefficients to be computed. The components of the global design matrix $\Psi_{q,f}$ are linked to their corresponding components in the local design matrices $\Psi_{i,j}^{e}$ by the connectivity array ($IEN$) with respect to the $e^{th}$ Bézier element: $q=IEN(i,e)$ and $f=IEN(j,e)$, with $q,\,f=1,\ldots,M$ and $i,\,j=1,\ldots,nb$ (more implementation details can be found in \citep{cottrell2009isogeometric,scott2011isogeometric}). Thus, the coefficients can be computed by solving the following global system of equations:
\begin{equation} \label{<eq_glob_coeff>}
\boldsymbol{\alpha}_{k}^{\ell}=\boldsymbol{\Psi}^{-1}\boldsymbol{\Delta}_{k}^{\ell} \in\mathbb{R}^{M\times 1}
\end{equation}
Once these coefficients are computed, the stochastic output response, given by Eq.\eqref{<eq_u_pod>}, can be locally approximated $\widehat{\mathbf{u}}^{e}$ in each B\'{e}zier element, which represents a part of the random parametric space, using  the online stage of the proposed POD-BSBEM. Thus, let us consider a new set of $n_{s}^{e}$ collocation points in the $e^{th}$ B\'{e}zier element, $\boldsymbol{\widehat{\xi}}^{e^{(i)}},\,i=1,\ldots,n_{s}^{e}$, with which the local B-splines basis functions $\left\lbrace \Psi_{j}^{e}\right\rbrace_{j=1}^{n_{b}}$ are evaluated to build the so-called local design matrix: $\boldsymbol{\widehat{\psi}}_{n_{s}^{e}\times nb}^{e}=\left\lbrace \Psi_{i,j}^{e}=\Psi_{j}^{e}\left( \boldsymbol{\widehat{\xi}}^{e\,(i)}\right),i=1,\ldots,n_{s}^{e},\,j=1,\ldots,nb \right\rbrace$. The local coefficients $\boldsymbol{\alpha}_{k}^{\ell,e}\in\mathbb{R}^{n_{b}\times 1}$ are extracted from the set of computed global coefficients, $\boldsymbol{\alpha}_{k}^{\ell}\in\mathbb{R}^{M\times 1}$, with the same connectivity array procedure mentioned above. The approximated local projection data associated with the $\mathcal{K}_{\ell}$ time-dependent modes can be expressed as:
\begin{equation} \label{<eq_delta_k_surro>}
\boldsymbol{\widehat{\Lambda}}^{\ell,e}=\left[\boldsymbol{\widehat{\psi}}^{e}\boldsymbol{\alpha}_{1}^{\ell,e}\mid\ldots\mid \boldsymbol{\widehat{\psi}}^{e}\boldsymbol{\alpha}_{k}^{\ell,e}\mid\ldots\mid\boldsymbol{\widehat{\psi}}^{e}\boldsymbol{\alpha}_{\mathcal{K}_{\ell}}^{\ell,e}\right]^{T} \in\mathbb{R}^{\mathcal{K}_{\ell}\times n_{s}^{e}}
\end{equation}
The components of the $\ell^{th}$ projection coefficient can be approximated locally in the $e^{th}$ B\'{e}zier element, by combining the time-dependent basis functions $\boldsymbol{X}^{\ell}\in\mathbb{R}^{N_{t}\times \mathcal{K}_{\ell}}$, as in Eq.\eqref{<eq_third_Pod>}, and the local projection data from Eq.\eqref{<eq_delta_k_surro>}: $\boldsymbol{\widehat{\beta}}^{\ell,e}=\boldsymbol{X}^{\ell}\boldsymbol{\widehat{\Lambda}}^{\ell,e}\in\mathbb{R}^{N_{t}\times n_{s}^{e}}$. The obtained data are then transformed from an $N_{t}\times n_{s}^{e}$ matrix to a $1\times n_{s}^{e}N_{t}$ row vector, i.e., $\boldsymbol{\widehat{\beta}}_{\ell}^{e}\in\mathbb{R}^{1\times n_{s}^{e}N_{t}}\equiv\boldsymbol{\widehat{\beta}}^{\ell,e}\in\mathbb{R}^{N_{t}\times n_{s}^{e}}$, and so the approximated projection coefficients associated with the $L$ modes reduced basis are obtained as:
\begin{equation} \label{<eq_L_B_surro>}
\boldsymbol{\mathcal{\widehat{B}}}^{e}=\begin{bmatrix} 
\boldsymbol{\widehat{\beta}}_{1}^{e}\\ 
\vdots \\
\boldsymbol{\widehat{\beta}}_{\ell}^{e}\\
\vdots \\
\boldsymbol{\widehat{\beta}}_{L}^{e}\\
 \end{bmatrix}\in\mathbb{R}^{L\times n_{s}^{e}N_{t}}\quad \textrm{with}\quad\boldsymbol{\widehat{\beta}}^{e}_{\ell}= \left[\widehat{\beta}_{\ell1}^{1} \cdots  \widehat{\beta}_{\ell N_{t}}^{1}\mid  \cdots \mid \widehat{\beta}_{\ell1}^{n_{s}^{e}}  \cdots  \widehat{\beta}_{\ell N_{t}}^{n_{s}^{e}} \right]\in\mathbb{R}^{1\times n_{s}^{e}N_{t}} 
\end{equation}
Thus, the approximation of the output response can be expressed by the local snapshot matrix constructed within the $e^{th}$ B\'{e}zier element:
\begin{equation} \label{<eq_U_e_surr>}
\boldsymbol{\mathcal{\widehat{U}}}^{e}= \boldsymbol{\Phi}\boldsymbol{\mathcal{\widehat{B}}}^{e}\in\mathbb{R}^{N_{e}\times n_{s}^{e}N_{t}}
\end{equation}
The approximation of the snapshot matrix over the whole random parametric domain is given by:
\begin{equation} \label{<eq_U_surr>}
\boldsymbol{\mathcal{\widehat{U}}}=\left[\boldsymbol{\mathcal{\widehat{U}}}^{1}\mid\ldots\mid \boldsymbol{\mathcal{\widehat{U}}}^{e}\mid\ldots\mid\boldsymbol{\mathcal{\widehat{U}}}^{N_{elt}}\right] \in\mathbb{R}^{N_{e}\times N_{s}N_{t}}
\end{equation}
The statistical moments can be estimated through the constructed surrogate model of the stochastic output response $\widehat{u}$ as follows:
\begin{equation} \label{<eq_pod_bsbem_mean>}
 \mu_{\widehat{u}}=\mathbb{E}\left[\widehat{u}\right]=\int_{\Xi}\widehat{u}\varrho_{\boldsymbol{\eta}}(\boldsymbol{\eta})d\boldsymbol{\eta} 
\end{equation} 
and
\begin{equation}\label{<eq_pod_bsbem_std>}
\sigma_{\widehat{u}}^{2}=\mathbb{E}\left[\widehat{u}^{2}\right]- \mu_{\widehat{u}}^{2} =\left[ \int_{\Xi}\left(\widehat{u}\right)^{2}\varrho_{\boldsymbol{\eta}}(\boldsymbol{\eta})d\boldsymbol{\eta}\right]-\mu_{\widehat{u}}^{2}  
\end{equation}

The integrals involved in the mean and variance expressions can be simplified by exploiting the correspondence between the CDF image domain and the parametric domain of the basis functions $[0,\,1] $. As indicated above, this CDF mapping, i.e., $ \xi_{i}=F_{i}(\eta_{i})$ and $ d\xi_{i}=\varrho_{i}(\eta_{i})d\eta_{i}$, allows a simplified evaluation of the integrals over each B\'{e}zier element definition domain. Thus, the statistical moments of the approximated stochastic output response of a quantity of interest over each computational node for a given simulation time are expressed as:

\begin{equation} \label{<eq_bsbem_mean_Ie>}
 \mu_{\widehat{u}}=\sum_{e=1}^{N_{elt}}\left[\int_{I_{e}}\widehat{u}^{e}d\boldsymbol{\xi}\right] 
\end{equation}
and
\begin{equation} \label{<eq_bsbem_std_Ie>}
\sigma_{\widehat{u}}^{2}=\sum_{e=1}^{N_{elt}}\left[\int_{I_{e}}\left( \widehat{u}^{e}\right)^{2}d\boldsymbol{\xi}\right]-\mu_{\widehat{u}}^{2}  
\end{equation}
where $\widehat{u}^{e}$ represents the spatio-temporal POD-BSBEM local approximation of the stochastic output response within the $e^{th}$ B\'{e}zier element formulated as: 
\begin{equation}\label{<eq_u_pod>}
 \widehat{u}^{e}(\mathbf{x},\boldsymbol{\eta},t)=\sum_{\ell=1}^{L}\left[\;\sum_{k=1}^{\mathcal{K}_{\ell}}\left[\;\sum_{j=1}^{n_{b}}\alpha_{j,k}^{\ell,e}\Psi_{j}^{e}\left(\boldsymbol{\widehat{\xi}}\right)\;\right]\chi_{k}^{\ell}(t)\;\right] \Phi_{\ell}(\mathbf{x})
\end{equation}

Based on the detailed steps presented above, the implementations of the offline and online stages of the proposed POD-BSBEM are summarized in Algorithm.\ref{alg:Offline_PODBSBEM} and Algorithm.\ref{alg:Online_PODBSBEM}, respectively.\\

\begin{algorithm}[!h]
\SetAlgoLined
\SetKw{KwInput}{\textbf{Input}: }
\SetKw{KwOutput}{\textbf{Output}: }
	
\KwInput{$p_{i}$, $nx_{i}$, $\mu_{\eta_{i}}$, $\sigma_{\eta_{i}}$, $\varrho_{i}(\eta_{i})$, $n_{s}^{e}=\prod_{i=1}^{m}(p_{i}+1)$, $N_{elt}=\prod_{i=1}^{m}nx_{i}$, $i=1,\ldots,m$}

\KwOutput{$\boldsymbol{\Phi}\in\mathbb{R}^{N_{e}\times L}$, $\boldsymbol{\alpha}^{\ell}\in\mathbb{R}^{M\times 1}$, $\boldsymbol{X}^{\ell}\in\mathbb{R}^{N_{t}\times \mathcal{K}_{\ell}}$, $\ell=1,\ldots,L$}

Generate the local collocation points: $\boldsymbol{\xi}^{e^{(r)}}\in I_{e},\,r=1,\ldots,n_{s}^{e};\,e=1,\ldots,N_{elt}$

Compute the high fidelity solutions: $\mathbf{u}(\boldsymbol{\eta}^{e^{(r)}},t_{j}),\,j=1,\ldots,N_{t}\;\textrm{with}\;\eta_{i}^{e\,(r)}=F_{i}^{-1}(\xi_{i}^{e\,(r)})$ 

Construct the time-trajectory snapshot matrices $\mathcal{U}_{s}\in\mathbb{R}^{N_{e}\times N_{t}}$ from Eq.\eqref{<eq_loc_snap_mat>} and the global snapshot matrix $\boldsymbol{\mathcal{U}}\in\mathbb{R}^{N_{e}\times N_{s}N_{t}}$ from Eq.\eqref{<eq_glob_snap_mat>}

Compress the time-trajectory basis: $\mathbb{T}_{s}=\left[T_{1}\mid\ldots\mid T_{L_{s}} \right]=POD\,(\mathcal{U}_{s},\epsilon_{t})$, with $s=1,\ldots,N_{s}$   

Extract the $L$ reduced basis functions: $\boldsymbol{\Phi}=\left[\Phi_{1}\mid\ldots\mid\Phi_{L} \right]=POD\left( \left[\mathbb{T}_{1}\mid\ldots\mid\mathbb{T}_{N_{s}}\right],\epsilon_{s}\right)$

Project the POD-associated coefficients: $\boldsymbol{\mathcal{B}}(\boldsymbol{\eta},t)=\boldsymbol{\Phi}^{T}\boldsymbol{\mathcal{U}}\in\mathbb{R}^{L\times N_{s}N_{t}}$ from Eq.\eqref{<eq_coeff_B_matrix>}

\For{$\ell=1:L$}{
 Transform the $\ell^{th}$ row vector to a matrix: $\boldsymbol{\beta}^{\ell}\in\mathbb{R}^{N_{t}\times N_{s}}\equiv\boldsymbol{\beta}_{\ell}\in\mathbb{R}^{1\times N_{s}N_{t}}$ from Eq.\eqref{<eq_trans_line_vec_to_matrix>}
 
 Extract the $\mathcal{K}_{\ell}$ time-dependent basis $\boldsymbol{X}^{\ell}(t)=\left[\boldsymbol{\chi}^{\ell}_{1}\mid\ldots\mid\boldsymbol{\chi}^{\ell}_{\mathcal{K}_{\ell}}\right]=POD\left(\boldsymbol{\beta}^{\ell},\epsilon_{s}\right)$   
 
 Perform the projection $\boldsymbol{\Lambda}^{\ell}=\boldsymbol{X}^{\ell^{T}}\boldsymbol{\beta}^{\ell}\in\mathbb{R}^{\mathcal{K_{\ell}}\times n_{s}^{e}N_{elt}}$  from Eq.\eqref{<eq_lamda_coeff>} 
 
 \For{$k=1:\mathcal{K}_{\ell}$}{
 Extract $\boldsymbol{\delta}_{k}^{\ell}=\left\lbrace\lambda_{ks};\,s=1,\ldots,N_{s} \right\rbrace^{T} \in\mathbb{R}^{N_{s}\times 1}\;\textrm{with}\;N_{s}=n_{s}^{e}N_{elt}$
 
 \For{$e=1:N_{elt}$}{
 Extract the local output projection data $\boldsymbol{\delta}_{k}^{\ell,e}\in\mathbb{R}^{n_{s}^{e}\times 1}\subset\boldsymbol{\delta}_{k}^{\ell}\in\mathbb{R}^{N_{s}\times 1}$
 
 Define the local design matrix $\boldsymbol{\Psi}^{e}=\boldsymbol{\psi}^{e^{T}}\boldsymbol{\psi}^{e}\in\mathbb{R}^{nb\times nb}$ and the local output vector $\boldsymbol{\Delta}^{\ell,e}_{k}=\boldsymbol{\psi}^{e^{T}}\boldsymbol{\delta}_{k}^{\ell,e}\in\mathbb{R}^{nb\times 1}$ according to Eq.\eqref{<eq_psi_loc>}, with $\boldsymbol{\psi}_{n_{s}^{e}\times nb}^{e}=\left\lbrace \Psi_{i,j}^{e}=\Psi_{j}^{e}(\boldsymbol{\xi}^{e\,(i)}),i=1,\ldots,n_{s}^{e},\,j=1,\ldots,nb \right\rbrace$
 
 Assemble the global design matrices: $\boldsymbol{\Psi}=\overset{N_{elt}}{\underset{e=1}{\mathbf{A}}}(\boldsymbol{\Psi}^{e}) \in\mathbb{R}^{M\times M}$ and global output vectors: $\boldsymbol{\Delta}_{k}^{\ell}=\overset{N_{elt}}{\underset{e=1}{\mathbf{A}}}(\boldsymbol{\Delta}^{\ell,e}_{k}) \in\mathbb{R}^{M\times 1}$ according to Eq.\eqref{<eq_Assembly>}  
 }
 Compute the global coefficients:
$\boldsymbol{\alpha}_{k}^{\ell}=\boldsymbol{\Psi}^{-1}\boldsymbol{\Delta}_{k}^{\ell} \in\mathbb{R}^{M\times 1} $ according to Eq.\eqref{<eq_glob_coeff>}
 }
 }
 \caption{Offline stage of the POD-BSBEM for stochastic time-dependent problems}\label{alg:Offline_PODBSBEM}
\end{algorithm}

\begin{algorithm}[!h]
\SetAlgoLined
\SetKw{KwInput}{\textbf{Input}: }
\SetKw{KwOutput}{\textbf{Output}: }
	
\KwInput{$p_{i}$, $nx_{i}$, $\left\lbrace \boldsymbol{\widehat{\xi}}^{e^{(i)}}\right\rbrace_{i=1}^{n_{s}^{e}}$, $\left\lbrace \Psi_{j}\right\rbrace_{j=1}^{n_{b}} $, $\boldsymbol{\Phi}$, $\boldsymbol{\alpha}^{\ell}$, $\boldsymbol{X}^{\ell}$, $\ell=1,\ldots,L$}

\KwOutput{ $\boldsymbol{\mathcal{\widehat{U}}}^{e}$, $\boldsymbol{\mathcal{\widehat{U}}}$, $\mu_{\widehat{u}}$, $\sigma_{\widehat{u}}^{2}$}

\For{$e=1:N_{elt}$}{
Generate a new set of local collocation points: $\boldsymbol{\widehat{\xi}}^{e^{(i)}}\in I_{e},\,i=1,\ldots,n_{s}^{e}$ 

Construct local design matrix $\boldsymbol{\widehat{\psi}}_{n_{s}^{e}\times nb}^{e}=\left\lbrace \Psi_{i,j}^{e}=\Psi_{j}^{e}\left( \boldsymbol{\widehat{\xi}}^{e\,(i)}\right),\,j=1,\ldots,nb \right\rbrace_{i=1}^{n_{s}^{e}}$
 
 \For{$\ell=1:L$}{
 Select the time-dependent matrix $\boldsymbol{X}^{\ell}\in\mathbb{R}^{N_{t}\times \mathcal{K}_{\ell}}$ 
 
 \For{$k=1:\mathcal{K}_{\ell}$}{
 
 Extract local coefficients from the global ones: $\boldsymbol{\alpha}^{\ell,e}_{k}\in\mathbb{R}^{n_{b}\times 1}\subset\boldsymbol{\alpha}^{\ell}_{k}\in\mathbb{R}^{M\times 1}$
 
 Construct the local regression: $\boldsymbol{\widehat{\delta}}_{k}^{\ell,e}=\boldsymbol{\widehat{\psi}}^{e}\boldsymbol{\alpha}^{\ell,e}_{k}\in\mathbb{R}^{n_{s}^{e}\times 1}$ following Eq.\eqref{<eq_bsbem>}   
 }
 Concatenate the approximated output data: $\boldsymbol{\widehat{\Lambda}}^{\ell,e}=\left[\boldsymbol{\widehat{\delta}}_{1}^{\ell,e}\mid\ldots\mid\boldsymbol{\widehat{\delta}}_{\mathcal{K}_{\ell}}^{\ell,e}\right]^{T} \in\mathbb{R}^{\mathcal{K}_{\ell}\times n_{s}^{e}}$
 
 Compute the data of the $\ell^{th}$ projection coefficient: $\boldsymbol{\widehat{\beta}}^{\ell,e}=\boldsymbol{X}^{\ell}\boldsymbol{\widehat{\Lambda}}^{\ell,e}\in\mathbb{R}^{N_{t}\times n_{s}^{e}}$
 
 Transform the obtained matrix to a row vector: $\boldsymbol{\widehat{\beta}}_{\ell}^{e}\in\mathbb{R}^{1\times n_{s}^{e}N_{t}}\equiv\boldsymbol{\widehat{\beta}}^{\ell,e}\in\mathbb{R}^{N_{t}\times n_{s}^{e}}$ 
 }
 Concatenate over the $L$-modes: $\boldsymbol{\mathcal{\widehat{B}}}^{e}=\left[\boldsymbol{\widehat{\beta}}_{1}^{e^{T}}\mid\ldots\mid\boldsymbol{\widehat{\beta}}_{\ell}^{e^{T}}\mid\ldots\boldsymbol{\widehat{\beta}}_{L}^{e^{T}}\right]^{T}\in\mathbb{R}^{L\times n_{s}^{e}N_{t}} $ 
 
 Compute the local approximated snapshot matrix: $\boldsymbol{\mathcal{\widehat{U}}}^{e}= \boldsymbol{\Phi}\boldsymbol{\mathcal{\widehat{B}}}^{e}\in\mathbb{R}^{N_{e}\times n_{s}^{e}N_{t}}$ 
 
 Compute the statistical moments according to Eqs.\eqref{<eq_bsbem_mean_Ie>} and \eqref{<eq_bsbem_std_Ie>} 
 
 }
Concatenate over the $N_{elt}$ B\'{e}zier elements to approximate the global snapshot matrix: $\boldsymbol{\mathcal{\widehat{U}}}=\left[\boldsymbol{\mathcal{\widehat{U}}}^{1}\mid\ldots\mid \boldsymbol{\mathcal{\widehat{U}}}^{e}\mid\ldots\mid\boldsymbol{\mathcal{\widehat{U}}}^{N_{elt}}\right] \in\mathbb{R}^{N_{e}\times N_{s}N_{t}}$  
 \caption{Online stage of the POD-BSBEM for stochastic time-dependent problems}
 \label{alg:Online_PODBSBEM}
\end{algorithm}
In addition to the sampling approaches such as Monte Carlo (MC) and Latin Hypercube Sampling (LHS), the performance and accuracy of the proposed POD-BSBEM are also compared to the non-intrusive reduced-order model-based artificial neural network (POD-ANN) and the full-order polynomial chaos expansion (Full-PCE). It should be noted that a detailed discussion of the development and implementation steps of these methods is beyond the scope of this work, and the reader is encouraged to consult the corresponding references \citep{hesthaven2018non,wang2019non,JACQUIER2021109854} and \citep{abdedou2020uncertainty,hosder2007efficient,hosder2010non}, respectively, where modeling details can be found.
\section{Results and discussion}\label{res}
This section presents the numerical results for two benchmark test cases, a steady-state Ackley function and a time-dependent Burgers equation, along with the results from a hypothetical dam break with real terrain data. These results are then discussed in terms of statistical moments and probability density function profiles of the quantities of interest, as well as the resulting probabilistic inundation mapping.
\subsection{Ackley function test case} \label{Ack_sol}
The so-called Ackley function is considered as a first test function to showcase the performance of the proposed approach. It presents a challenging candidate due to its high irregularity and non-linearity aspects. This function is defined over a 2D space domain $-5\leq x,y \leq 5$, which is discretized with a mesh of $N_{e}=160\times 160$ nodes \citep{raisee2015non}. The stochastic version of the Ackley function is formulated as follows \citep{sun2019non,raisee2015non}:
\begin{equation}\label{<Ack_analy>}
\begin{split}
 u(x,y,\mathbf{\boldsymbol{\xi}})=-20(1+0.1\xi_{3})\left(exp\left[-0.2(1+0.1\xi_{2})\sqrt{0.5(x^{2}+y^{2})} \right] \right)\\
 -exp\left(0.5\left[ cos(2\pi(1+0.1\xi_{1})x+ cos(2\pi(1+0.1\xi_{1})y)  \right] \right)+20+e
\end{split}
\end{equation}
where $\boldsymbol{\xi}=\left[\xi_{1},\xi_{2},\xi_{3} \right]$ denotes a random vector defined by three uncertain input parameters assumed to be described by uniform distributions over $\left[-1,1 \right]^{3}$ and $e\approx\;2.718281$.\\

\begin{figure}[ht!]
  \centering
    \begin{subfigure}[b]{0.49\textwidth}
      \centering
        \includegraphics[width=\textwidth]{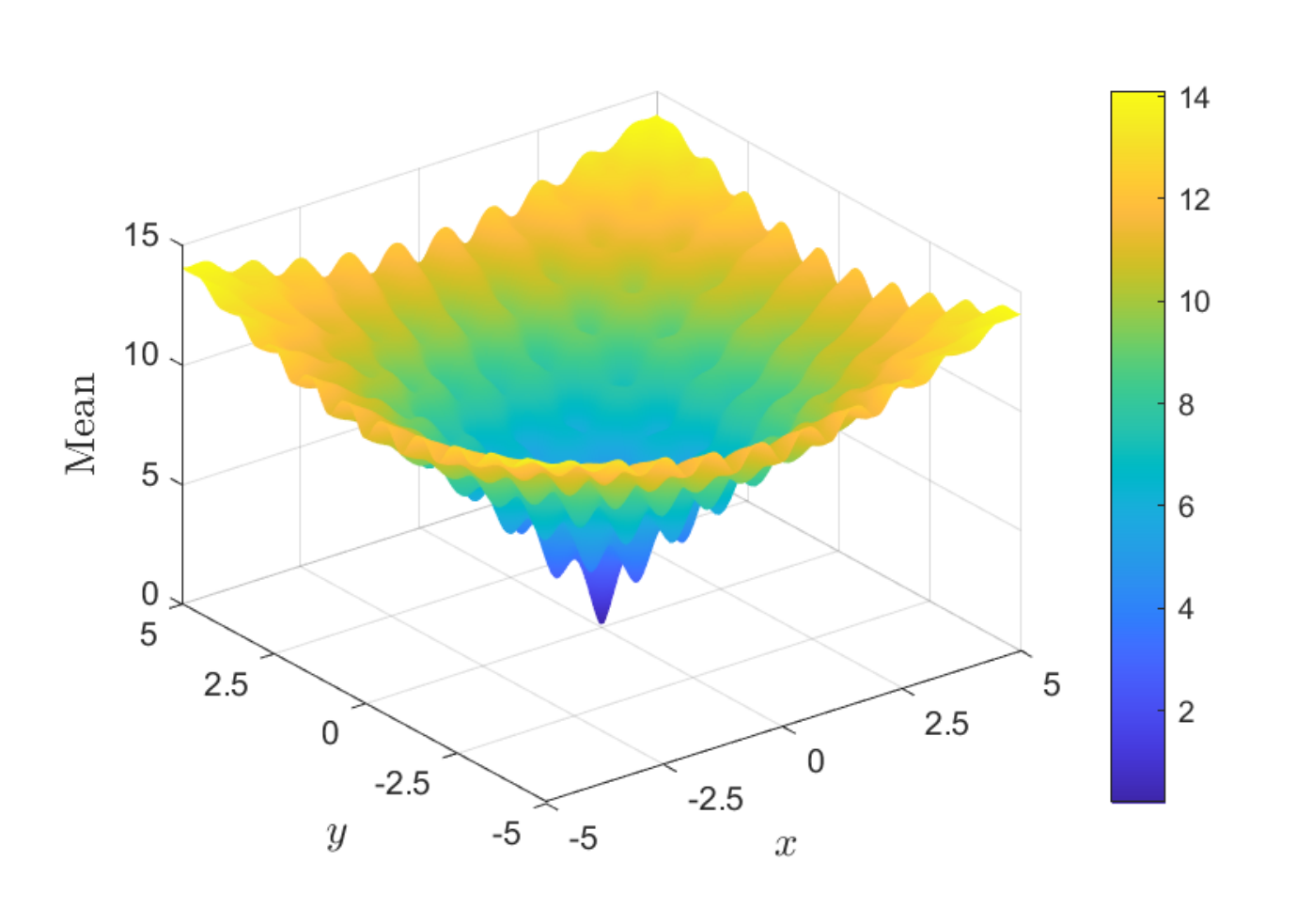}
         \caption{}
         \label{fig:Ack_3D_Mean_lhs}
    \end{subfigure}  
  \hfill
    \begin{subfigure}[b]{0.49\textwidth}
      \centering
        \includegraphics[width=\textwidth]{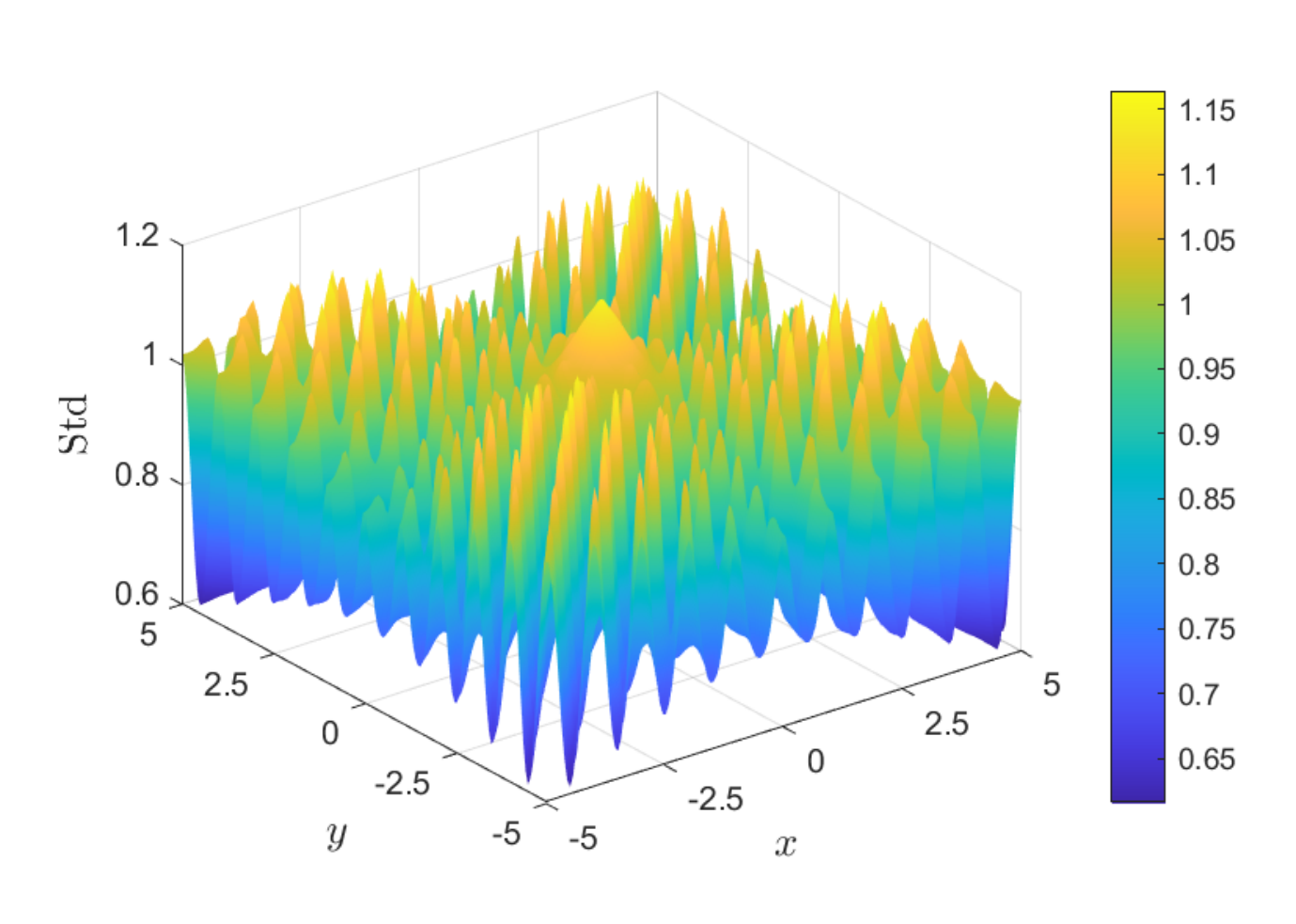}
         \caption{}
         \label{fig:Ack_3D_Std_lhs}
     \end{subfigure}
     
     \begin{subfigure}[b]{0.49\textwidth}
      \centering
        \includegraphics[width=\textwidth]{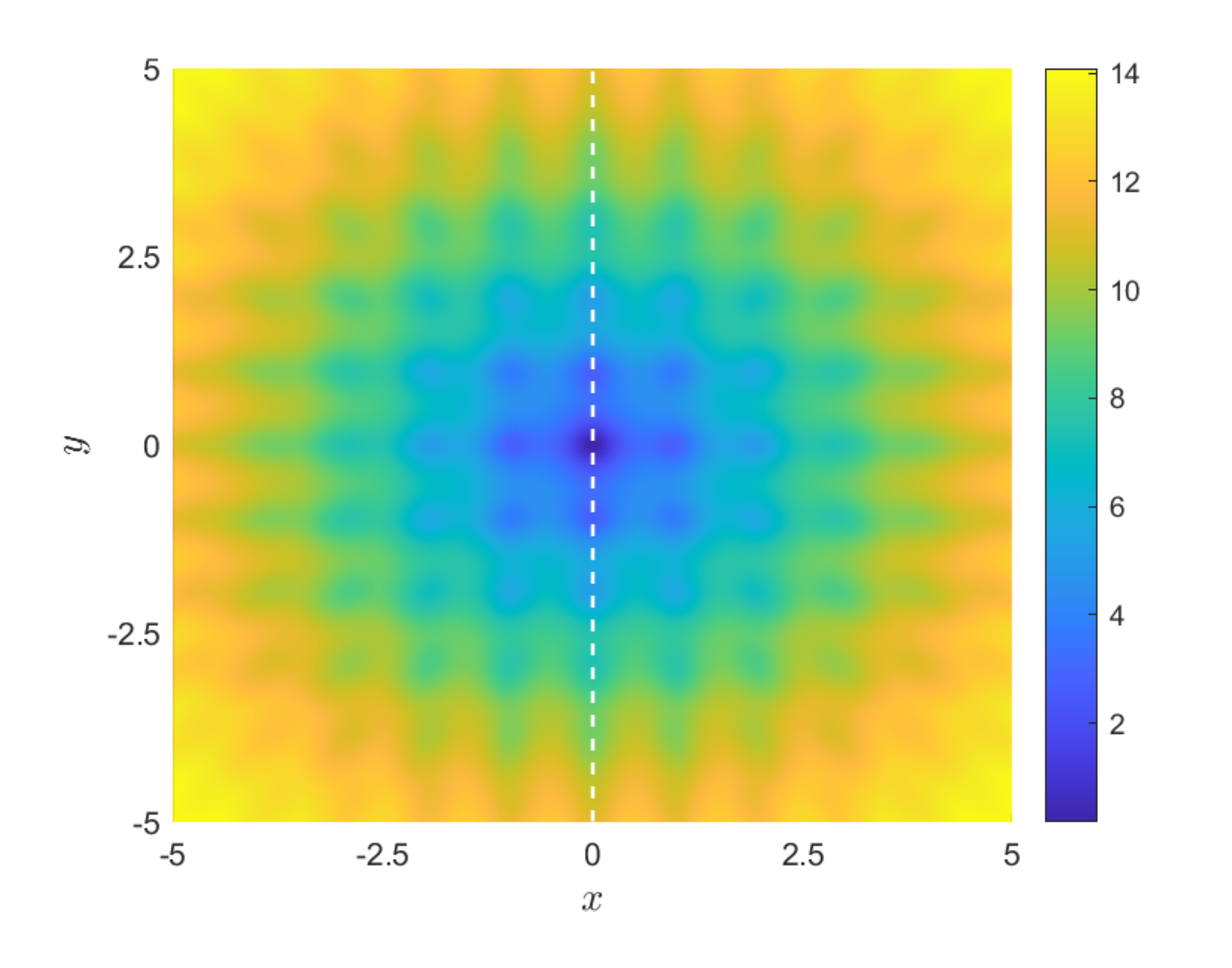}
         \caption{}
         \label{fig:Ack_Countour_Mean_lhs}
    \end{subfigure}  
  \hfill
    \begin{subfigure}[b]{0.49\textwidth}
      \centering
        \includegraphics[width=\textwidth]{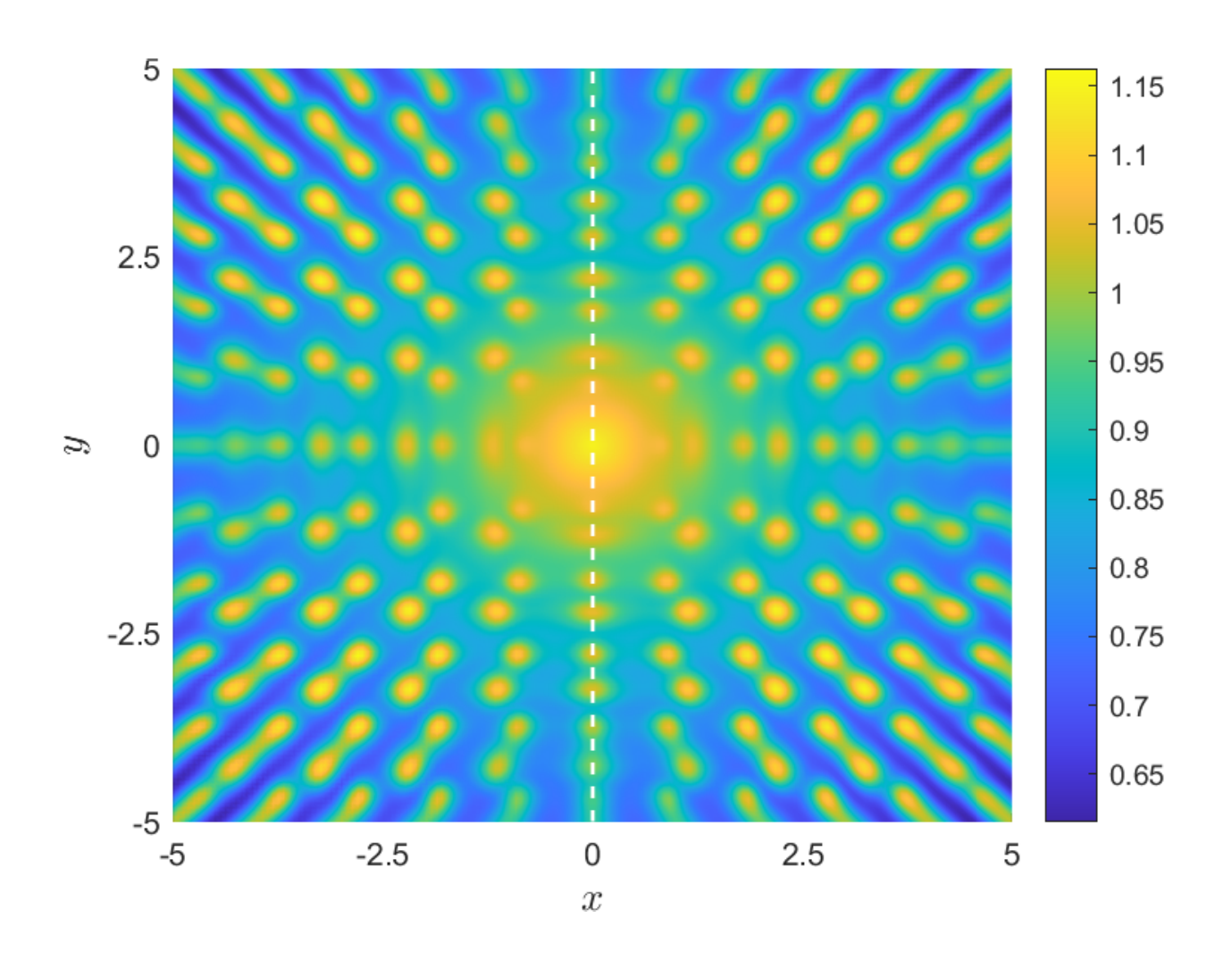}
         \caption{}
         \label{fig:Ack_Countour_Std_lhs}
     \end{subfigure}
   \caption{Mean and standard deviation fields for the Ackley function obtained from $N_{s}=50\;000$ LHS realizations. First row: 3D plots, second row: contour plots. (For interpretation of the references to color  in this figure, please refer to the web version of this article).}
   \label{fig:Ack_Mean_Std_lhs}
\end{figure}

The results of the non-intrusive reduced order models (POD-BSBEM and POD-ANN) are compared with the solution obtained using  $50\,000$ realizations on a sample generated with the LHS algorithm, considered as a reference solution. Fig.\ref{fig:Ack_Mean_Std_lhs} presents the 3D view and contour plots of the statistics in order to show the high irregular aspect of the Ackley function. In the following, the results concerning the present test case are mainly  displayed and discussed in terms of the variation of the mean and the standard deviation profiles as a function of $y$ at $x=0$ (cross-section chosen arbitrarily).\\

\begin{figure}[ht!]
  \centering
    \begin{subfigure}[b]{0.49\textwidth}
      \centering
        \includegraphics[width=\textwidth]{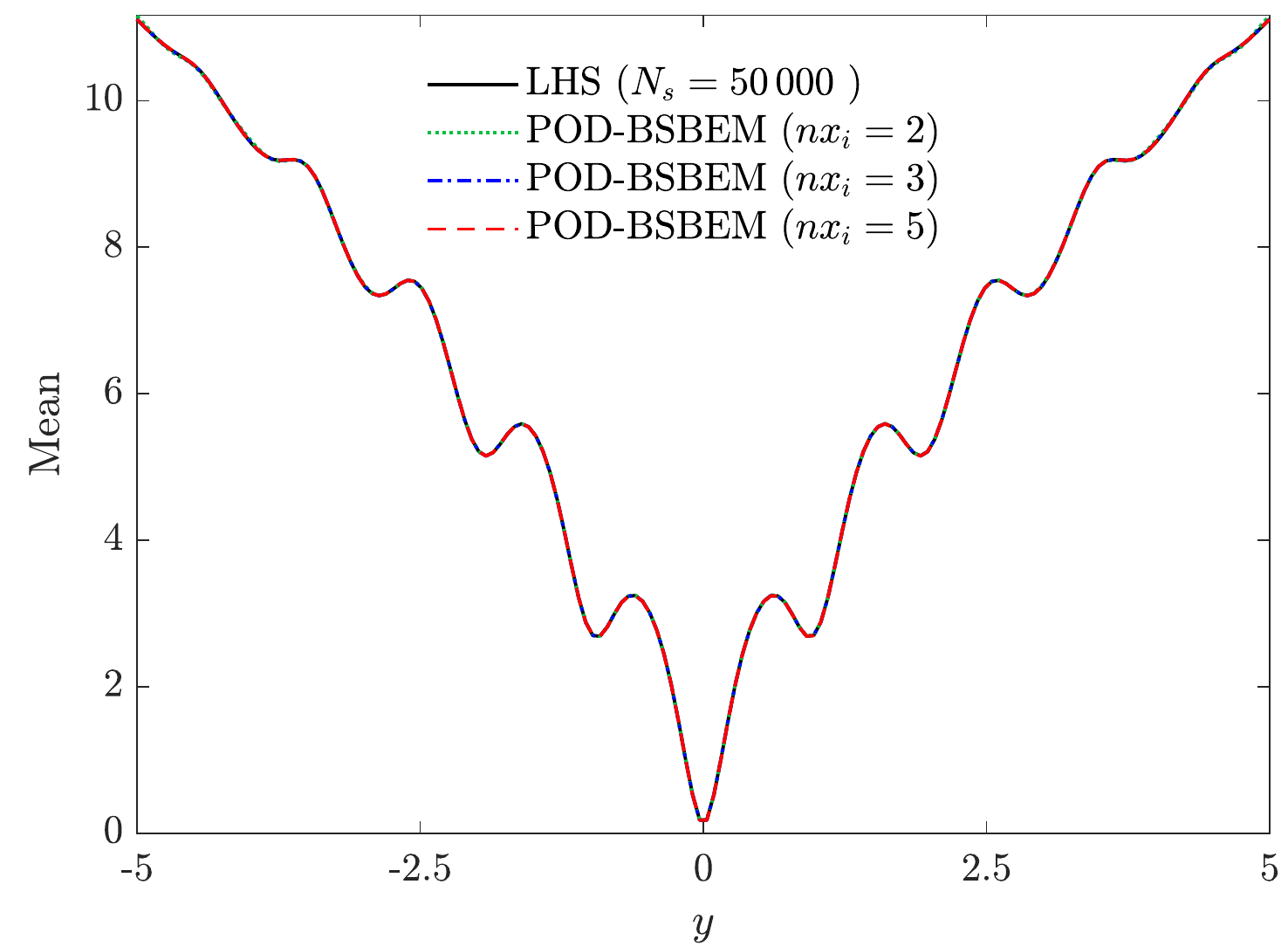}
         \caption{$\epsilon_{s}=10^{-5}\;(L=6\;modes)$}
         \label{fig:Ack_POD_BSBEM_mean_L6}
    \end{subfigure}  
  \hfill
    \begin{subfigure}[b]{0.49\textwidth}
      \centering
        \includegraphics[width=\textwidth]{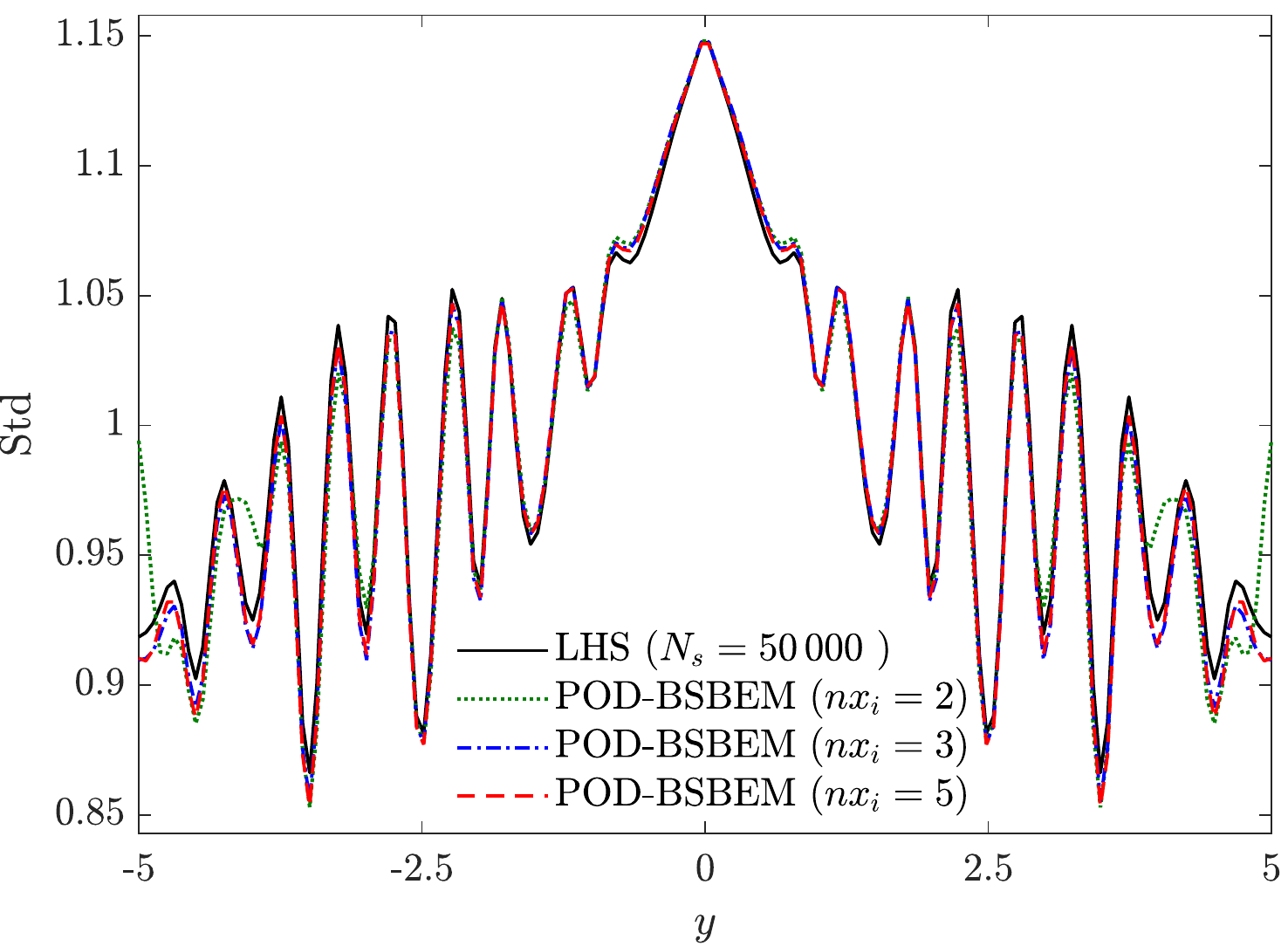}
         \caption{$\epsilon_{s}=10^{-5}\;(L=6\;modes)$}
         \label{fig:Ack_POD_BSBEM_std_L6}
     \end{subfigure}
\hfill
    \begin{subfigure}[b]{0.49\textwidth}
      \centering
        \includegraphics[width=\textwidth]{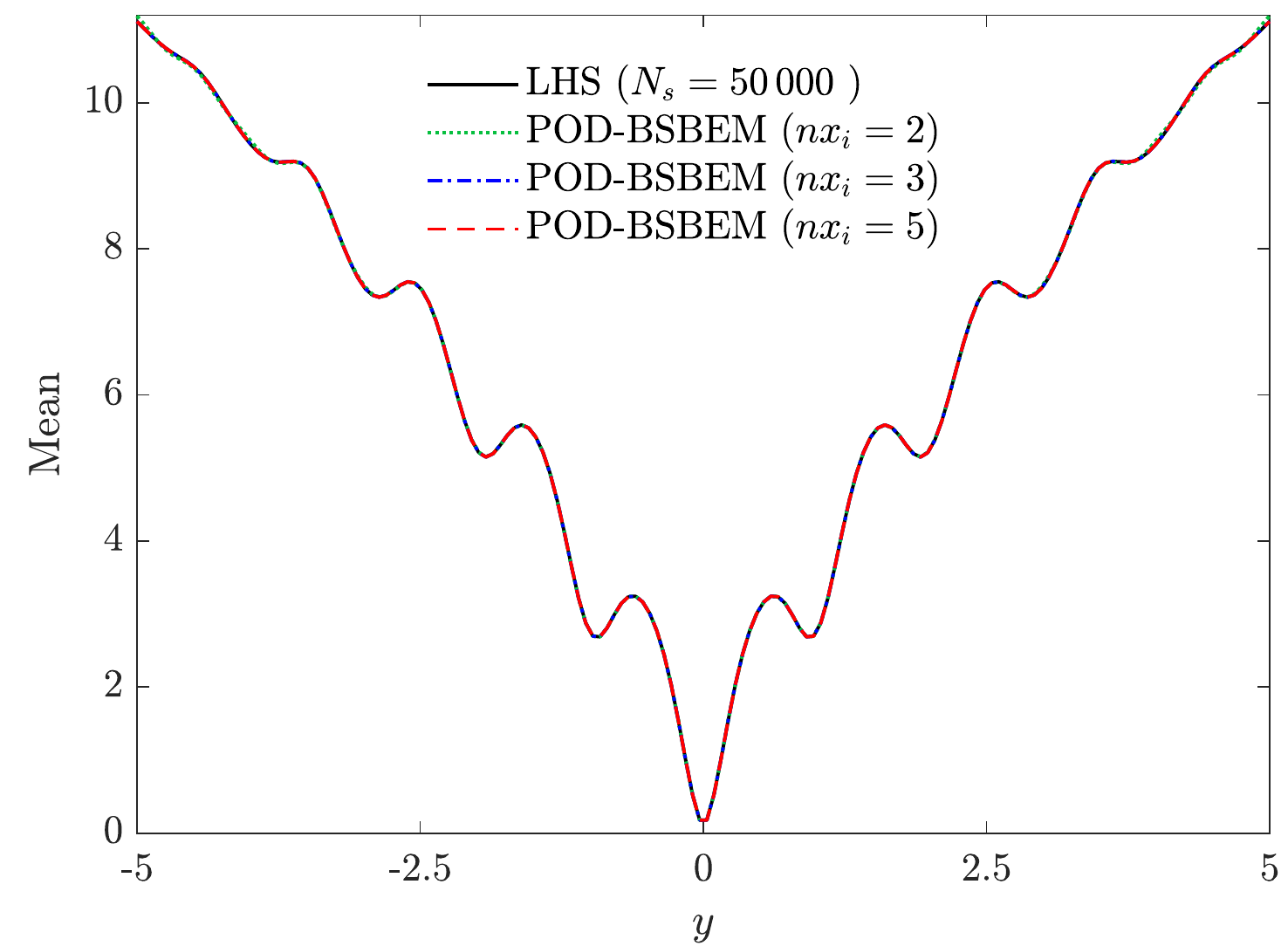}
         \caption{$\epsilon_{s}=10^{-10}\;(L=14\;modes)$}
         \label{fig:Ack_POD_BSBEM_mean_L14}
    \end{subfigure}  
  \hfill
    \begin{subfigure}[b]{0.49\textwidth}
      \centering
        \includegraphics[width=\textwidth]{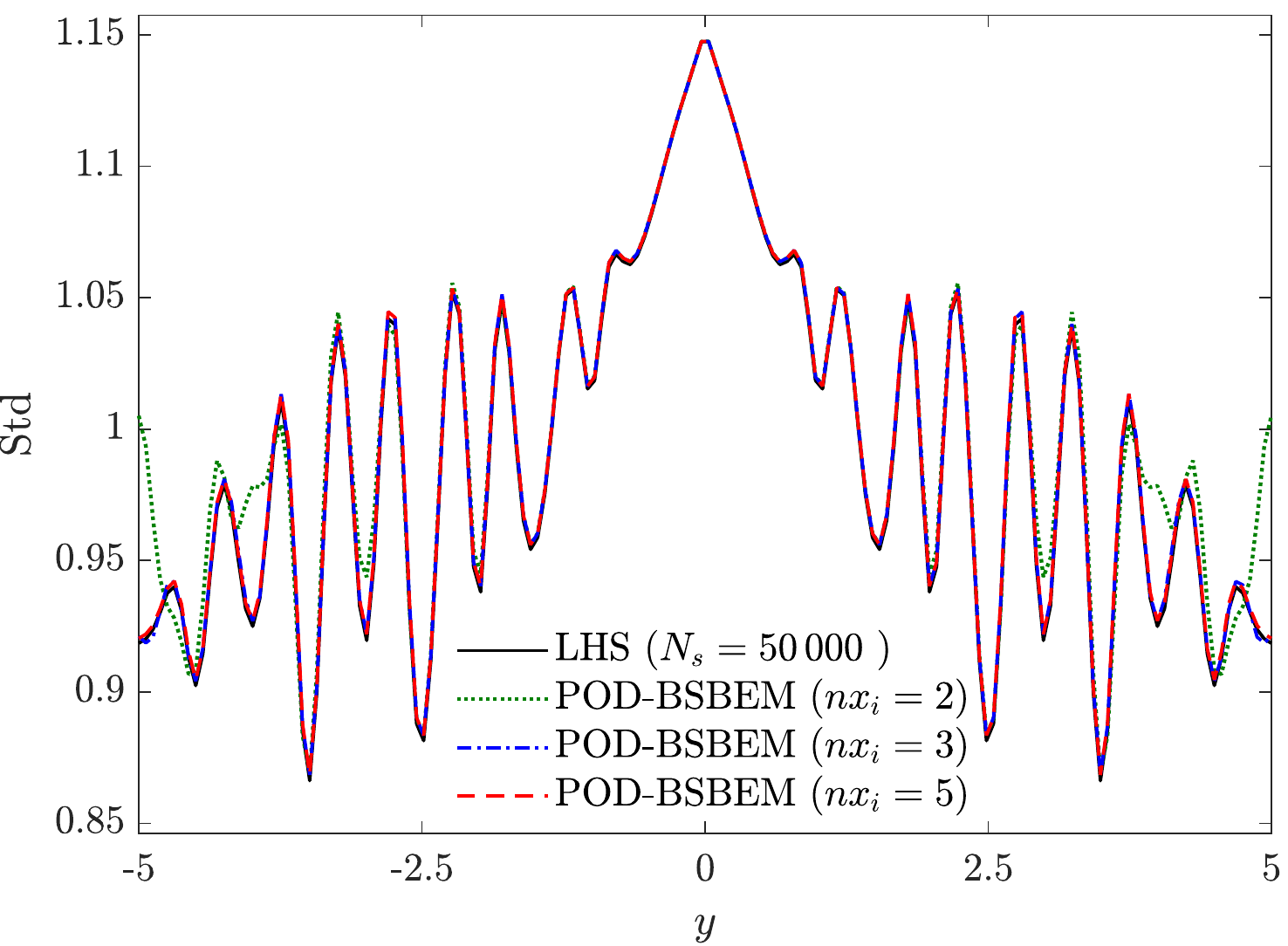}
         \caption{$\epsilon_{s}=10^{-10}\;(L=14\;modes)$}
         \label{fig:Ack_POD_BSBEM_std_L14}
     \end{subfigure}
   \caption{Statistical moments as a function of $y$ (at $x=0$) obtained using the non-intrusive POD-BSBEM with $p_{i\in\{1,2,3\}}=2$ and $nx_{i\in\{1,2,3\}}=2,\;3,\;5$: $L=6$ modes (a-b) and $L=14$ modes (c-d), compared to the LHS method (with $N_{s}=50\,000$ realizations).}
   \label{fig:Ack_POD_BSBEM}
\end{figure}

Fig.\ref{fig:Ack_POD_BSBEM} presents the variation of the mean and the standard deviation as a function of $y$ (at cross-section $x=0$) obtained by the non-intrusive reduced order model POD-BSBEM with polynomial order $p_{i}=2$ and various values of the number of B\'{e}zier elements ($nx_{i}=2,3$ and $5$, with $i\in\{1,2,3\}$). The results are displayed for two values of the energy reconstruction level, $\epsilon_{s}=10^{-5}$ and $10^{-10}$, corresponding to the number of modes $L=6$ and $14$, respectively. These plots show that the results obtained by the proposed reduced model (POD-BSBEM) for the mean statistics accurately reproduce the reference LHS solution (with $50\,000$ realizations), independently of the number of B\'{e}zier elements and the number of POD modes. The standard deviation profile reveals a clear deviation between the solution obtained by the POD-BSBEM approach and that from the LHS method for low values of the number of B\'{e}zier elements ($nx_{i}=2$) and of the number of modes $L=6$ (Fig.\ref{fig:Ack_POD_BSBEM_std_L6}). This deviation decreases as the number of B\'{e}zier elements increases, where a satisfactory matching can be noted with the LHS reference solution, more particularly for $nx_{i}=5$ and $L=14$, as shown in Fig.\ref{fig:Ack_POD_BSBEM_std_L14}.\\

\begin{figure}[ht!]
  \centering
    \begin{subfigure}[b]{0.49\textwidth}
      \centering
        \includegraphics[width=\textwidth]{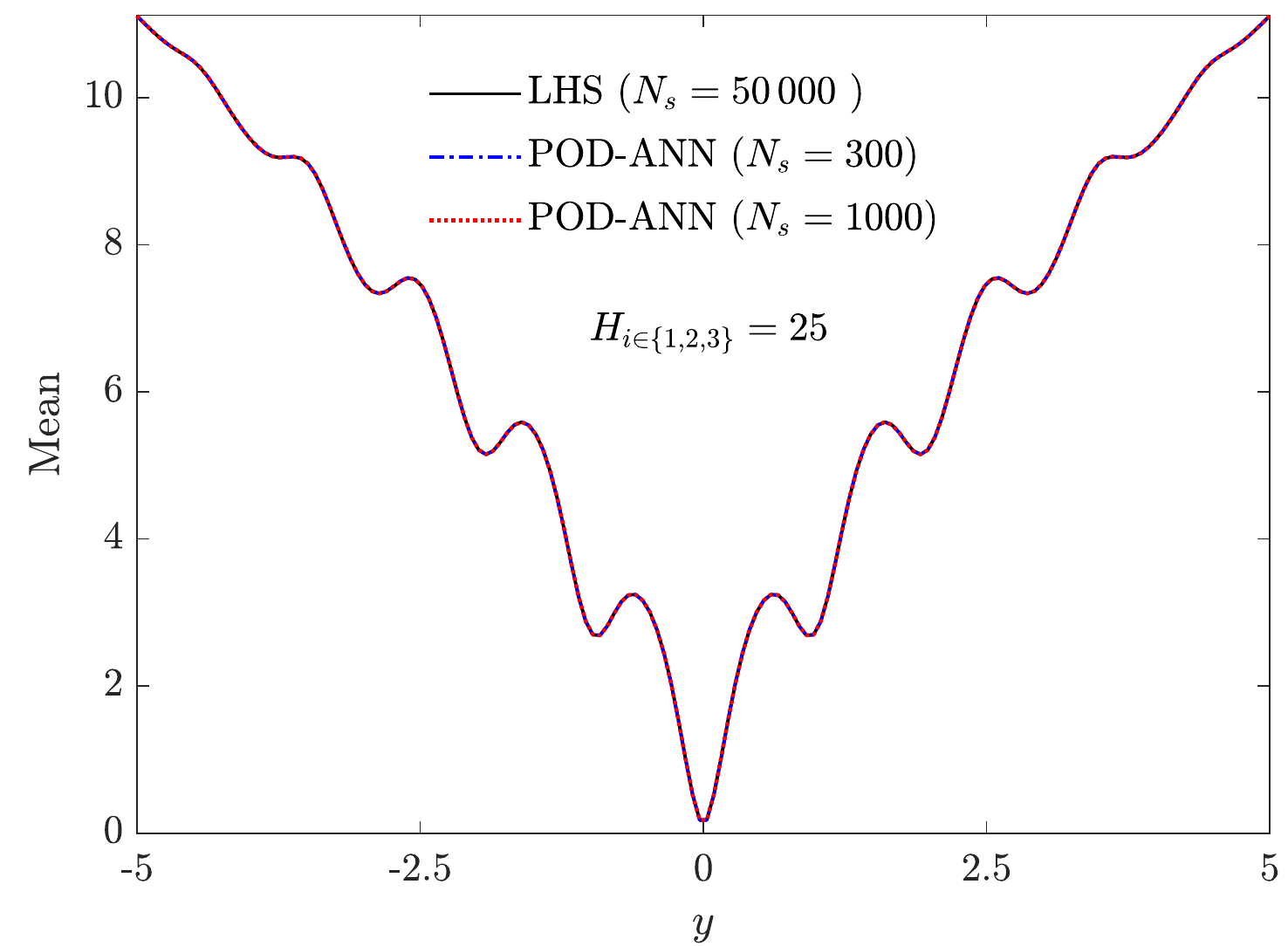}
         \caption{$\epsilon_{s}=10^{-5}\;(L=6\;modes)$}
         \label{fig:Ack_POD_ANN_mean_L6}
    \end{subfigure}  
  \hfill
    \begin{subfigure}[b]{0.49\textwidth}
      \centering
        \includegraphics[width=\textwidth]{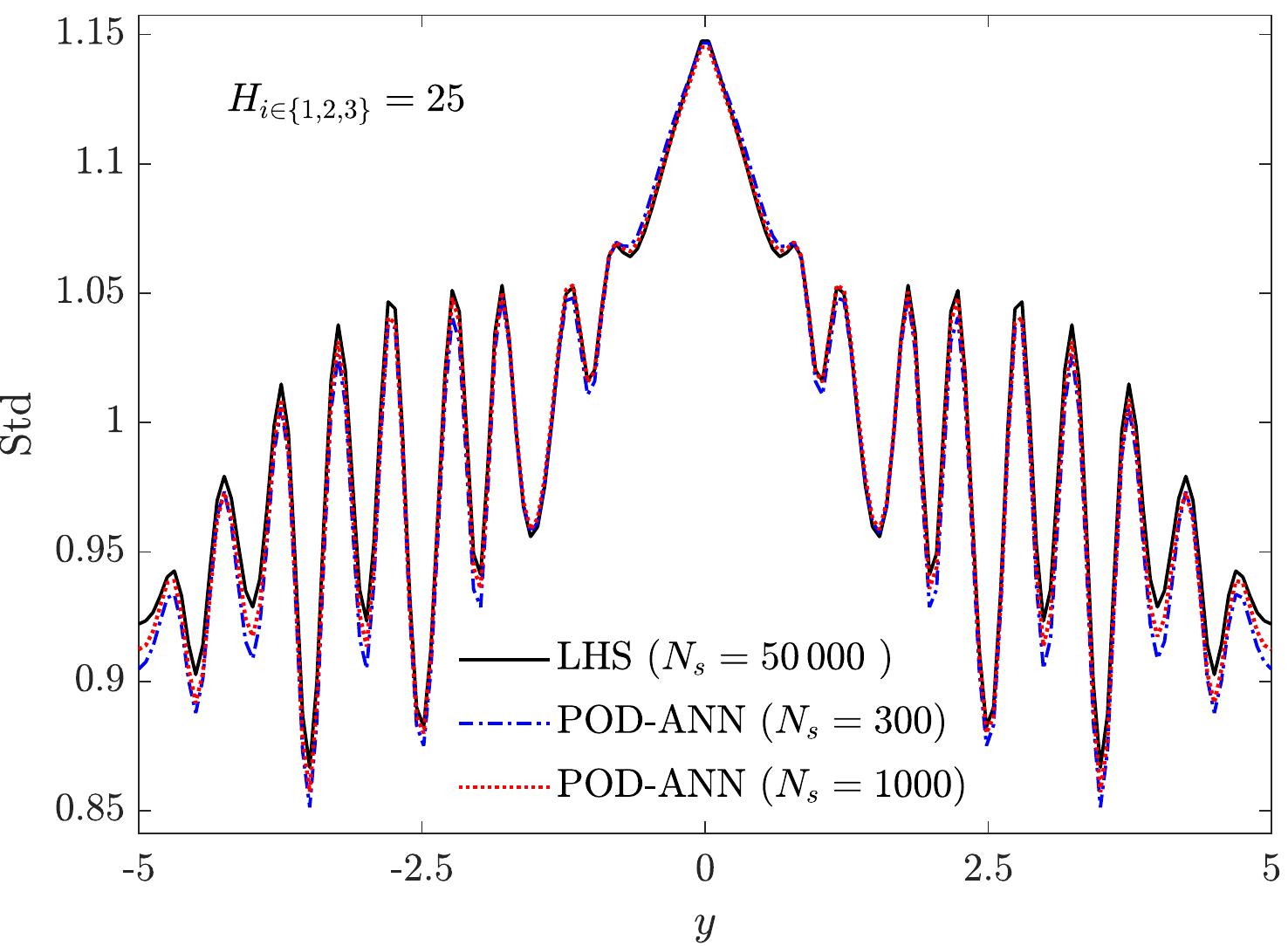}
         \caption{$\epsilon_{s}=10^{-5}\;(L=6\;modes)$}
         \label{fig:Ack_POD_ANN_std_L6}
     \end{subfigure}
\hfill
    \begin{subfigure}[b]{0.49\textwidth}
      \centering
        \includegraphics[width=\textwidth]{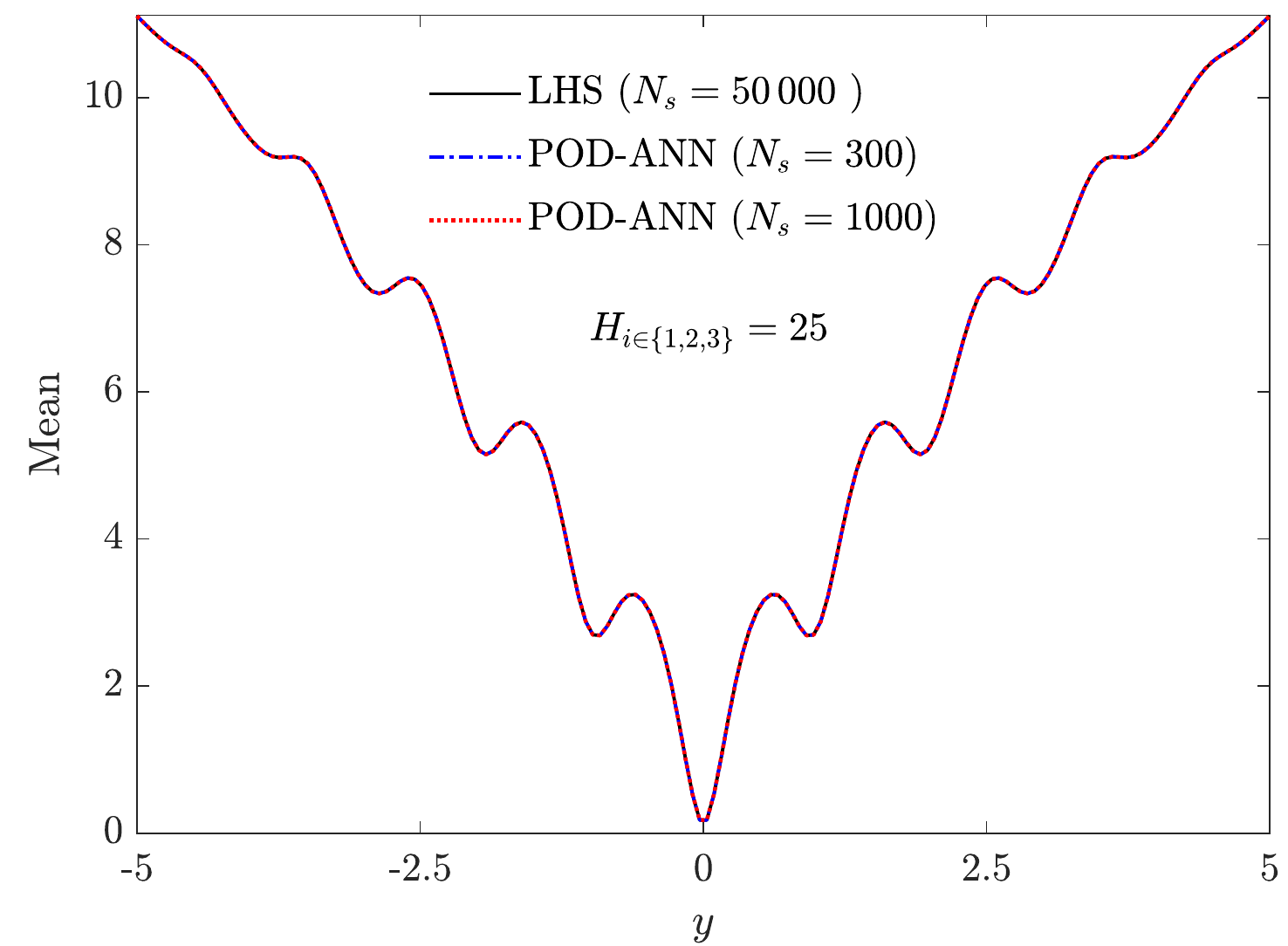}
         \caption{$\epsilon_{s}=10^{-10}\;(L=14\;modes)$}
         \label{fig:Ack_POD_ANN_mean_L14}
    \end{subfigure}  
  \hfill
    \begin{subfigure}[b]{0.49\textwidth}
      \centering
        \includegraphics[width=\textwidth]{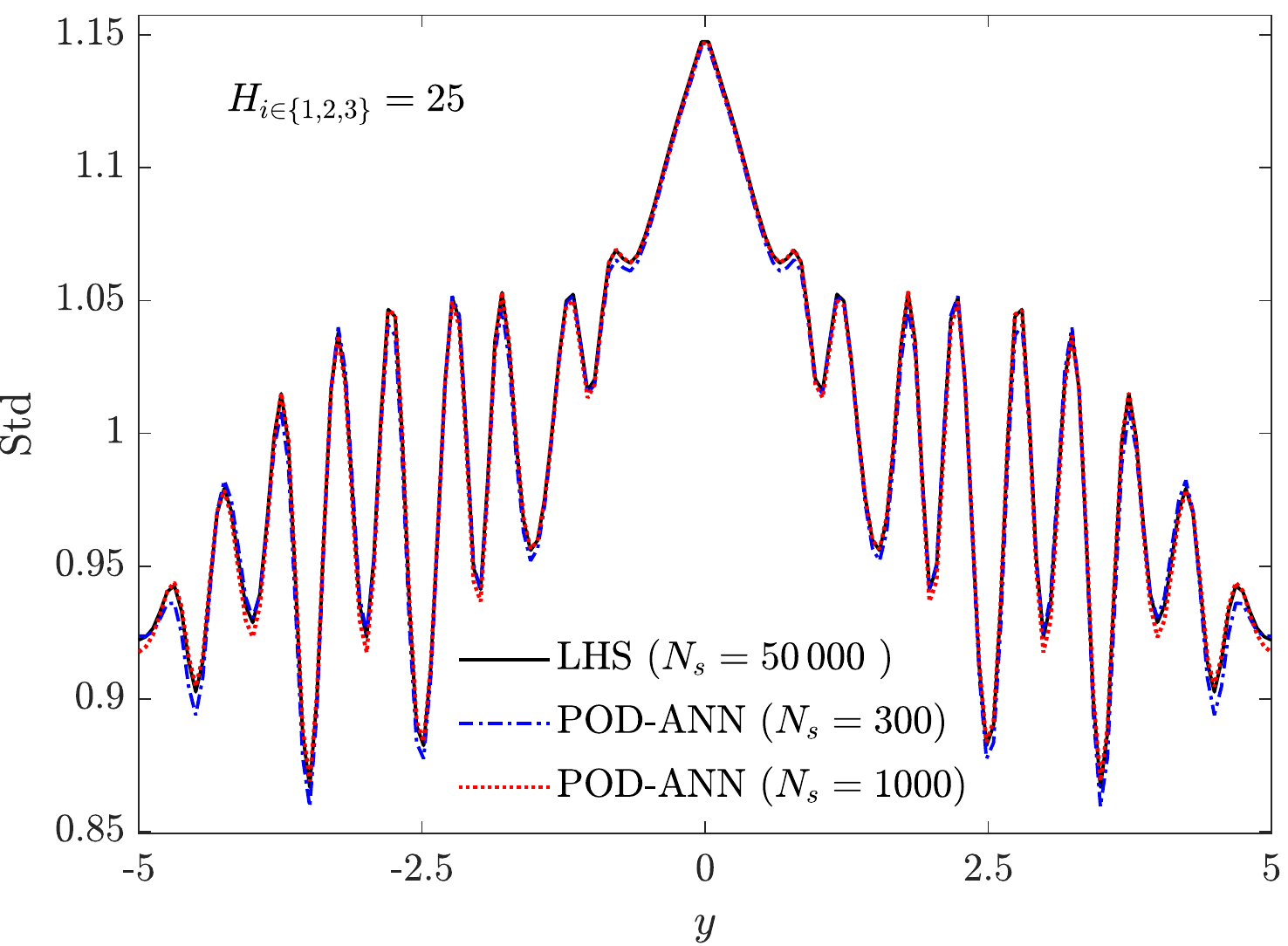}
         \caption{$\epsilon_{s}=10^{-10}\;(L=14\;modes)$}
         \label{fig:Ack_POD_ANN_std_L14}
     \end{subfigure}
   \caption{Statistical moments as a function of $y$ (at $x=0$) obtained using the non-intrusive POD-ANN with $H_{i\in\{1,2,3\}}=25$ and $N_{s}=300,\;1000$: $L=6$ modes (a-b) and $L=14$ modes (c-d), compared to the LHS method (with $N_{s}=50\,000$ realizations).}
   \label{fig:Ack_POD_ANN}
\end{figure}

The variation of the statistical moments of the Ackley function as a function of $y$ at cross-section $x=0$ obtained using the reduced model based on artificial neural networks (POD-ANN) with three hidden size layers $H_{1}=H_{2}=H_{3}=25$ is shown in Fig.\ref{fig:Ack_POD_ANN}. It should be mentioned that the neural network topology adopted in the present test case was selected after a series of numerical tests by varying the number of hidden layers and the width of each layer, in which a trade-off was made between the network performance (based on the mean squared error) and the learning time. As mentioned above, the Bayesian regularization of the training algorithm is employed in the supervised learning process to iteratively adjust the regularizing hyperparameter. Two sets of the uncertain input parameters ($N_{s}=300$ and $1000$) are randomly selected within the parameter domain using the LHS approach. The generated sampling data are all scaled to the same range using the mean normalization \citep{wang2019non}, i.e., $\widehat{\eta}_{i}=\frac{\eta_{i}-\tilde{\eta}}{\eta_{max}-\eta_{min}},\,i=1\ldots,N_{s}$, where $\tilde{\eta}$, $\eta_{min}$, and $\eta_{max}$ represent the mean, minimum and maximum of the input data, respectively. Their corresponding outputs, which represent the snapshot matrix, are obtained from $N_{s}$ runs of the high-fidelity solution (Ackley function). The results in Fig.\ref{fig:Ack_POD_ANN} show that for a relatively low number of eigenmodes ($L=6$), the standard deviation profiles predicted by POD-ANN present a visible deviation with those of the reference LHS solution for both sizes of the learning sets (Fig.\ref{fig:Ack_POD_ANN_std_L6}). When the number of modes increases to $L=14$ (obtained with an energy construction level of $\epsilon_{s}=10^{-10}$), a remarkable improvement in the predicted standard deviation profile can be observed, where the POD-ANN solution agrees well with the LHS reference solution, particularly for $N_{s}=1000$, as shown in Fig.\ref{fig:Ack_POD_ANN_std_L14}.\\

\begin{figure}[ht!]
  \centering
    \begin{subfigure}[b]{0.49\textwidth}
      \centering
        \includegraphics[width=\textwidth]{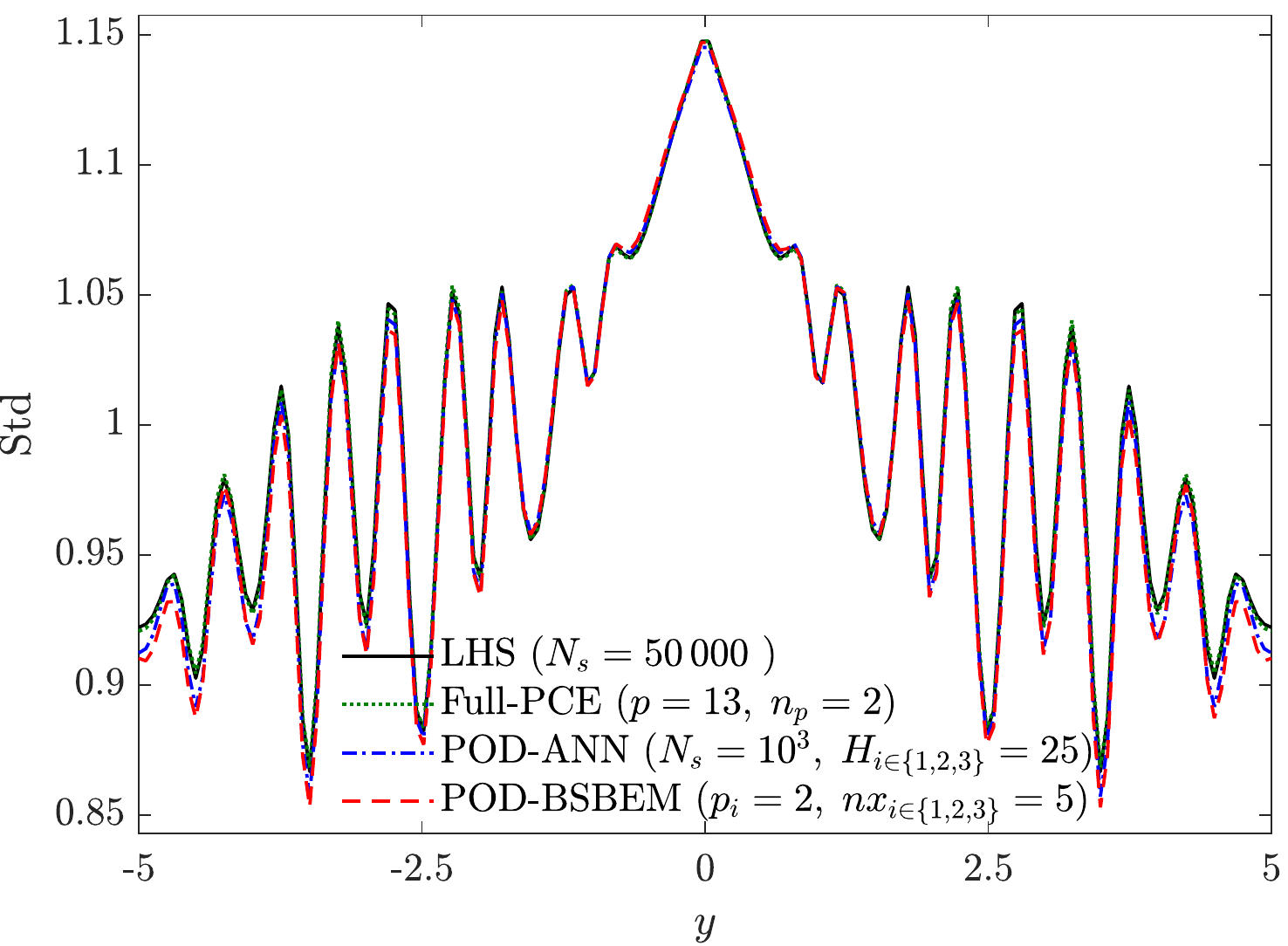}
         \caption{$\epsilon_{s}=10^{-5}\;(L=6\;modes)$}
         \label{fig:ACK_PCE_POD_ANN_BSBEM_L6}
    \end{subfigure}  
  \hfill
    \begin{subfigure}[b]{0.49\textwidth}
      \centering
        \includegraphics[width=\textwidth]{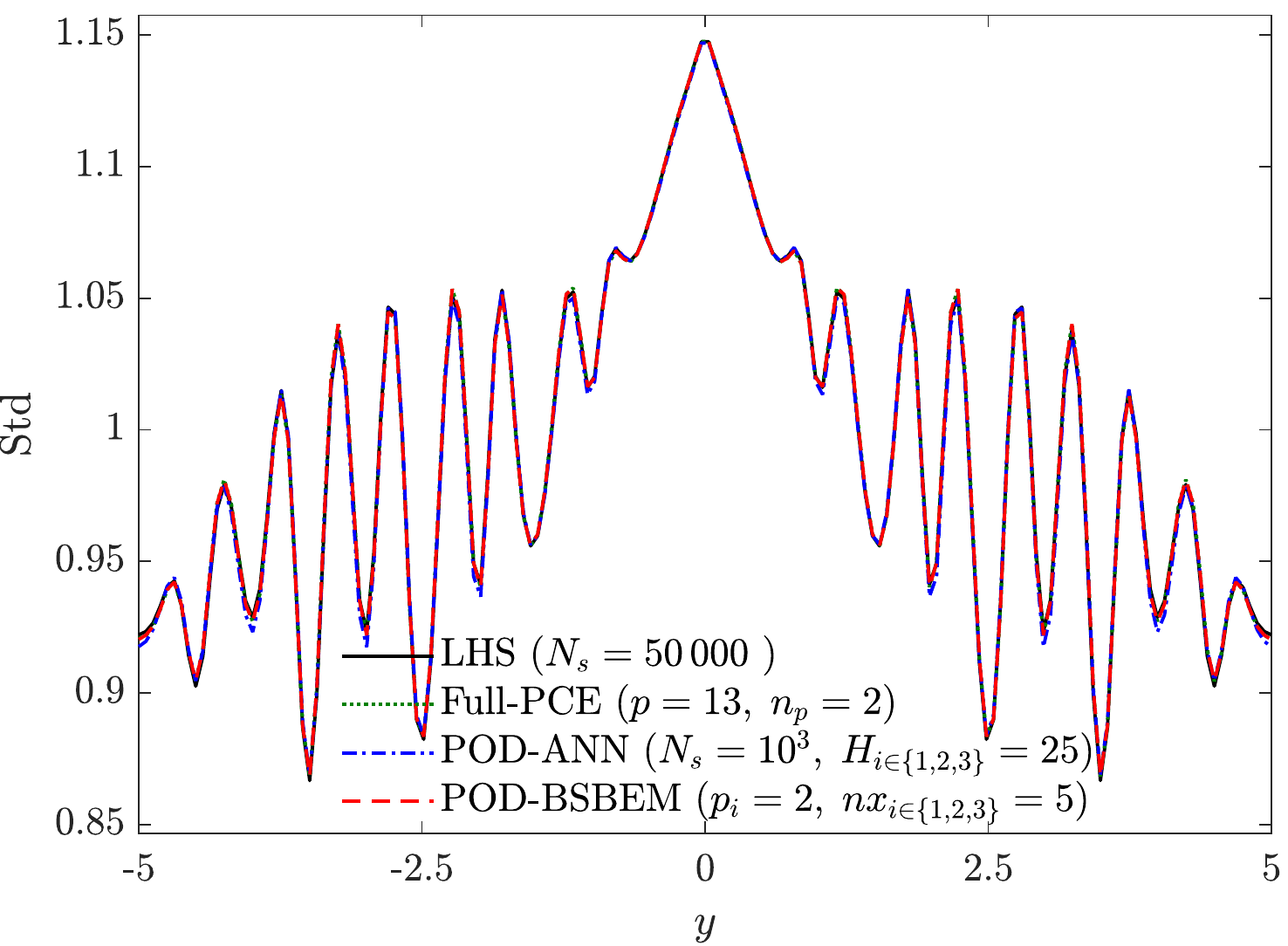}
         \caption{$\epsilon_{s}=10^{-10}\;(L=14\;modes)$}
         \label{fig:ACK_PCE_POD_ANN_BSBEM_L14}
     \end{subfigure}      
   \caption{Comparison of the standard deviation profiles obtained with the Full-PCE, POD-ANN and POD-BSBEM with that of the LHS reference solution (with $50\,000$ realizations). (a): $\epsilon_{s}=10^{-5}$ ($L=6$ modes); (b): $\epsilon_{s}=10^{-10}$ ($L=14$ modes).}
   \label{fig:ACK_PCE_POD_ANN_BSBEM}
\end{figure}

To further investigate the ability of the proposed non-intrusive POD-BSBEM to accurately approximate the output stochastic response of quantities of interest, its results are compared with those from the Full-PCE and POD-ANN approaches, and to the results of $50\,000$ realizations of the LHS (reference) method, as shown in Fig.\ref{fig:ACK_PCE_POD_ANN_BSBEM}. These results are presented solely in terms of the variation of the standard deviation profile, as the mean is well reproduced by all four methods and for any value of their hyper parameters, as shown in Figs.\ref{fig:Ack_POD_BSBEM}-\ref{fig:Ack_POD_ANN}. The results from the Full Polynomial chaos expansion  model are provided in order to show the benefits of introducing the reduced order model concept in the estimation of the stochastic output response. The Full-PCE solution, obtained with a polynomial order of $p=13$ and oversampling ratio $n_{p}=2$, as suggested in \citep{raisee2015non}, shows a good agreement with the LHS solution. It should be emphasized that the Full-PCE method is implemented by adopting the classical polynomial chaos decomposition such that the expansion's coefficients are computed by regression for each of the $N_{e}=25\,600$ nodes representing the computational mesh considered in the present test case. As illustrated in Fig.\ref{fig:ACK_PCE_POD_ANN_BSBEM_L6}, the results obtained by the reduced-order models POD-BSBEM and POD-ANN, built with $L=6$ significant modes (corresponding to a threshold $\epsilon_{s}=10^{-5}$), show a slight deviation in the std  profile compared with those from the Full-PCE and LHS methods. Fig.\ref{fig:ACK_PCE_POD_ANN_BSBEM_L14} shows a significant improvement in the prediction of the std  profile obtained by the reduced-order models as the number of dominant modes increases ($L=14$), particularly with the POD-BSBEM, where a faithful reproduction of the statistical moment can be observed, compared to that of the LHS reference approach, and with much less computational effort.\\

For a more profound insight, the accuracy of POD-BSBEM in the approximation of output statistics of quantities of interest is investigated by evaluating the relative error in the $L^{2}$-norm, expressed as follow:
\begin{equation}\label{<L2_Err>}
Err_{L^{2},\Phi}^{Surr}=\sqrt{\frac{\sum_{i=1}^{N_{e}}(\Phi_{i,Surr}-\Phi_{i,LHS})^2}{\sum_{i=1}^{N_{e}}(\Phi_{i,LHS})^2}}
\end{equation}
with $\Phi$ denotes either the mean or the standard deviation, and $Surr$ stands for Full-PCE, POD-ANN and POD-BSBEM. Errors are computed by evaluating the statistical moments obtained by both full and reduced order surrogate models (Full-PCE, POD-ANN and POD-BSBEM) and comparing them with those of the LHS reference solution over the $N_{e}$ total number of nodes representing the computational domain. It is obvious from the results reported in Table \ref{tab:TabFull-PCE_ANN_BSBEM_epsi} that the decrease in the tolerance $\epsilon_{s}$ (which results in an increase in the POD dimension) leads to a decrease in the $L^{2}$-error for both reduced-order models (POD-BSBEM and POD-ANN). This indicates that the proposed POD-BSBEM approximates the mean and standard deviation with the same accuracy as that of the Full-PCE model. However, it is worth mentioning that expansions' construction costs vary significantly between both models, as only $14$ decompositions are required for POD-BSBEM (which corresponds to $\epsilon_{s}=10^{-10}$) as a function of the random input parameters, where a vector of global coefficients has to be computed for each mode, while the Full-PCE must build a decomposition for each node of the $N_{e}=25\,600$ mesh points representing the entire computational domain. Hence, the proposed reduced-order models offers a very real gain in terms of reducing the computational efforts required for  the prediction of the statistical moments without loss of accuracy compared to the Full-PCE, particularly for cases that require a high number of mesh nodes, such as for dam-break flow applications, where the computational domain may contain millions of nodes.\\

\begin{table}[h!]
\caption{Effect of $\epsilon_{s}$ on the relative error in the $L^{2}$-norm for the mean and standard deviation obtained from Full-PCE ($p=13$, $n_{p}=2$), POD-ANN ($H_{i\in\{1,2,3\}}=25$, $N_{s}=10^{3}$) and POD-BSBEM ($p_{i\in\{1,2,3\}}=2$, $nx_{i\in\{1,2,3\}}=5$). Errors are computed with respect to the LHS reference solution (with $N_{s}=50\,000$ realizations).}
\centering
\begin{tabular*}{0.85\textwidth}{@{\extracolsep{\fill}}l l  llll }
\hline        
& & \multicolumn{2}{c}{$Err_{L^{2},\,Mean}^{Surr}$}& \multicolumn{2}{c}{$Err_{L^{2},\,Std}^{Surr}$}\\
\cline{3-4}
\cline{5-6}
\raisebox{2ex}{$ \epsilon_{s} $} & \raisebox{2ex}{$ L $} &POD-ANN& POD-BSBEM& POD-ANN& POD-BSBEM\\
\hline        
 $10^{-3}$& $3$& $9.1591$E-05& $1.9050 $E-04& $0.0673$& $0.072232$\\
$10^{-5}$& $6$& $6.7062$E-05& $7.8300$E-05& $0.009251$& $0.007650$\\
$10^{-8}$& $11$& $5.3319$E-05& $3.3817$E-05& $0.004231$& $0.002541$\\
$10^{-10}$& $14$& $1.1563$E-05& $2.6996$E-05& $0.003692$& $0.002540$\\

\hline        
\multicolumn{2}{l}{Full-PCE}&\multicolumn{2}{c}{$1.7572$E-05}&\multicolumn{2}{c}{$0.002552$}\\
\hline 
\end{tabular*}
\label{tab:TabFull-PCE_ANN_BSBEM_epsi}
\end{table}

To gain a better understanding of the proposed non-intrusive reduced-order model performance, the effect of the number of B\'{e}zier elements ($nx_{i}$) on the $L^{2}$-error norm of the statistical moments for the tolerance $\epsilon_{s}=10^{-10}$ is investigated. As with the former case, the results, reported in Table \ref{tab:TabFull-PCE_ANN_BSBEM_nxi}, are compared with those from the POD-ANN and the Full-PCE in terms of computational efficiency. The number of B\'{e}zier elements represents a key controlling parameter in the POD-BSBEM concept; it defines the number of elements within which local basis functions with compact support are used to approximate the stochastic output response. It also contributes to defining the total number of collocation points (number of numerical solver calls) given by: $N_{s}=N_{elt}.nb$ where $N_{elt}=\prod_{i=1}^{m}{nx_{i}}$ and $nb=\prod_{i=1}^{m}{(p_{i}+1)}$, where $m$ denotes the number of random input parameters and $p_{i}$ is the polynomial order of the basis functions. It is worth noting that comparing the computational costs of both reduced order-models (POD-ANN and POD-BSBEM) is especially interesting given that their accuracies are almost the same. Thus, for each value of $nx_{i}$, the aforementioned resulting number of collocation points $N_{s}$ is considered as a size of the learning set to evaluate the $L^{2}$-error in the mean and standard deviation of the POD-ANN. The reported results show that the POD-BSBEM approaches the LHS reference solution with the same accuracy as the Full-PCE as the number of B\'{e}zier elements increases, and with much lower decomposition costs. It is also worth noting that predictions from the POD-BSBEM are slightly more accurate than those from the reduced-order model-based POD-ANN, particularly when $N_{s}$ takes relatively high values.\\          

\begin{table}[h!]
\caption{Effect of the number of B\'{e}zier elements $nx_{i}$ on the relative error in the $L^{2}$-norm for the mean and standard deviation obtained from the Full-PCE ($p=13$, $n_{p}=2$), POD-ANN ($H_{i\in\{1,2,3\}}=25$) and POD-BSBEM ($p_{i\in\{1,2,3\}}=2$, $nx_{i\in\{1,2,3\}}=5$) with $\epsilon_{s}=10^{-10}$ ($L=14$). Errors are computed with respect to the LHS reference solution (with $N_{s}=50\,000$ realizations).}
\centering
\begin{tabular*}{0.85\textwidth}{@{\extracolsep{\fill}}l l  llll }
\hline        
& & \multicolumn{2}{c}{$Err_{L^{2},\,Mean}^{Surr}$}& \multicolumn{2}{c}{$Err_{L^{2},\,Std}^{Surr}$}\\
\cline{3-4}
\cline{5-6}
\raisebox{2ex}{$nx_{i\in\{1,2,3\}}$} & \raisebox{2ex}{$ N_{s} $} &POD-ANN& POD-BSBEM& POD-ANN& POD-BSBEM\\
\hline        
 $2$& $216$& $2.9978$E-05& $0.001143$& $0.005805$& $0.007084$\\
$3$& $729$& $4.9126$E-05& $8.1399$E-05& $0.002548$& $0.003234$\\
$4$& $1\,728$& $1.8819$E-05& $4.9126$E-05& $0.006231$& $0.002578$\\
$5$& $3\,375$& $3.3262$E-05& $2.6996$E-05& $0.004934$& $0.002540$\\

\hline        
\multicolumn{2}{l}{Full-PCE}&\multicolumn{2}{c}{$1.7572$E-05}&\multicolumn{2}{c}{$0.002552$}\\
\hline 
\end{tabular*}
\label{tab:TabFull-PCE_ANN_BSBEM_nxi}
\end{table}

The computational costs of the presented approaches are summarized in Table \ref{tab:Tab_comp_cpu_Ack}. All the programs involved in this test case were implemented serially and the computations were performed using a personal computer: Intel (R) Xeon (R) CPU E3-125 v6 @3.30 GHz with  16 GB of memory. It is worth mentioning that the offline phase involves the collection of snapshots from the high-fidelity solutions and building the POD reduced basis, learning, and computing coefficients by regression. The proposed reduced-order model (POD-BSBEM) has a lower computational cost in both its offline and online phases in comparison with the POD-ANN technique. It also presents a much better performance than the Full-order model polynomial chaos expansion (Full-PCE).

\begin{table}[h!]
\caption{Comparison of the computational cost required by the LHS ($N_{s}=50\,000$ realizations), Full-PCE ($p=13$, $n_{p}=2$), POD-ANN ($H_{i\in\{1,2,3\}}=25$) and POD-BSBEM ($p_{i\in\{1,2,3\}}=2$, $nx_{i\in\{1,2,3\}}=5$) with $\epsilon_{s}=10^{-10}$ ($L=14$). The time cost unit is a second ($s$).}
\centering
\begin{tabular*}{0.85\textwidth}{@{\extracolsep{\fill}}l r  rrr }
\hline        
           &POD-ANN& POD-BSBEM& Full-PCE& LHS\\
\hline        
 Offline& $3\,988.51$& $141.02$& -& -\\
Online& $264.50$& $20.45$& $478.94$& $330.15$\\
\hline 
\end{tabular*}
\label{tab:Tab_comp_cpu_Ack}
\end{table}

\subsection{Burgers' equation test case} \label{Burg_sol}
The second test case concerns the well-known viscous Burgers' equation, which describes the standard nonlinear time-dependent advection-diffusion problems in one spatial dimension. Its dimensionless form with a parameterized diffusion coefficient is given by \citep{burgers1948mathematical,guo2019data}:
\begin{equation}\label{<Burg_equ>}
\frac{\partial{u}}{\partial{t}}+u\frac{\partial{u}}{\partial{x}}=\frac{1}{Re}\frac{\partial^{2}{u}}{\partial{x^{2}}}\;,\qquad x\in\left[ 0,\;1\right],\qquad t\in\left[ 0,\;1\right]
\end{equation}
with the following initial and Dirichlet boundary conditions: $u(x,0)=\frac{x}{1+exp\left( \frac{Re}{16}(4x^{2}-1)\right)}$ and $u(0,t)=u(1,t)=0$, respectively \citep{ahmed2020long}. $Re$ denotes the Reynolds number, considered as an uncertain input parameter describing the variability surrounding the viscous effects. The aforementioned partial differential equation admits an exact analytical solution of the unknown field variable $u(x,t)$ expressed as \citep{SAN2019271,maleewong2011line}: 
\begin{equation}\label{<Burg_analy>}
u(x,t)=\frac{\frac{x}{t+1}}{1+\sqrt{\frac{t+1}{t_{0}}}exp(Re\frac{x^{2}}{4t+4})}
\end{equation}
with $t_{0}=exp(\frac{Re}{8})$. The 1D computational space domain $\Omega=\left[0,1\right]$ is discretized with a mesh of $N_{e}=1000$ nodes in which the above analytical solution is evaluated. The discretization over the time domain $\mathcal{T}=\left[0,1\right]$ is performed with $N_{t}=50$ time instances, leading to a uniform time step $\Delta t=0.02$.\\

The test case described above represents a challenging time-dependent problem for uncertainty propagation analysis. The analytical solution of the viscous Burgers' equation, given by Eq.\eqref{<Burg_analy>}, is adopted as a deterministic model through which uncertainties, traduced in terms of the variability surrounding the diffusion coefficient, are propagated into the output quantity of interest. Thus, the Reynolds number is considered as an uncertain input parameter described by uniform distribution with two nominal values $\mu_{Re}=200$ and $\mu_{Re}=800$ and a coefficient of variation of $cv=25\,\%$, which represent the standard deviation values of $50$ and $200$, respectively. The values of the input random Reynolds number are selected within its plausible variability ranges $Re_{\mu=200,\,\sigma=50}\in\mathcal{U}\left[114,\;287\right]$ and $Re_{\mu=800,\,\sigma=200}\in\mathcal{U}\left[454,\;1\,146\right]$, which correspond to a smooth and a very steep solution, respectively, of the Burgers' equation. This test case is of particular interest because it allows an efficient analysis of the efficacy of the proposed stochastic non-intrusive reduced order model  compared to the other reduced and full-order models for both smooth and steep output responses with even discontinuous shock. The most relevant results are presented below and discussed in terms of statistical moments, $L^{2}$-error, and probability density function profiles.\\

\begin{figure}[ht!]
  \centering
    \begin{subfigure}[b]{0.49\textwidth}
      \centering
        \includegraphics[width=\textwidth]{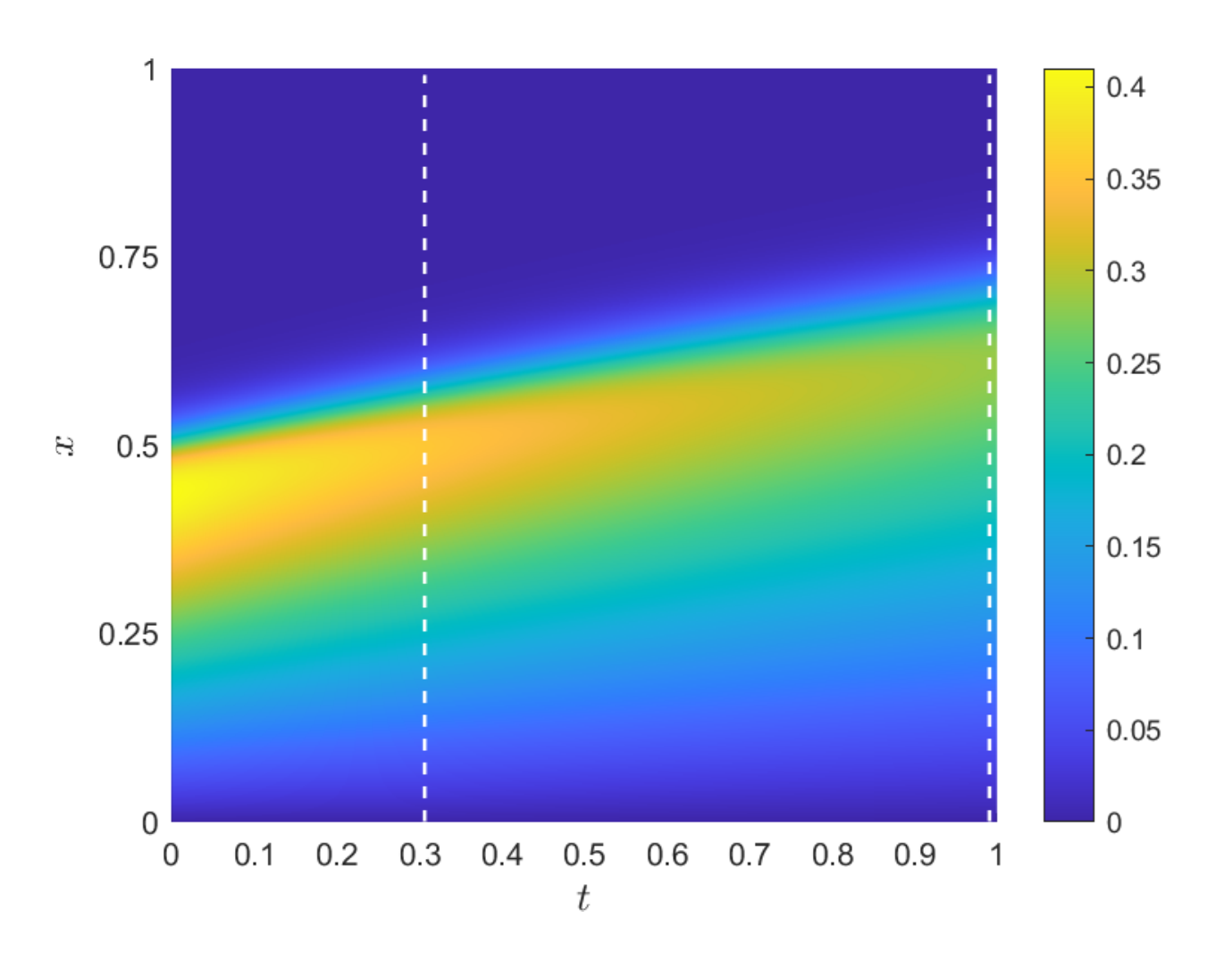}
         \caption{Mean ($\mu_{Re}=200,\,cv=25\,\%$)}
         \label{fig:Burg_mean_2d_Re_200}
    \end{subfigure}  
  \hfill
    \begin{subfigure}[b]{0.49\textwidth}
      \centering
        \includegraphics[width=\textwidth]{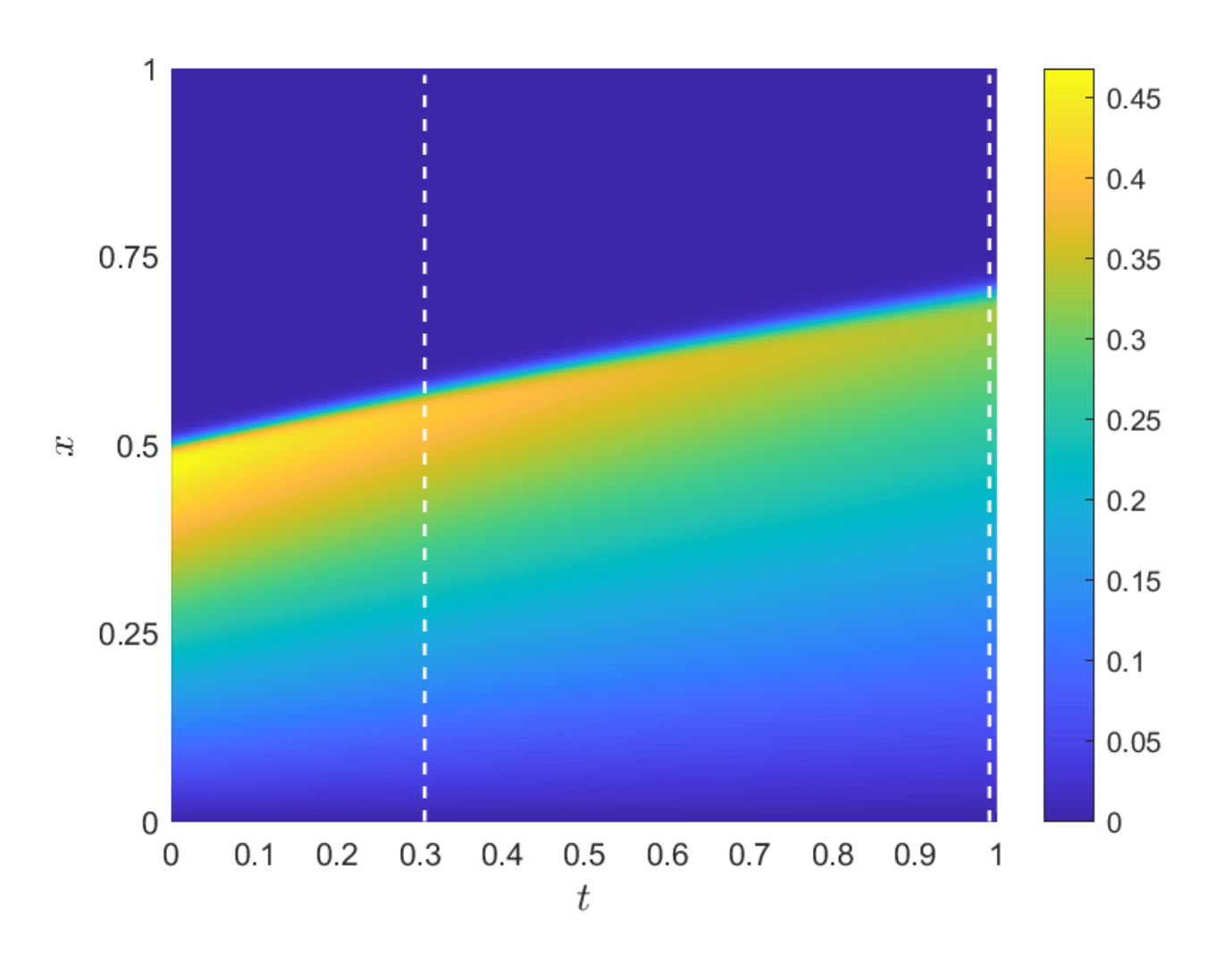}
         \caption{Mean ($\mu_{Re}=800,\,cv=25\,\%$)}
         \label{fig:Burg_mean_2d_Re_800}
     \end{subfigure} 
     \hfill
    \begin{subfigure}[b]{0.49\textwidth}
      \centering
        \includegraphics[width=\textwidth]{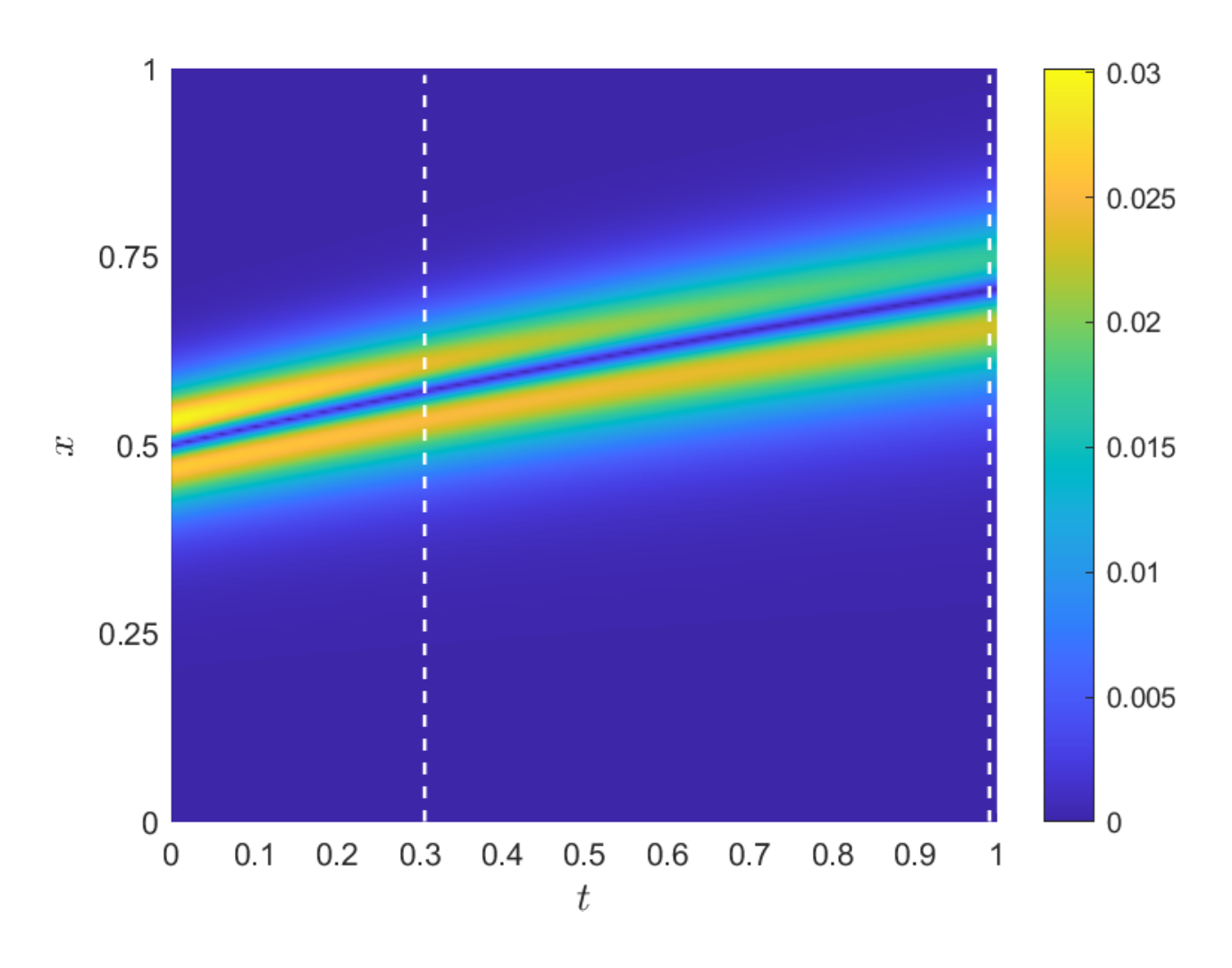}
         \caption{Std ($\mu_{Re}=200,\,cv=25\,\%$)}
         \label{fig:Burg_std_2d_Re_200}
     \end{subfigure} 
     \hfill
    \begin{subfigure}[b]{0.49\textwidth}
      \centering
        \includegraphics[width=\textwidth]{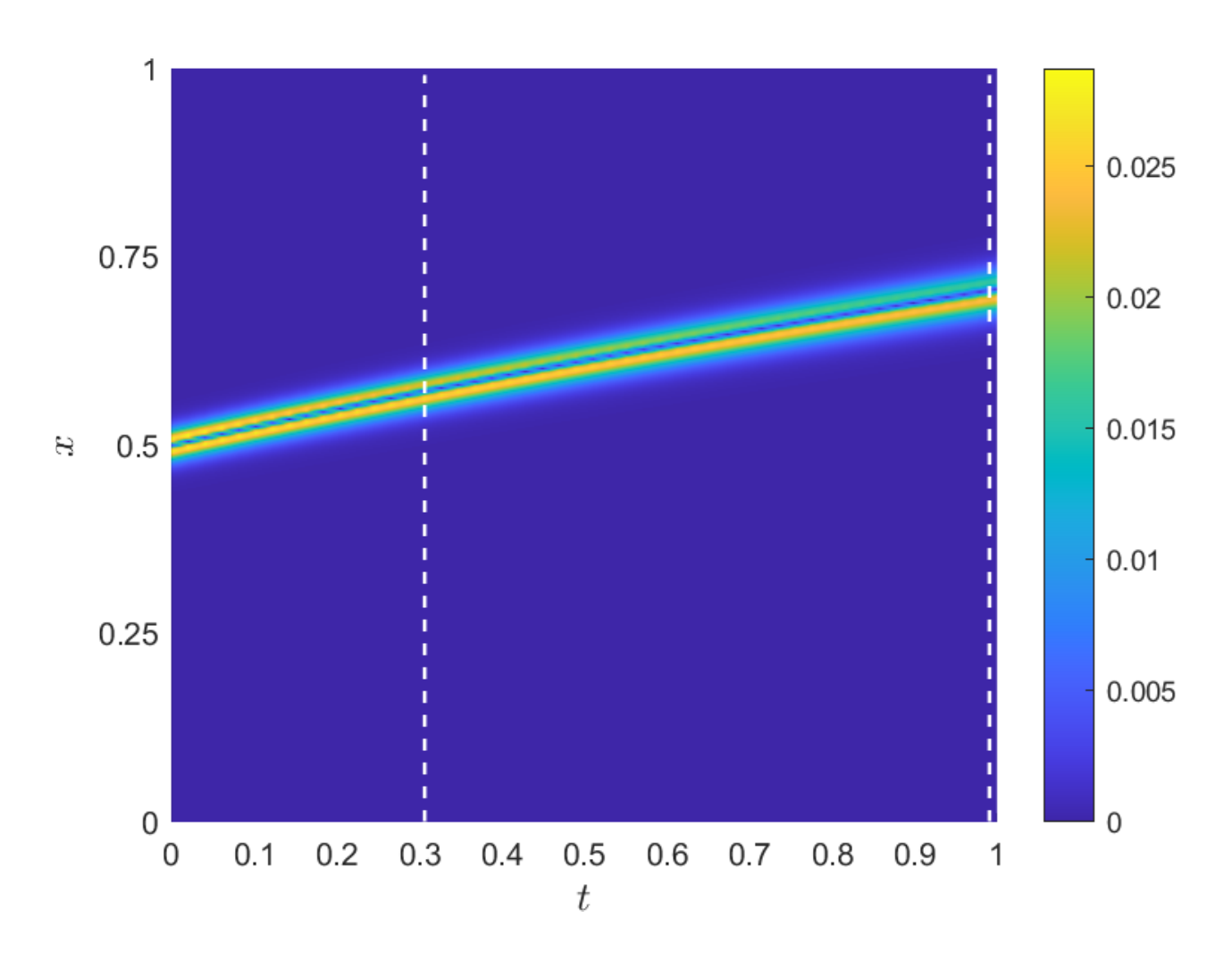}
         \caption{Std ($\mu_{Re}=800,\,cv=25\,\%$)}
         \label{fig:Burg_std_2d_Re_800}
     \end{subfigure}
   
   \caption{Contour plots of the mean (top row) and standard deviations (bottom row) of the Burgers' analytical solution obtained from $10^{6}$ Monte Carlo realizations. (a and c): $Re=200$, (b and d): $Re=800$. The two vertical white lines indicate the times for which the statistics of the Full-PCE, POD-ANN and POD-BSBEM are compared to those of the MC reference solutions. (For interpretation of the references to color in this figure, the reader is referred to the web version of this article).}
   \label{fig:Burg_2d_MC_Re_200_800}
\end{figure}

Fig.\ref{fig:Burg_2d_MC_Re_200_800} shows the contour plots of the mean (top row) and standard deviation (bottom row) of the Burgers' solution in the space-time plane, obtained with the Monte Carlo method for two variability ranges of the uncertain input parameter  corresponding to $\mu_{Re}=200$ (first column) and $\mu_{Re}=800$ (second column) with the same coefficient of variation of $cv=25\,\%$. These contour plots are presented for a better overall visualization of the statistical moments corresponding to smooth and step variation of the $u$-velocity, as can be observed in Figs.\ref{fig:Burg_mean_2d_Re_200} and \ref{fig:Burg_mean_2d_Re_800}, respectively. The two vertical dashed lines denote the $t$ locations for which the variation of the statistical moments' profiles obtained from the POD-BSBEM, POD-ANN, and Full-PCE models as a function of the $x$-space coordinate are compared to those from the MC method with $N_{s}=10^{6}$ realizations, which is considered as a high-fidelity reference solution. It is worth mentioning that the feed-forward neural network architecture retained in the POD-ANN approach to deal with the Burgers' time-dependent test case is constituted by three hidden layers, each of which has $50$ neurons ($H_{1,2,3}=50$). The constructed networks were trained using the scaled conjugate gradient backpropagation algorithm within the MATLAB neural network training toolbox (\textit{nntraintool}) with its available default options.\\

\begin{figure}[ht!]
  \centering
    \begin{subfigure}[b]{0.49\textwidth}
      \centering
        \includegraphics[width=\textwidth]{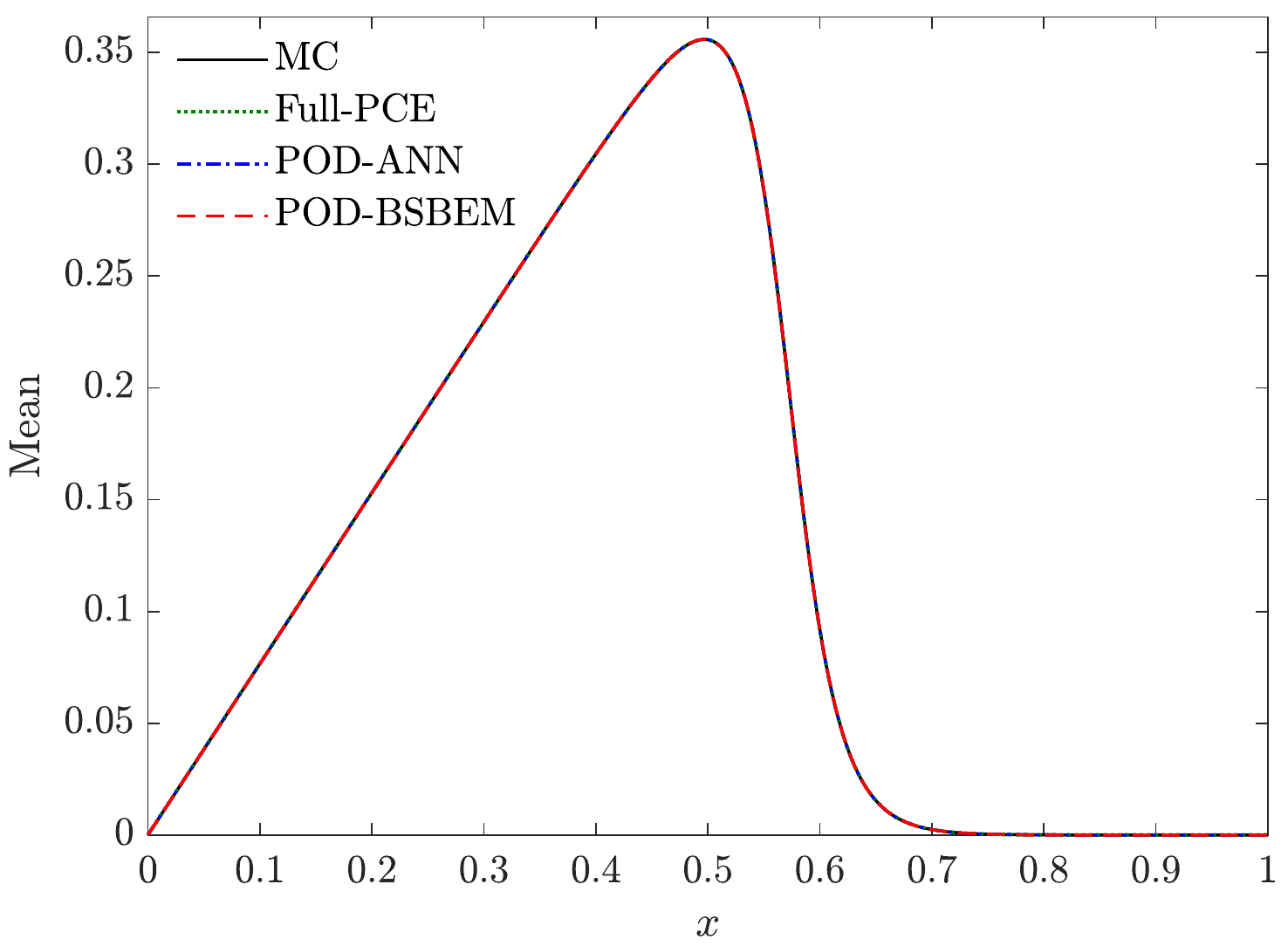}
         \caption{$t\approx0.3$}
         \label{fig:Burg_mean_t_0_3_Re_200}
    \end{subfigure}  
  \hfill
    \begin{subfigure}[b]{0.49\textwidth}
      \centering
        \includegraphics[width=\textwidth]{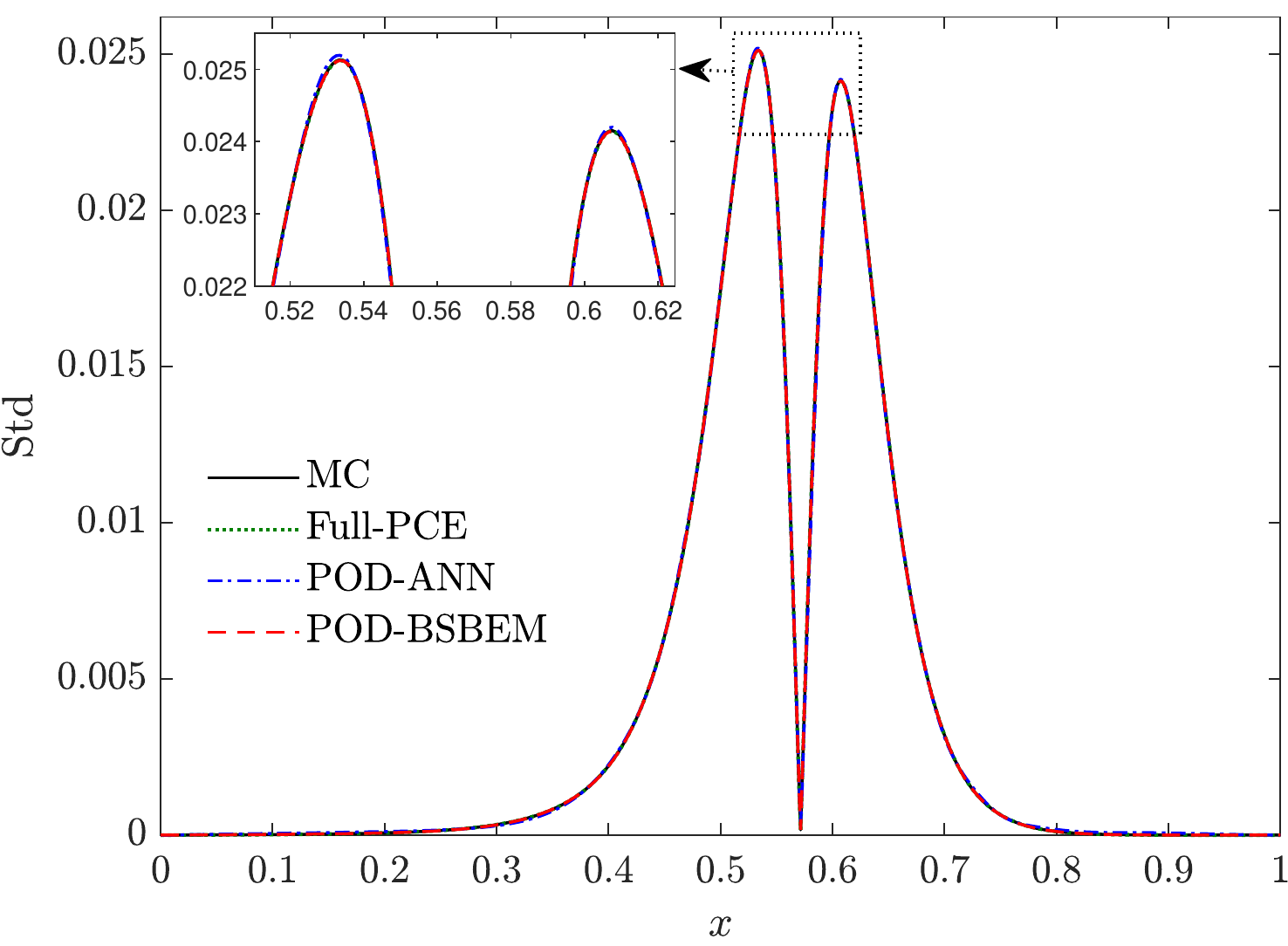}
         \caption{$t\approx0.3$}
         \label{fig:Burg_std_t_0_3_Re_200}
     \end{subfigure} 
     \hfill
    \begin{subfigure}[b]{0.49\textwidth}
      \centering
        \includegraphics[width=\textwidth]{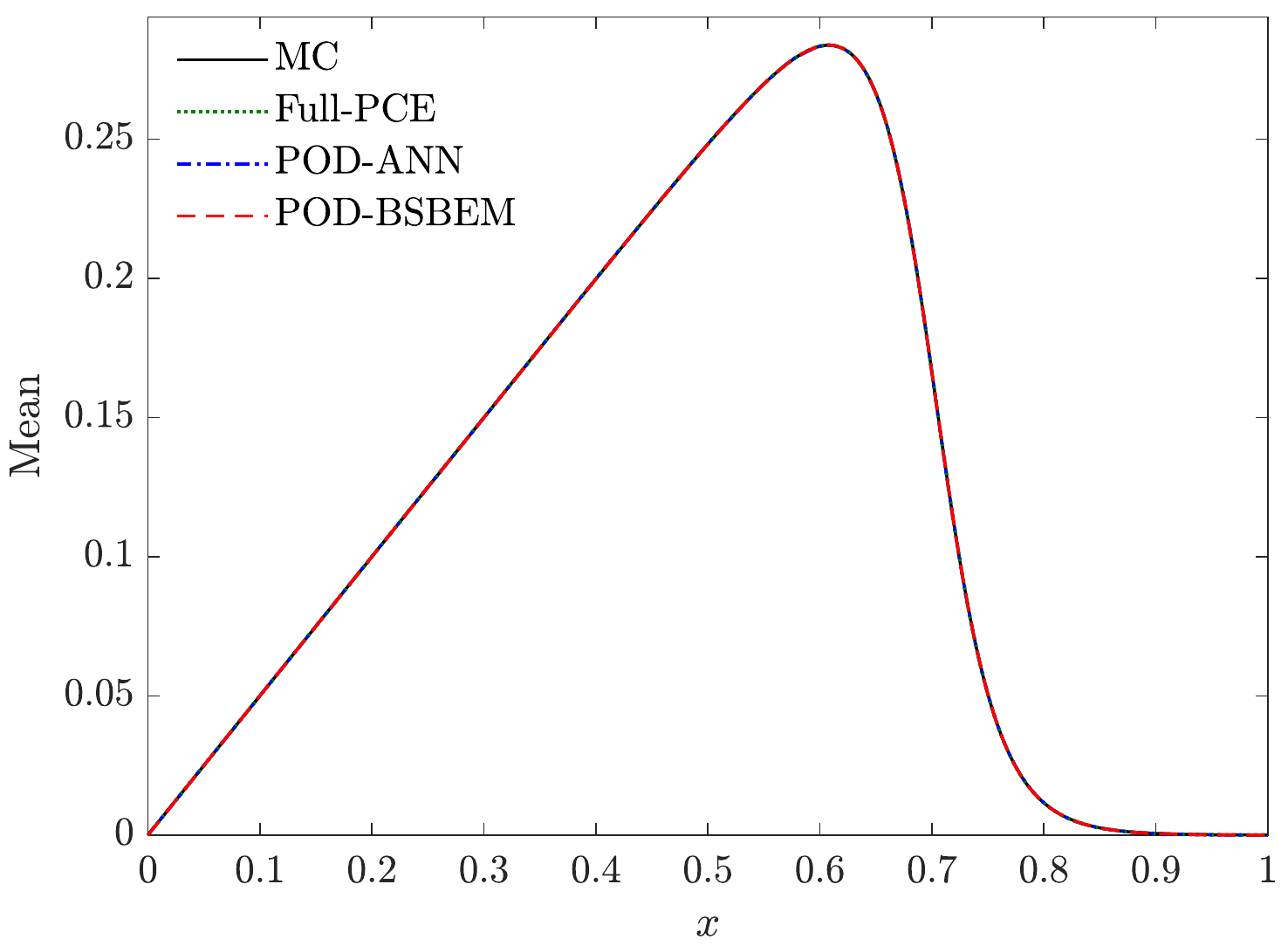}
         \caption{$t\approx1$}
         \label{fig:Burg_mean_t_1_Re_200}
     \end{subfigure} 
     \hfill
    \begin{subfigure}[b]{0.49\textwidth}
      \centering
        \includegraphics[width=\textwidth]{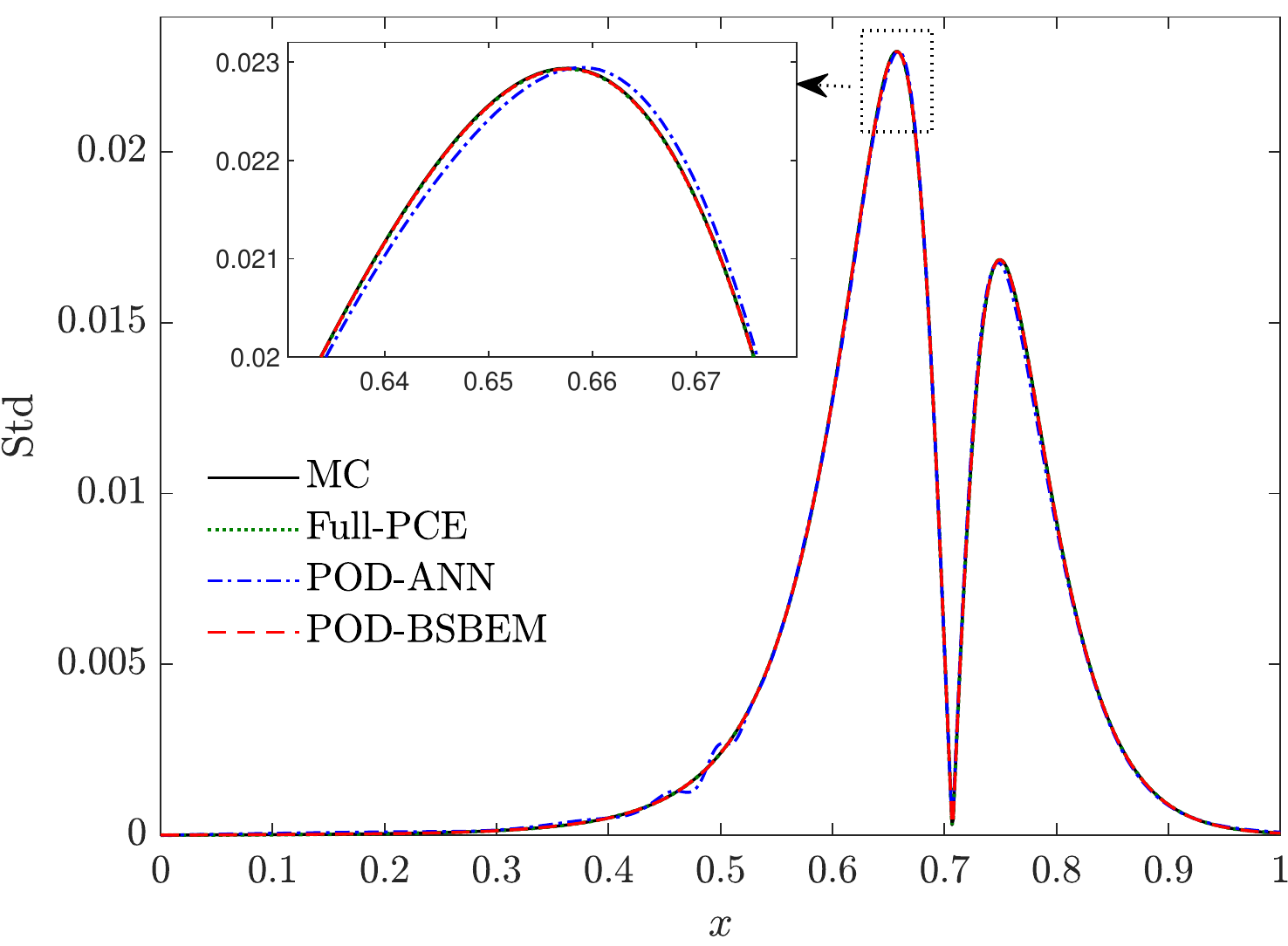}
         \caption{$t\approx1$}
         \label{fig:Burg_std_t_1_Re_200}
     \end{subfigure}    
   \caption{The mean (first column) and standard deviations (second column) of the Burgers' equation  obtained with the POD-BSBEM ($p=2,\;nx=10$), POD-ANN ($H_{1,2,3}=50,\;N_{s}=300$) and Full-PCE ($p=6,\;n_{p}=2$) compared to those of the MC solution (with $10^{6}$ realizations) for $\mu_{Re=200}$ and $cv=25\,\%$ at different times. (a and b): $t\approx0.3$ and (c and d): $t\approx1$.}
   \label{fig:Burg_mean_std_comp_time_Re_200}
\end{figure}

A comparison of the mean and standard deviation profiles as a function of the $x$-coordinate at different times ($t\approx0.3,\,1$) obtained with POD-BSBEM (with $p=2,\,nx=10$), POD-ANN (with $H_{1,2,3}=50,\,N_{s}=300$), and Full-PCE (with $p=6,\,n_{p}=2$) with those from the MC reference solution (with $N_{s}=10^{6}$ realizations) for $\mu_{Re}=200$ and $cv=25\,\%$ is depicted in Fig.\ref{fig:Burg_mean_std_comp_time_Re_200}. In addition to the high-fidelity MC method, the Full-PCE approach, which is considered as a stochastic expansion full-order model, is introduced in order to assess the contribution of the reduced-order aspect of the proposed POD-BSBEM in comparison with that of the POD-ANN. It can be observed that the statistical moments' distributions of the $u$-velocity obtained by the proposed non-intrusive reduced order model (POD-BSBEM) are in excellent agreement with the distributions from the full-order PCE and the high-fidelity MC reference solution. The variability of the relatively smooth output response of the Burgers' solution is localized around the moving shock wave, whose amplitude decreases as a function of time. Slight deviations in this output response can be observed in the POD-ANN solution as time evolves, conversely to the POD-BSBEM predictions, which present profiles nearly indistinguishable from those of the reference MC solution, as shown by the close-up views of the standard deviation profiles (Figs.\ref{fig:Burg_std_t_0_3_Re_200} and \ref{fig:Burg_std_t_1_Re_200}).\\

\begin{figure}[ht!]
  \centering
    \begin{subfigure}[b]{0.49\textwidth}
      \centering
        \includegraphics[width=\textwidth]{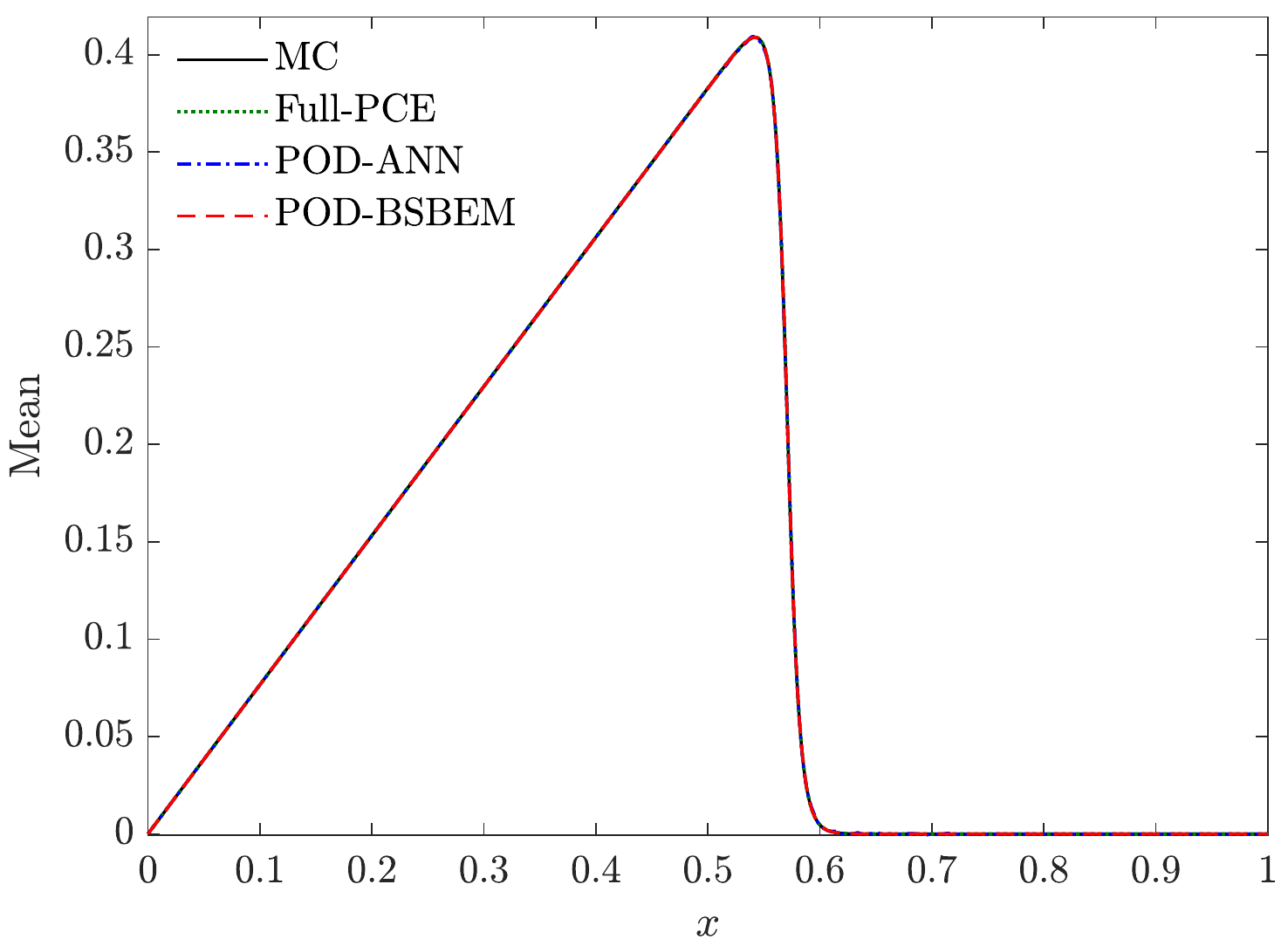}
         \caption{$t\approx0.3$}
         \label{fig:Burg_mean_t_0_3_Re_800}
    \end{subfigure}  
  \hfill
    \begin{subfigure}[b]{0.49\textwidth}
      \centering
        \includegraphics[width=\textwidth]{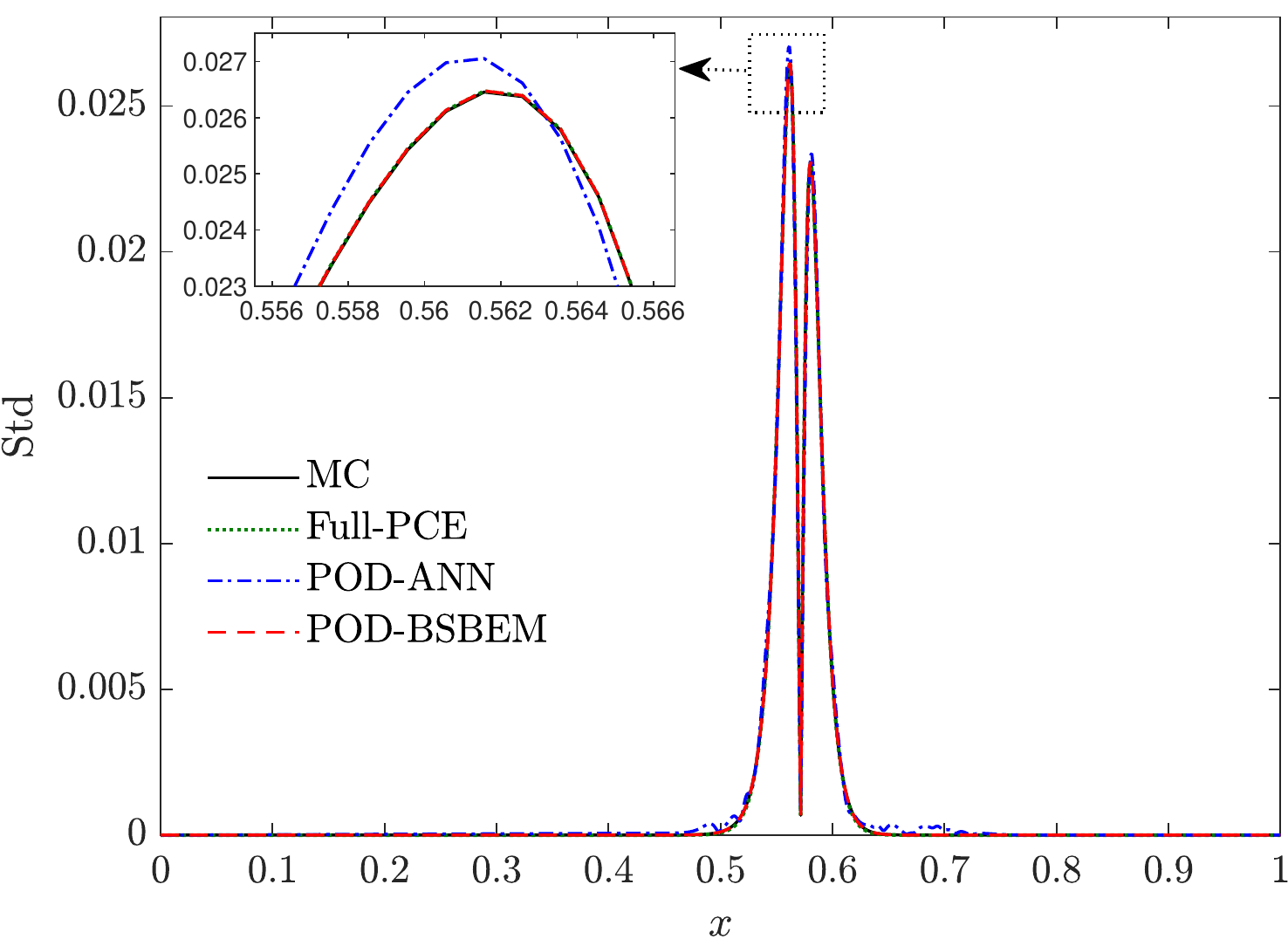}
         \caption{$t\approx0.3$}
         \label{fig:Burg_std_t_0_3_Re_800}
     \end{subfigure} 
     \hfill
    \begin{subfigure}[b]{0.49\textwidth}
      \centering
        \includegraphics[width=\textwidth]{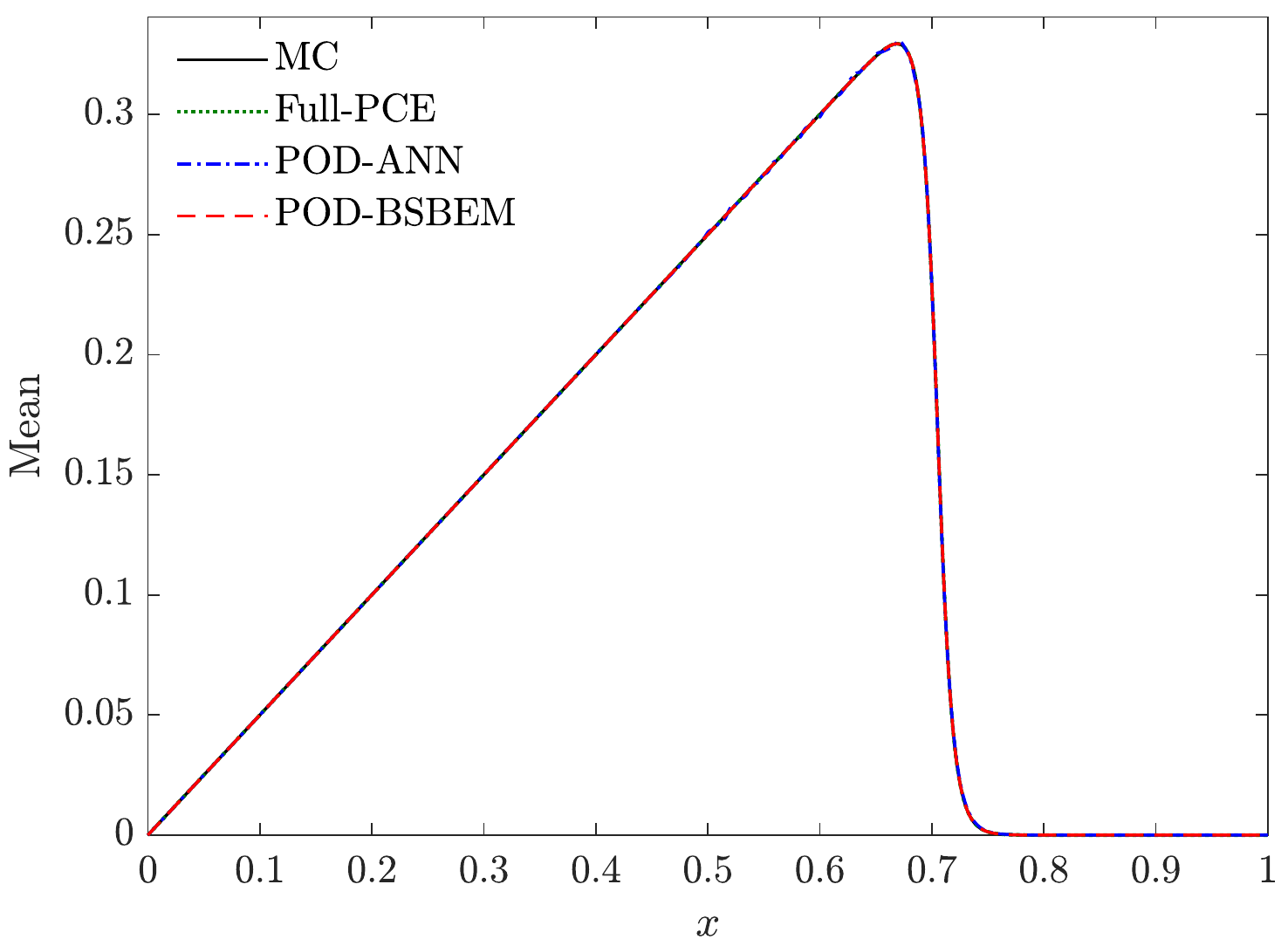}
         \caption{$t\approx1$}
         \label{fig:Burg_mean_t_1_Re_800}
     \end{subfigure} 
     \hfill
    \begin{subfigure}[b]{0.49\textwidth}
      \centering
        \includegraphics[width=\textwidth]{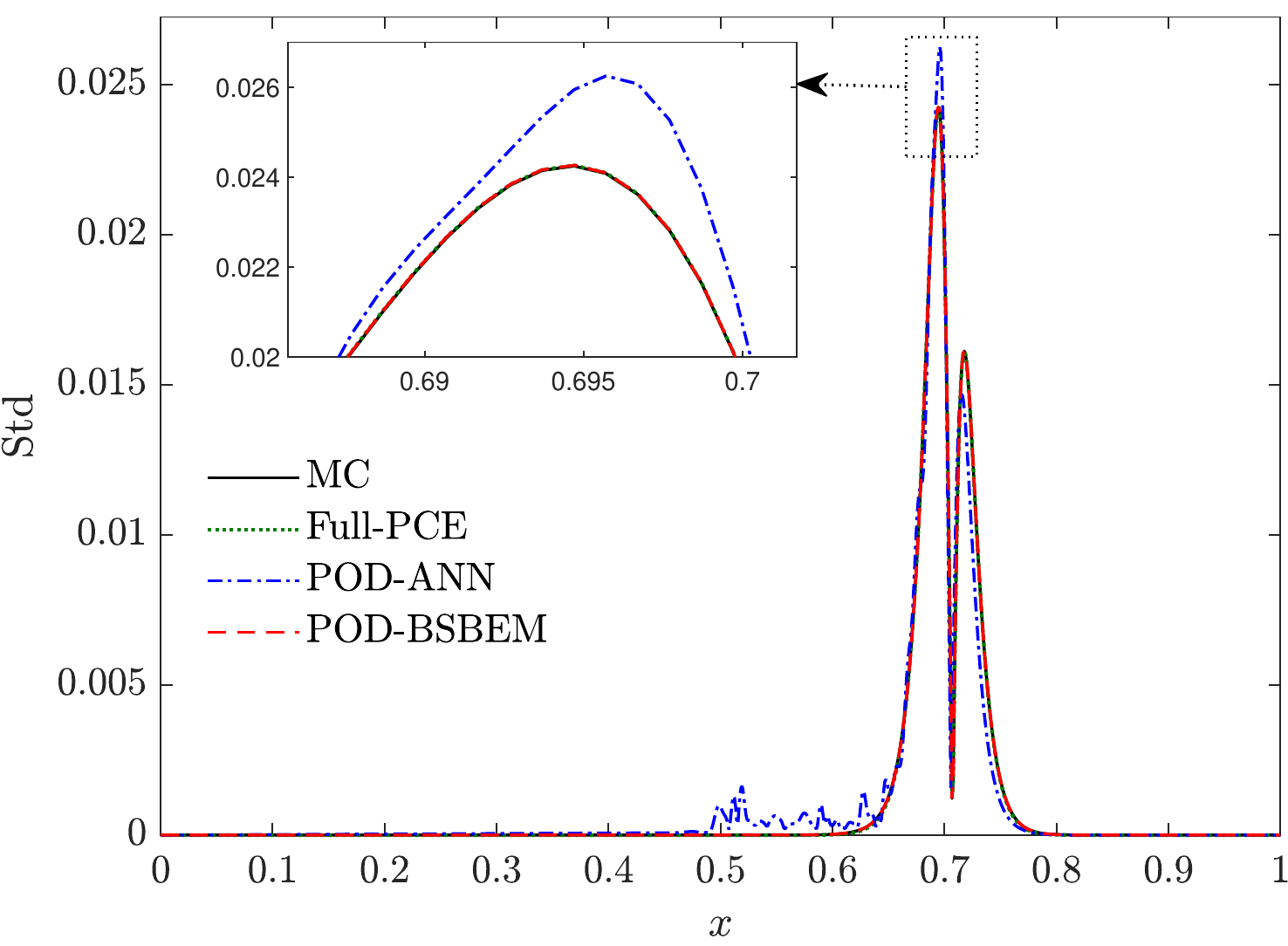}
         \caption{$t\approx1$}
         \label{fig:Burg_std_t_1_Re_800}
     \end{subfigure}    
   \caption{The mean (first column) and standard deviations (second column) of the Burgers' equation obtained with the POD-BSBEM ($p=2,\;nx=10$), POD-ANN ($H_{1,2,3}=50,\;N_{s}=300$) and Full-PCE ($p=6,\;n_{p}=2$), compared to those of the MC solution (with $10^{6}$ realizations) for $\mu_{Re=800}$ and $cv=25\,\%$ at different times. (a and b): $t\approx0.3$, and (c and d): $t\approx1$.}
   \label{fig:Burg_mean_std_comp_time_Re_800}
\end{figure}

A similar comparison is presented in Fig.\ref{fig:Burg_mean_std_comp_time_Re_800} for relatively high values of the uncertain input parameter with $\mu_{Re}=800$. This increase in the mean value of the Reynolds number associated with a wide variability range ($cv=25\,\%$) leads to a visible transformation of the output response in the Burgers' solution, where a more pronounced steep slope can be observed throughout  the velocity front wave  propagation. These conditions represent a real challenge for the reduced-order models to faithfully reproduce the output response of the quantities of interest. While both reduced-order models (POD-BSBEM and POD-ANN) accurately estimate the mean of the Burgers' solution where the profiles are in  good accordance with those from Full-PCE and MC, as depicted in the first column of Fig.\ref{fig:Burg_mean_std_comp_time_Re_800}, it is clear that the POD-ANN approach shows deviations in the estimation of the peak value of the standard deviation profile where oscillations, which gain in intensity over  time, can be observed. In contrast, the predictions obtained by the proposed POD-BSBEM are in excellent accordance with those from the Full-PCE and the high-fidelity MC reference solution, as shown by the close-up views in the second column of Fig.\ref{fig:Burg_mean_std_comp_time_Re_800}.\\

\begin{figure}[ht!]
  \centering
    \begin{subfigure}[b]{0.49\textwidth}
      \centering
        \includegraphics[width=\textwidth]{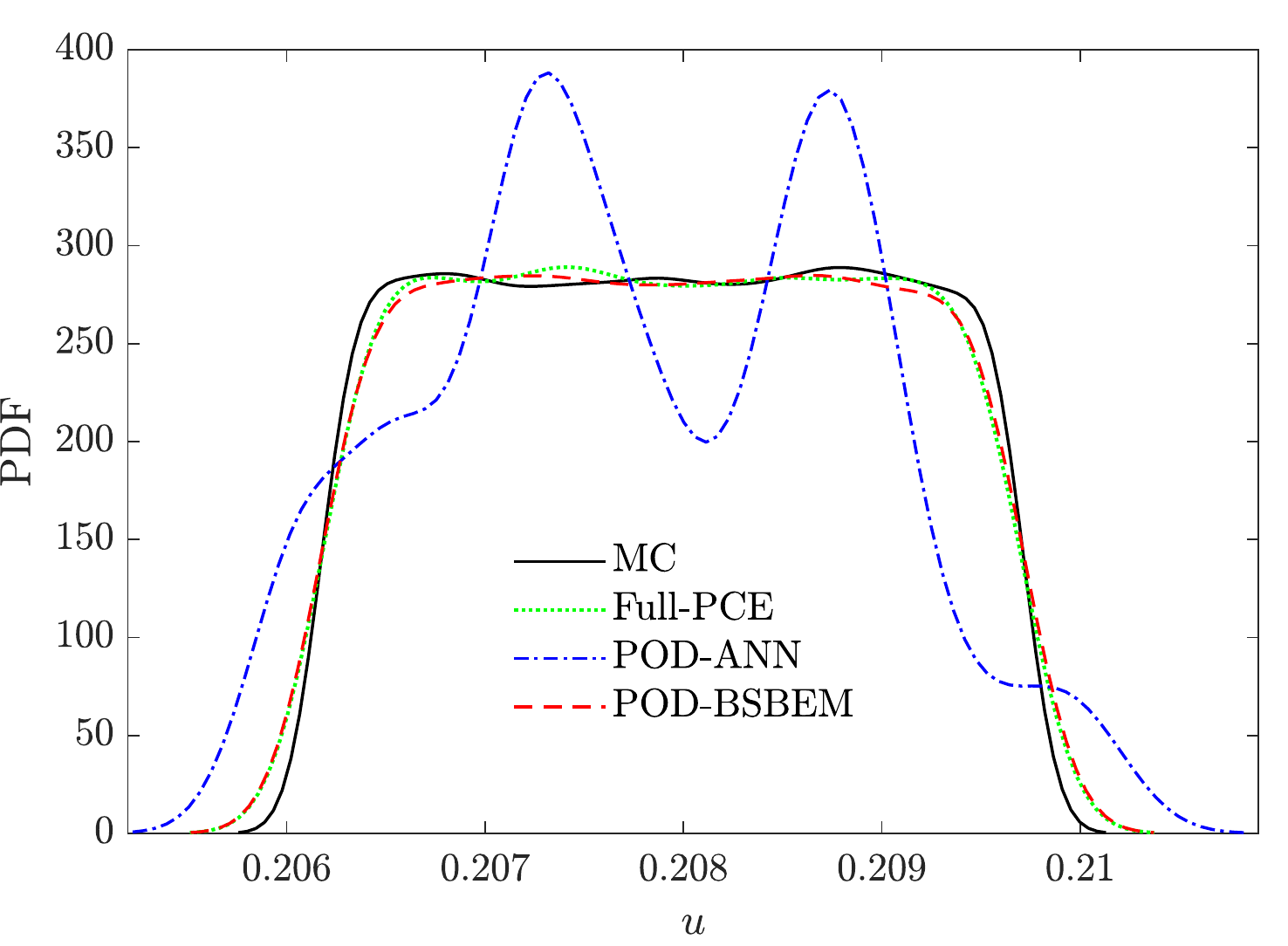}
         \caption{$\mu_{Re}=200,\,t\approx0.3,\,x\approx0.57$}
         \label{fig:Burg_pdf_Re_200_t_0_3_x_0_57}
    \end{subfigure}  
  \hfill
    \begin{subfigure}[b]{0.49\textwidth}
      \centering
        \includegraphics[width=\textwidth]{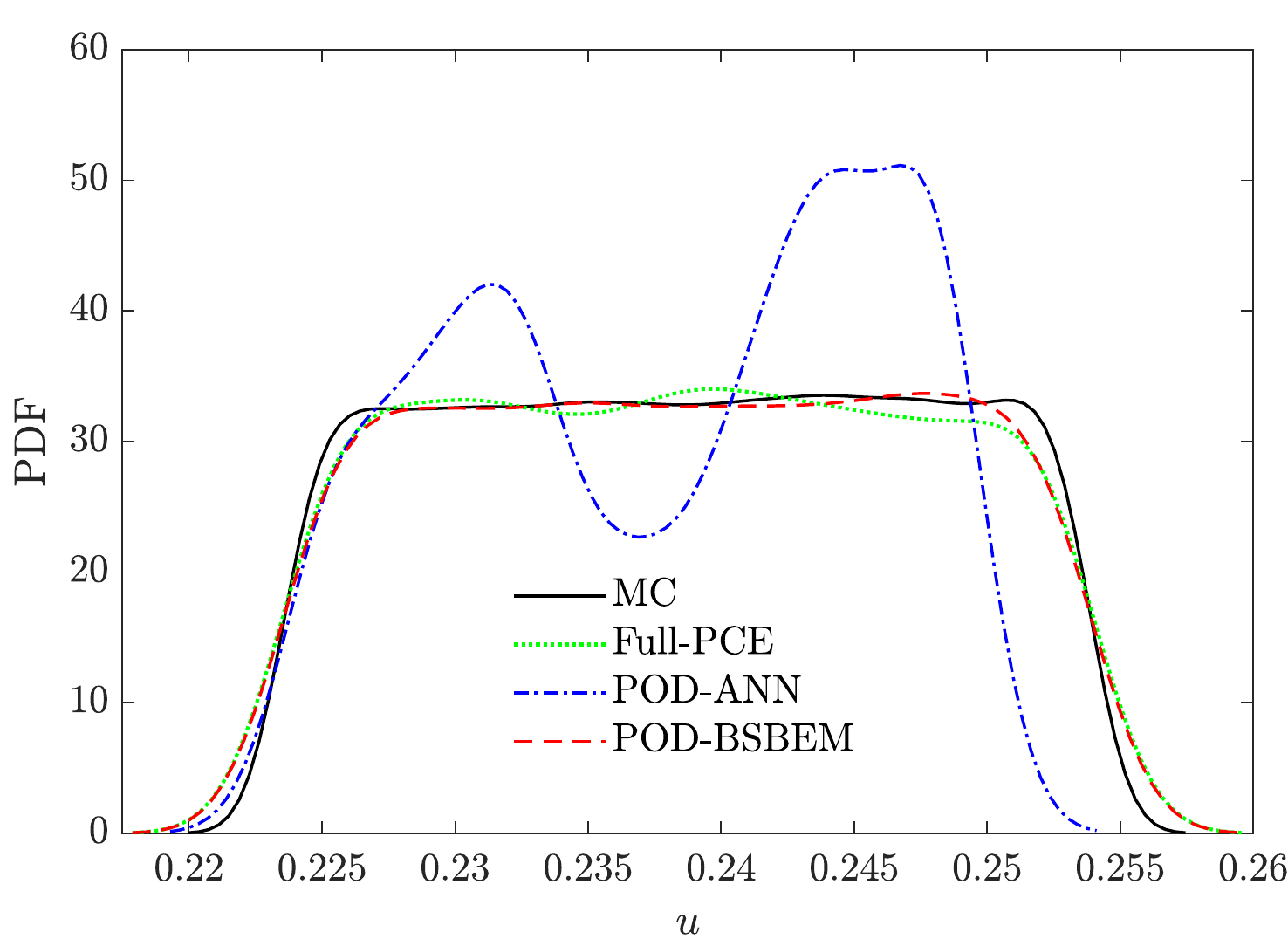}
         \caption{$\mu_{Re}=800,\,t\approx0.3,\,x\approx0.57$}
         \label{fig:Burg_pdf_Re_800_t_0_3_x_0_57}
     \end{subfigure} 
    \hfill
    \begin{subfigure}[b]{0.49\textwidth}
      \centering
        \includegraphics[width=\textwidth]{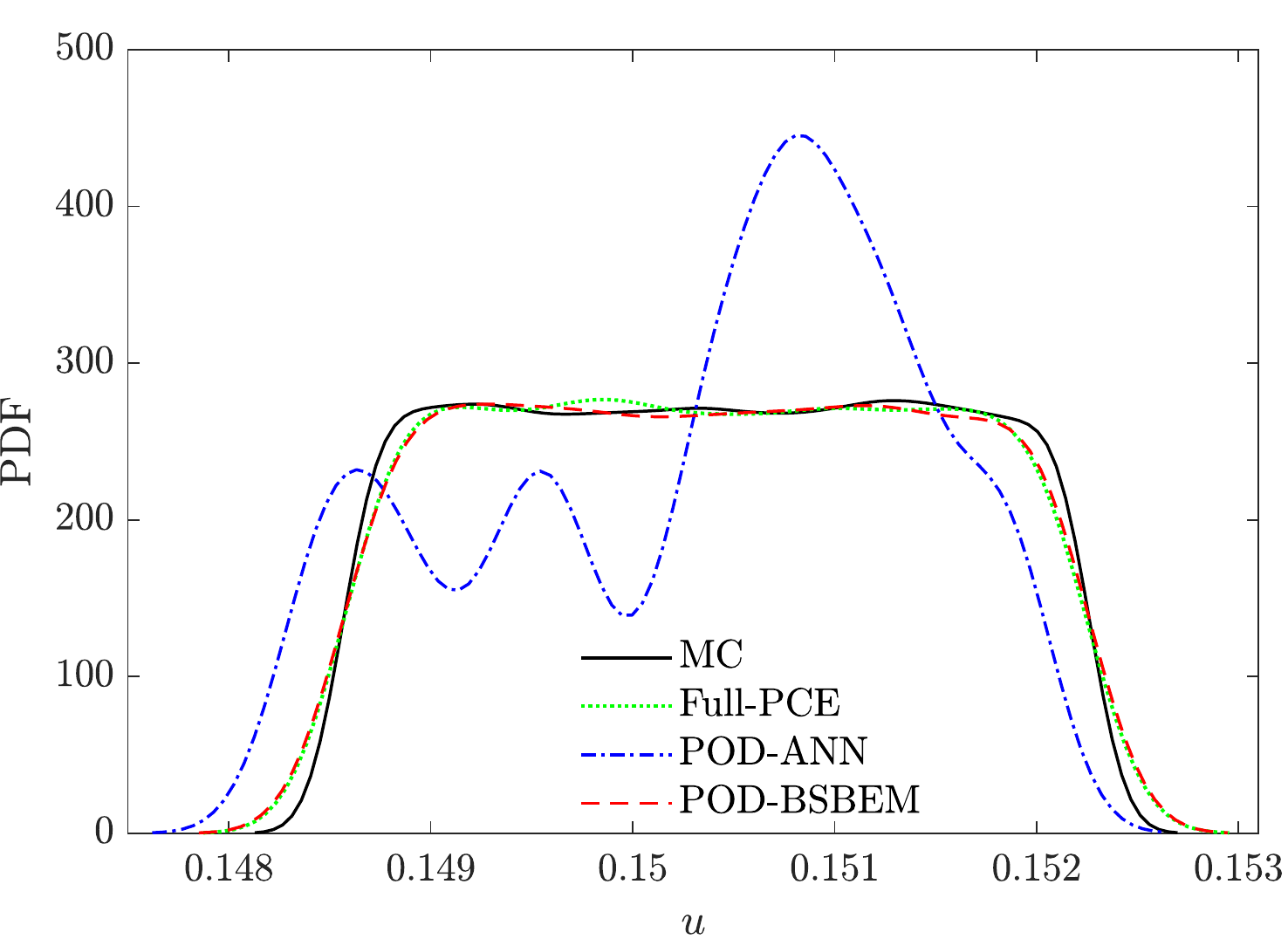}
         \caption{$\mu_{Re}=200,\,t\approx1,\,x\approx0.7$}
         \label{fig:Burg_pdf_Re_200_t_1_x_0_7}
     \end{subfigure} 
     \hfill
    \begin{subfigure}[b]{0.49\textwidth}
      \centering
        \includegraphics[width=\textwidth]{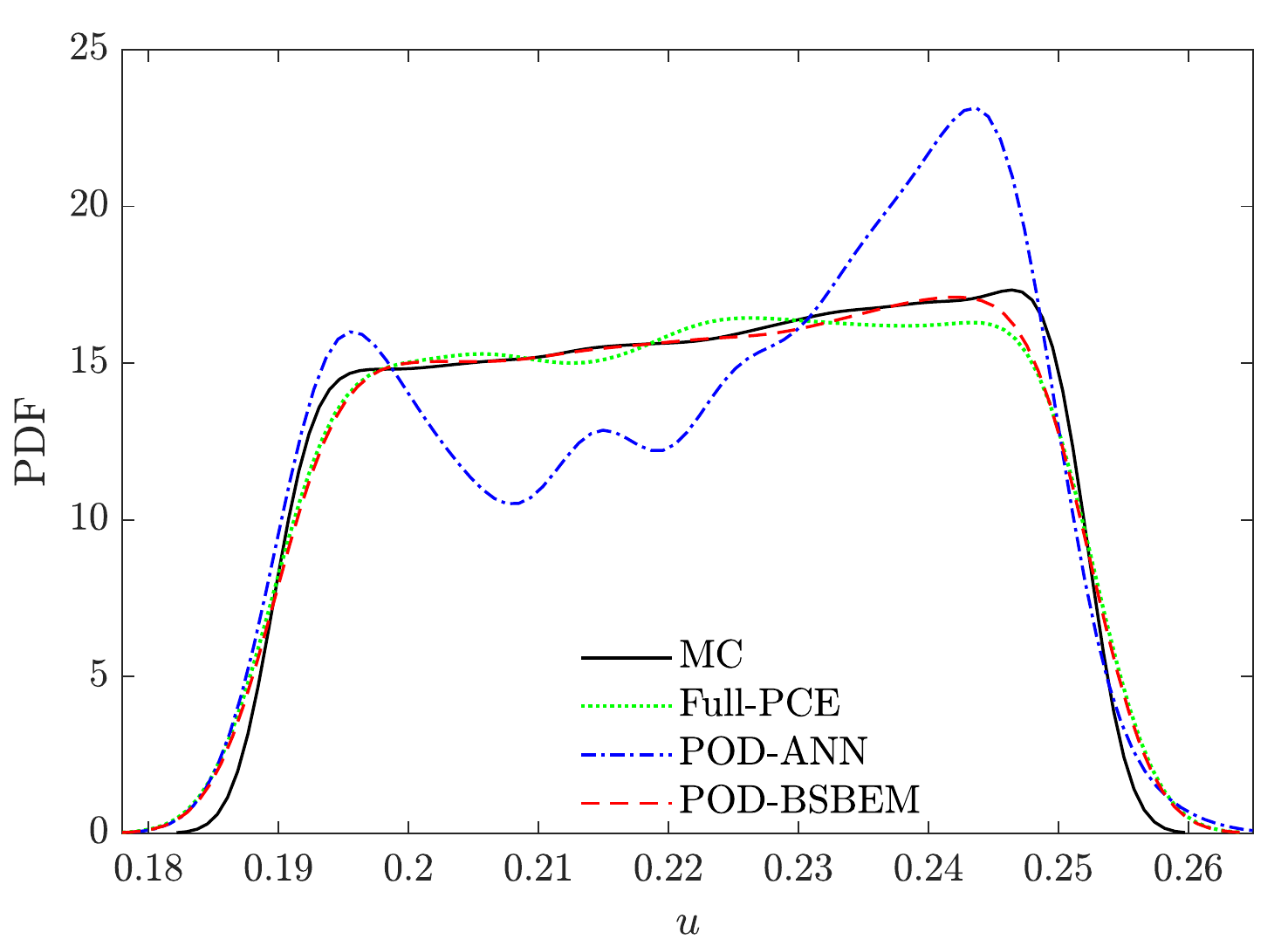}
         \caption{$\mu_{Re}=800,\,t\approx1,\,x\approx0.7$}
         \label{fig:Burg_pdf_Re_800_t_1_x_0_7}
     \end{subfigure}    
   \caption{Kernel density estimated probability density function (PDF) of Burgers' solutions obtained with POD-BSBEM ($p=2,\;nx=10$), POD-ANN ($H_{1,2,3}=50,\;N_{s}=300$ and Full-PCE ($p=6,\;n_{p}=2$) compared to the MC reference solution (with $10^{6}$ realizations) for $\mu_{Re=200}$ (first column) and $\mu_{Re=800}$ (second column). (a and b): $t\approx0.3,\,x\approx0.57$, and (c and d): $t\approx1,\,x\approx0.7$.}
   \label{fig:Burg_pdfs_Re_200_800}
\end{figure}

Another meaningful aspect for which the efficacy of the proposed POD-BSBEM is assessed, in addition to the estimation of the first two statistical moments, is its ability to accurately estimate the probability density function of the output response of quantities of interest comparison to the estimations of the other aforementioned techniques. To visualize this aspect, the distributions of the estimated kernel density function obtained by POD-BSBEM ($p=2,\,nx=10$), POD-ANN ($H_{1,2,3}=50,\,N_{s}=300$) and Full-PCE ($p=6,\,n_{p}=2$) using the constructed surrogate models are compared with those from the $10^{5}$ Monte Carlo realizations at different spatio-temporal locations for $\mu_{Re}=200$ and $800$ with $cv=25\,\%$, are depicted in Fig.\ref{fig:Burg_pdfs_Re_200_800}. It should be noted that the spatial locations where the pdfs are estimated correspond approximately to the positions of the front wave of the Burgers’ solution for a given $t$-value. From Fig.\ref{fig:Burg_pdfs_Re_200_800}, it is obvious  that all the pdf plots resulting from the proposed POD-BSBEM follow nearly the same trend of the pdfs obtained from the Full-PCE and the MC high-fidelity reference solution, whose shapes are almost uniform, unlike the POD-ANN approach which presents clear discrepancies in its pdf plots for all the represented $x$ and $t$ locations. The effect of increasing the mean value of the input Reynolds number from $\mu_{Re}=200$ to $\mu_{Re}=800$ with a coefficient of variation of $cv=25\,\%$, which is traduced physically by a steep solution of the Burgers' equation, is visible in the variability range of the $u$-output response that widens as time evolves, as can be seen in the second column of Fig.\ref{fig:Burg_pdfs_Re_200_800}.\\

\begin{figure}[ht!]
  \centering
    \begin{subfigure}[b]{0.49\textwidth}
      \centering
        \includegraphics[width=\textwidth]{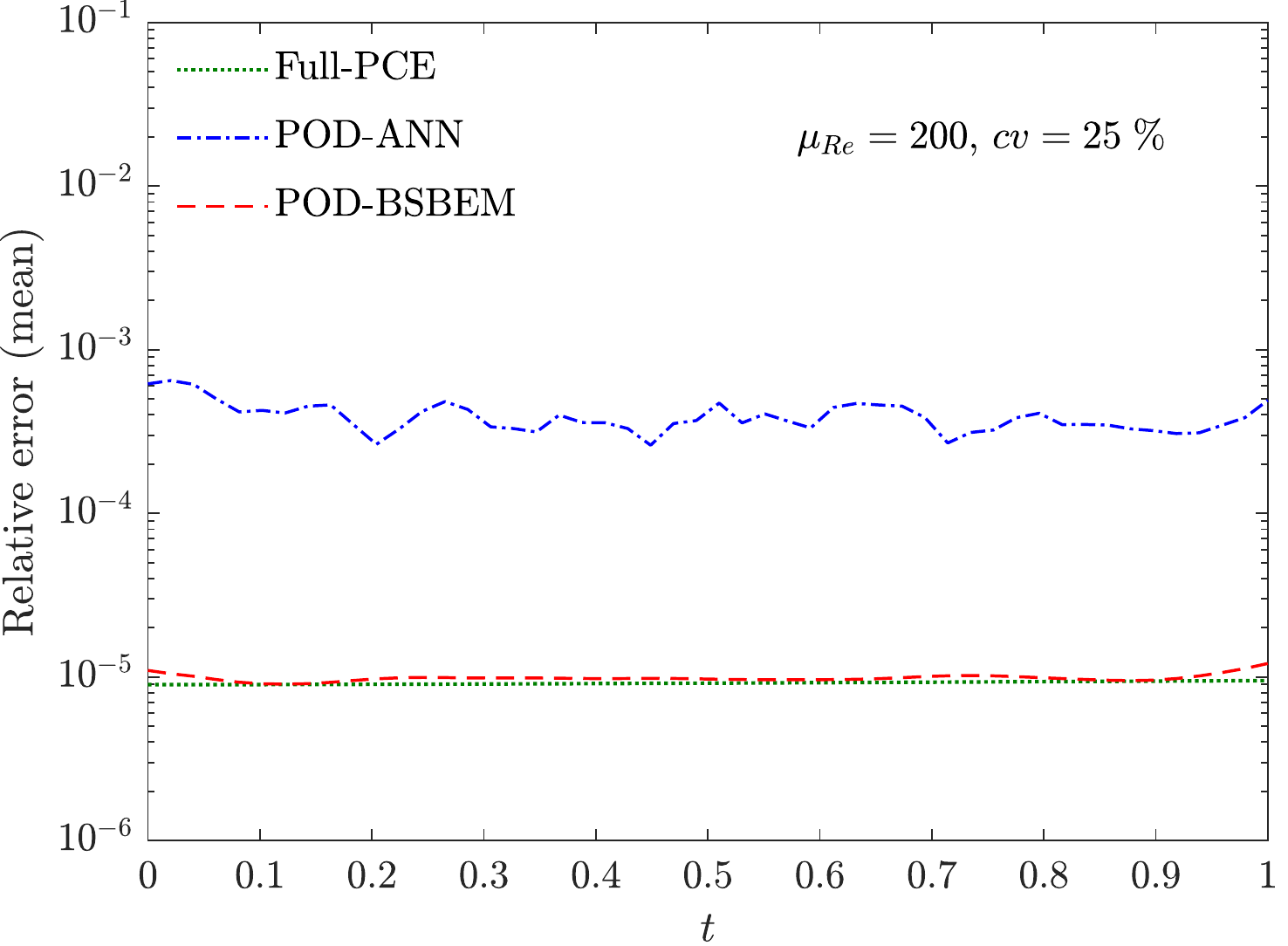}
         \caption{$\mu_{Re}=200,\;cv=25\;\%$}
         \label{fig:Burg_L2_mean_Re_200}
    \end{subfigure}  
  \hfill
    \begin{subfigure}[b]{0.49\textwidth}
      \centering
        \includegraphics[width=\textwidth]{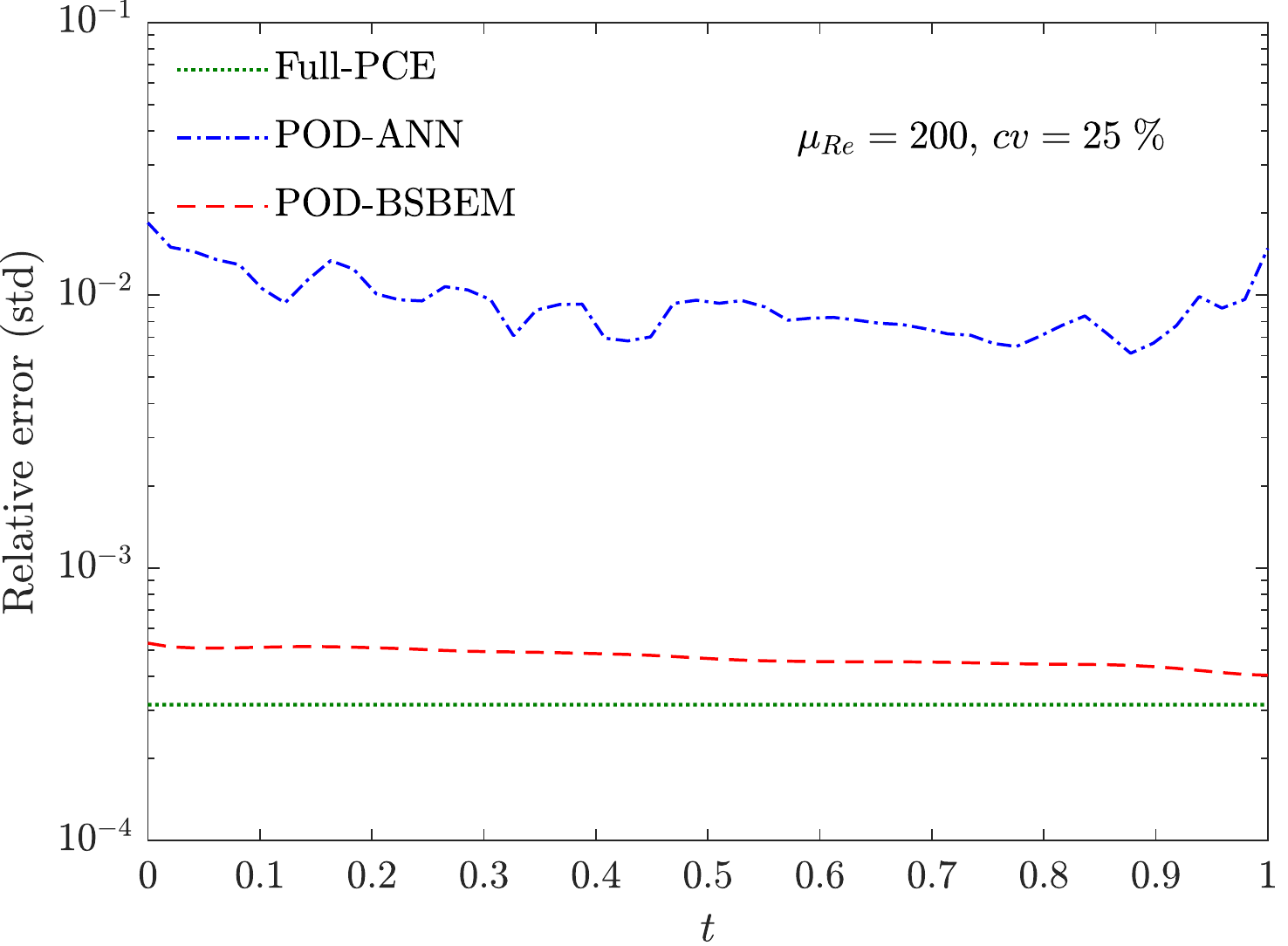}
         \caption{$\mu_{Re}=200,\;cv=25\;\%$}
         \label{fig:Burg_L2_std_Re_200}
     \end{subfigure} 
     \hfill
    \begin{subfigure}[b]{0.49\textwidth}
      \centering
        \includegraphics[width=\textwidth]{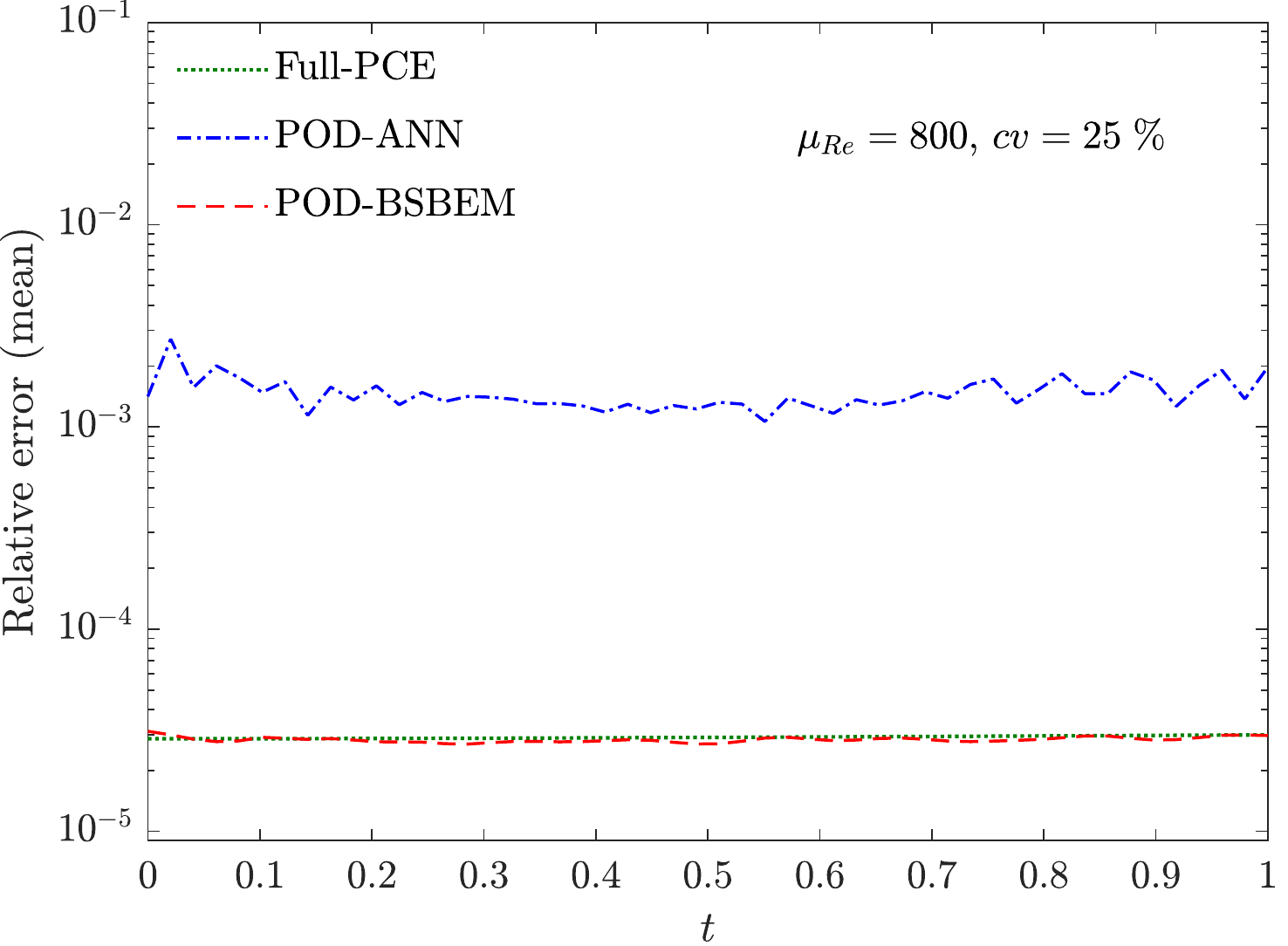}
         \caption{$\mu_{Re}=800,\;cv=25\;\%$}
         \label{fig:Burg_L2_mean_Re_800}
     \end{subfigure} 
     \hfill
    \begin{subfigure}[b]{0.49\textwidth}
      \centering
        \includegraphics[width=\textwidth]{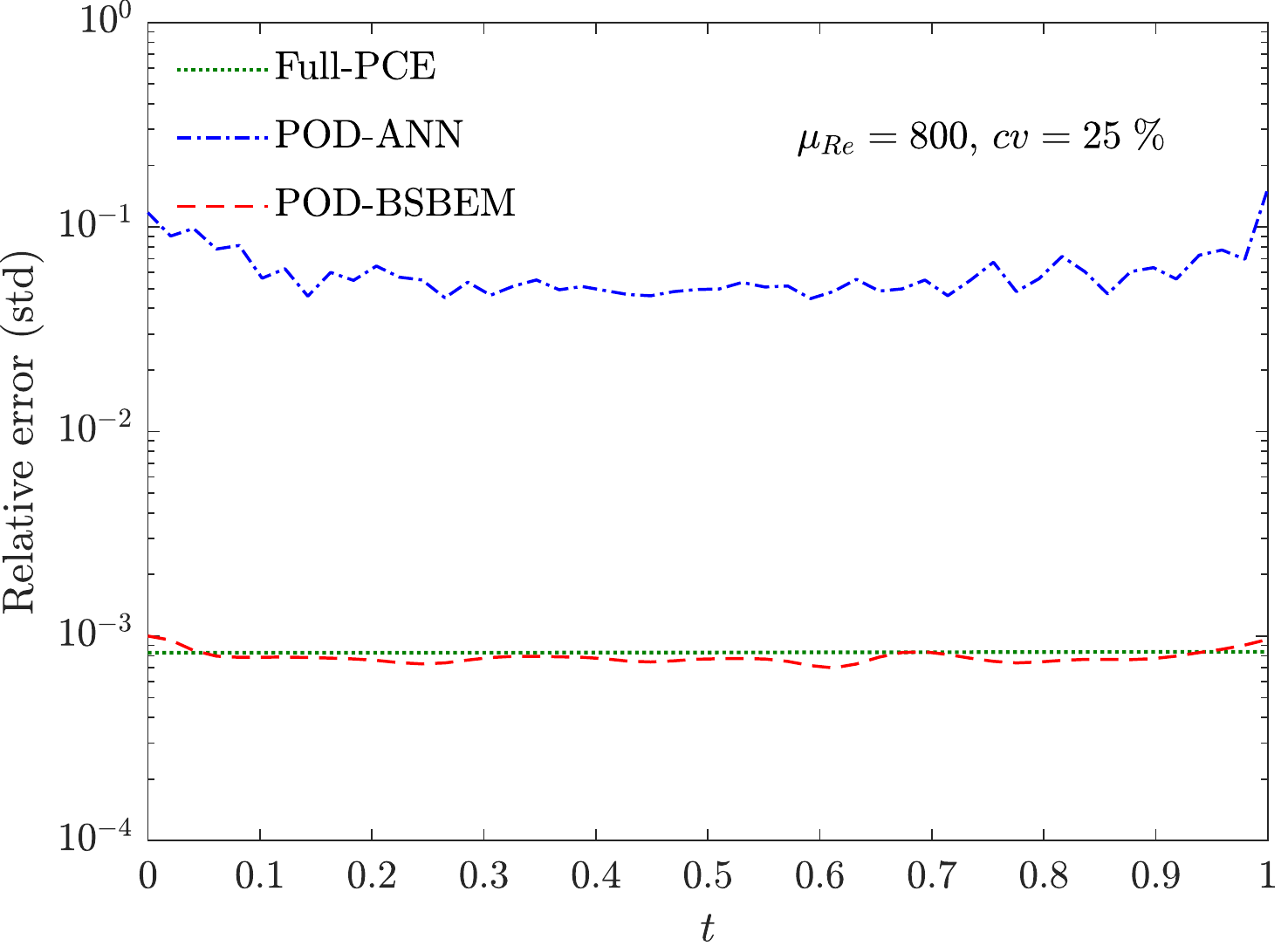}
         \caption{$\mu_{Re}=800,\;cv=25\;\%$}
         \label{fig:Burg_L2_std_Re_800}
     \end{subfigure}  
   \caption{Comparison of the relative $L^{2}$-error profiles of the mean (first column) and standard deviations (second column) as a function of time, obtained with the Full-PCE ($p=6,\;n_{p}=6$), POD-ANN ($H_{1,2,3}=50,\;N_{s}=300$) and POD-BSBEM ($p=2,\;nx=10$). Errors are computed with respect to the MC reference solution (with $10^{6}$ realizations). (a and b): $\mu_{Re}=200,\;cv=25\,\%$, (c and d): $\mu_{Re}=800,\;cv=25\,\%$.}
   \label{fig:Burg_L2_err_Re_200_800}
\end{figure}

Fig.\ref{fig:Burg_L2_err_Re_200_800} presents a complementary quantitative assessment of the efficiency of the proposed POD-BSBEM technique was performed via the $L^{2}-$error analysis. The variation of the $L^{2}-$error profiles of the mean and standard deviation, whose mathematical expression in its time-dependent form is given by the right term of Eq.\eqref{<L2_max>}, resulting from the POD-BSBEM, POD-ANN and Full-PCE for $\mu_{Re}=200$ and $\mu_{Re}=800$ with $cv=25\,\%$. The relative error in the $L^{2}$ norm is computed with respect to the MC high fidelity reference solution for each time step and over  the $N_{e}$ nodes that contains the computational domain. It can be seen that the predicted errors of the POD-BSBEM are small and have the same order of magnitude ($\approx10^{-6}$ and $10^{-3}$ for the mean and standard deviation, respectively) as the Full-PCE technique, unlike the predicted errors from the POD-ANN, particularly for the standard deviation, which are on the order of $10^{-1}$ when the input Reynolds number takes relatively high values. This indicates that the proposed POD-BSBEM yields accurate predictions of statistical moments of both smooth and steep time-dependent Burgers' solutions.\\

The effect of the number of POD modes on the performance of the proposed POD-BSBEM approach is assessed through the maximum value of the time-dependent relative error in the $L^{2}$-norm, retained as a criterion for evaluating and comparing the accuracy of the different models. Its mathematical formulation is expressed as follows:
\begin{equation}\label{<L2_max>}
Err_{L^{2}}^{max}=max\left(Err_{L^{2},\Phi}^{Surr}(t)= \sqrt{\frac{\sum_{i=1}^{N_{e}}\left(\Phi_{i,Surr}(t)-\Phi_{i,MC}(t)\right)^2}{\sum_{i=1}^{N_{e}}(\Phi_{i,MC}(t))^2}}\right) 
\end{equation}
The $L^{2}$-error mentioned above is evaluated with respect to the MC reference solution over all the computational nodes for each time step. The variation of the $Err_{L^{2}}^{max}$ error as a function of the number of POD modes $L$ for two variability ranges of the uncertain input Reynolds number is reported in Tables \ref{tab:max_L2error_Re_200_epsi} and \ref{tab:max_L2error_Re_800_epsi}, respectively. It should be mentioned that the $Err_{L^{2}}^{max}$ of the Full-PCE approach is not affected by the variation of the POD modes' number since it represents a full-order model; this variation was included in the aforementioned tables to serve  as a reference with which to compare the accuracy of the reduced-order models. \\

Thus, for relatively low values of the input Reynolds number, the results reported in Table \ref{tab:max_L2error_Re_200_epsi} show that errors in both statistical moments decrease as the number of POD modes increases (as a result of decreasing the threshold $\epsilon_{s}=\epsilon_{t}$), where very low values can be observed, particularly in the mean error. However, a clear discrepancy can be noted in the standard deviation errors obtained with $L=34$ modes (which corresponds to $\epsilon_{s}=\epsilon_{t}=10^{-10}$) between the POD-ANN and POD-BSBEM techniques, which are on the order of $10^{-2}$ and $10^{-4}$, respectively. The same trend can be observed for relatively high values of the Reynolds number with a wide variability range, where the number of POD modes, for the same threshold value, is higher than in the former case (low Reynolds number values). The predicted errors obtained from the POD-BSBEM, which are of the same degree of magnitude as those of the Full-PCE, are much lower than the predicted errors of the POD-ANN approach for both the mean and the standard deviation, as reported in Table \ref{tab:max_L2error_Re_800_epsi}. This quantitative analysis clearly illustrates the ability of the proposed model to  accurately estimate the time-dependent stochastic output response of quantities of interest.\\

\begin{table}[h!]
\caption{Effect of $\epsilon_{s}$ and $\epsilon_{t}$ on the maximum relative error over time in the $L^{2}$-norm for the mean and standard deviation obtained from Full-PCE ($p=6$, $n_{p}=2$), POD-ANN ($H_{i\in\{1,2,3\}}=50$, $N_{s}=300$) and POD-BSBEM ($p=2$, $nx=10$). Errors are computed with respect to the MC reference solution (with $N_{s}=10^{6}$ realizations) for $Re_{\mu=200,\,\sigma=50}\in\mathcal{U}\left[114,\;287\right]$.}
\centering
\begin{tabular*}{0.85\textwidth}{@{\extracolsep{\fill}}l l  llll }
\hline        
& & \multicolumn{2}{c}{$Err_{L^{2},\,Mean}^{max}$}& \multicolumn{2}{c}{$Err_{L^{2},\,Std}^{max}$}\\
\cline{3-4}
\cline{5-6}
\raisebox{2ex}{$ \epsilon_{s},\;\epsilon_{s} $} & \raisebox{2ex}{$ L $} &POD-ANN& POD-BSBEM& POD-ANN& POD-BSBEM\\
\hline        
 $10^{-3}$& $7$& $0.008882$& $0.009098 $& $0.211378$& $0.409364$\\
$10^{-5}$& $13$& $8.4564$E-04& $5.3330$E-04& $0.024265$& $0.023550$\\
$10^{-8}$& $24$& $6.1322$E-04& $1.1664$E-04& $0.017521$& $7.9154$E-04\\
$10^{-10}$& $34$& $8.9175$E-04& $1.2082$E-05& $0.021630$& $5.2967$E-04\\

\hline        
\multicolumn{2}{l}{Full-PCE}&\multicolumn{2}{c}{$1.3006$E-05}&\multicolumn{2}{c}{$5.4106$E-04}\\
\hline 
\end{tabular*}
\label{tab:max_L2error_Re_200_epsi}
\end{table}

\begin{table}[h!]
\caption{Effect of $\epsilon_{s}$ and $\epsilon_{t}$ on the maximum relative error over time in the $L^{2}$-norm for the mean and standard deviation obtained from Full-PCE ($p=6$, $n_{p}=2$), POD-ANN ($H_{i\in\{1,2,3\}}=50$, $N_{s}=300$) and POD-BSBEM ($p=2$, $nx=10$). Errors are computed with respect to the MC reference solution (with $N_{s}=10^{6}$ realizations) for $Re_{\mu=800,\,\sigma=200}\in\mathcal{U}\left[454,\;1\,146\right]$.}
\centering
\begin{tabular*}{0.85\textwidth}{@{\extracolsep{\fill}}l l  llll }
\hline        
& & \multicolumn{2}{c}{$Err_{L^{2},\,Mean}^{max}$}& \multicolumn{2}{c}{$Err_{L^{2},\,Std}^{max}$}\\
\cline{3-4}
\cline{5-6}
\raisebox{2ex}{$ \epsilon_{s},\;\epsilon_{s} $} & \raisebox{2ex}{$ L $} &POD-ANN& POD-BSBEM& POD-ANN& POD-BSBEM\\
\hline        
 $10^{-3}$& $12$& $0.020304$& $0.019633$& $0.653037$& $0.707108$\\
$10^{-5}$& $26$& $0.002284$& $0.001225$& $0.109855$& $0.003519$\\
$10^{-8}$& $53$& $0.002546$& $4.5725$E-05& $0.118874$& $0.002541$\\
$10^{-10}$& $76$& $0.002233$& $3.0940$E-05& $0.100247$& $0.001056$\\

\hline        
\multicolumn{2}{l}{Full-PCE}&\multicolumn{2}{c}{$3.0616$E-05}&\multicolumn{2}{c}{$0.001110$}\\
\hline 
\end{tabular*}
\label{tab:max_L2error_Re_800_epsi}
\end{table}

As with the Ackley test case, the estimation of the computational costs of the proposed approach are compared to  those of other techniques for this time-dependent test case; these are presented in Table \ref{tab:Tab_comp_cpu_Burg}. The computations were carried out on a personal computer: an Intel (R) Xeon (R) CPU E3-125 v6 @3.30 GHz with 16 GB of memory. The results reveal that the online estimation of the statistical moments of the output response of quantities of interest by the POD-BSBEM can be made  very quickly. Indeed, the associated speed-up factors of the present time-dependent test case, evaluated with respect to the time required by the crude Monte Carlo method, are about $885$ for POD-BSBEM, $164$ for Full-PCE, and $39$ for POD-ANN, thus demonstrating the efficiency of the proposed approach compared to other techniques.\\

\begin{table}[h!]
\caption{Comparison of the time costs of the Full-PCE ($p=6$, $n_{p}=2$), POD-ANN ($H_{1,2,3}=50$), POD-BSBEM ($p=2$, $nx=10$) and the MC (with $N_{s}=10^{6}$ realizations), with $\epsilon_{s}=\epsilon_{t}=10^{-10}$ ($L=76$) and $Re_{\mu=800,\,\sigma=200}\in\mathcal{U}\left[454,\;1\,146\right]$. The unit of the time cost is second ($s$).}
\centering
\begin{tabular*}{0.85\textwidth}{@{\extracolsep{\fill}}l r  rrr }
\hline        
           &POD-ANN& POD-BSBEM& Full-PCE& MC\\
\hline        
 Offline& $4\,513.07$& $5.73$& -& -\\
Online& $98.93$& $4.32$& $23.29$& $3\,824.83$\\
\hline 
\end{tabular*}
\label{tab:Tab_comp_cpu_Burg}
\end{table}

\subsection{Application to a hypothetical dam-break with real terrain data} \label{dam_break}
\begin{figure}[ht!]
  \centering
    \begin{subfigure}[b]{0.49\textwidth}
      \centering
        \includegraphics[width=\textwidth]{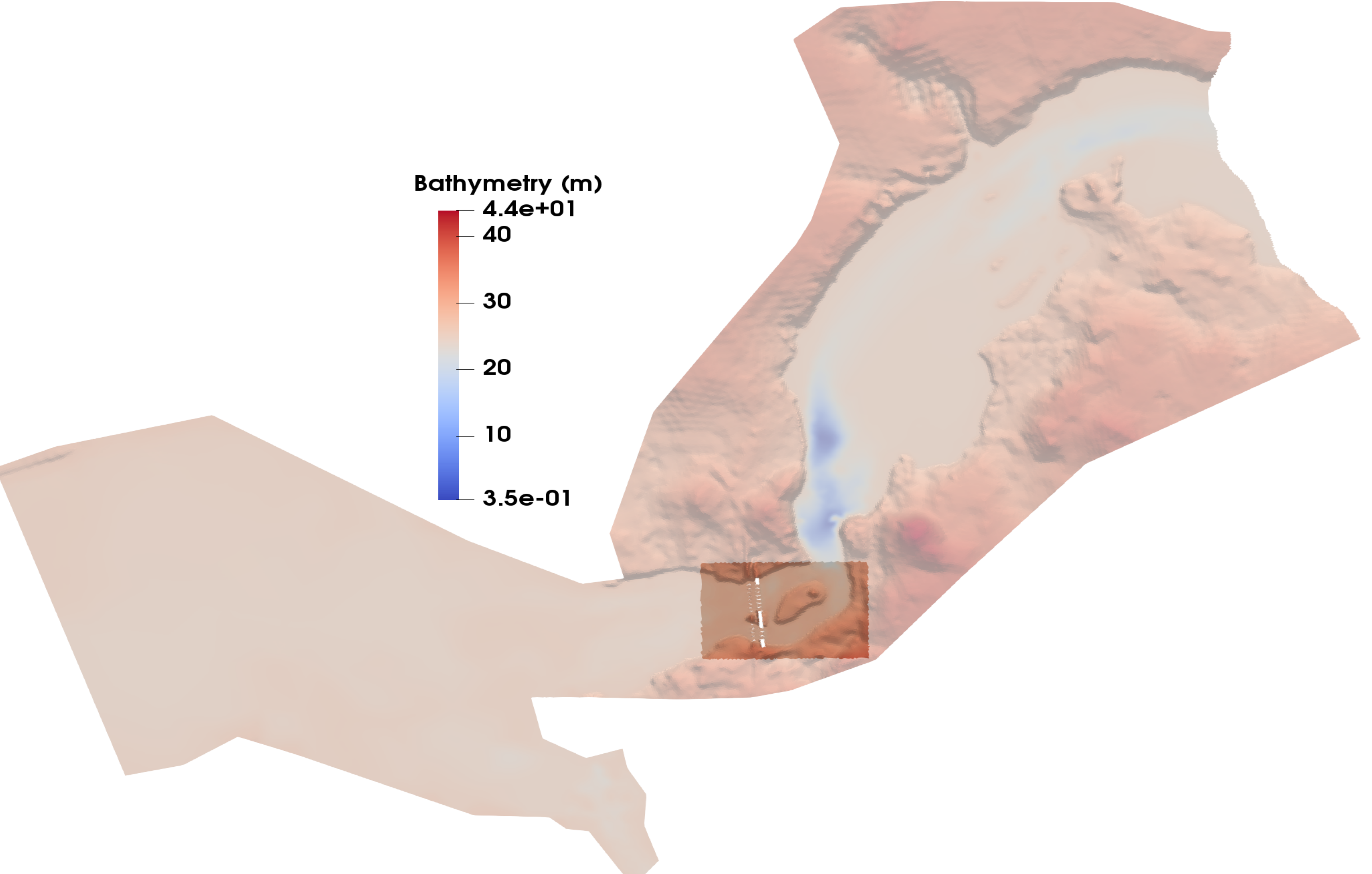}
         \caption{Reach of the Mille-Iles river}
         \label{fig:M_Iles_reach_bathy}
    \end{subfigure}  
  \hfill
    \begin{subfigure}[b]{0.49\textwidth}
      \centering
        \includegraphics[width=\textwidth]{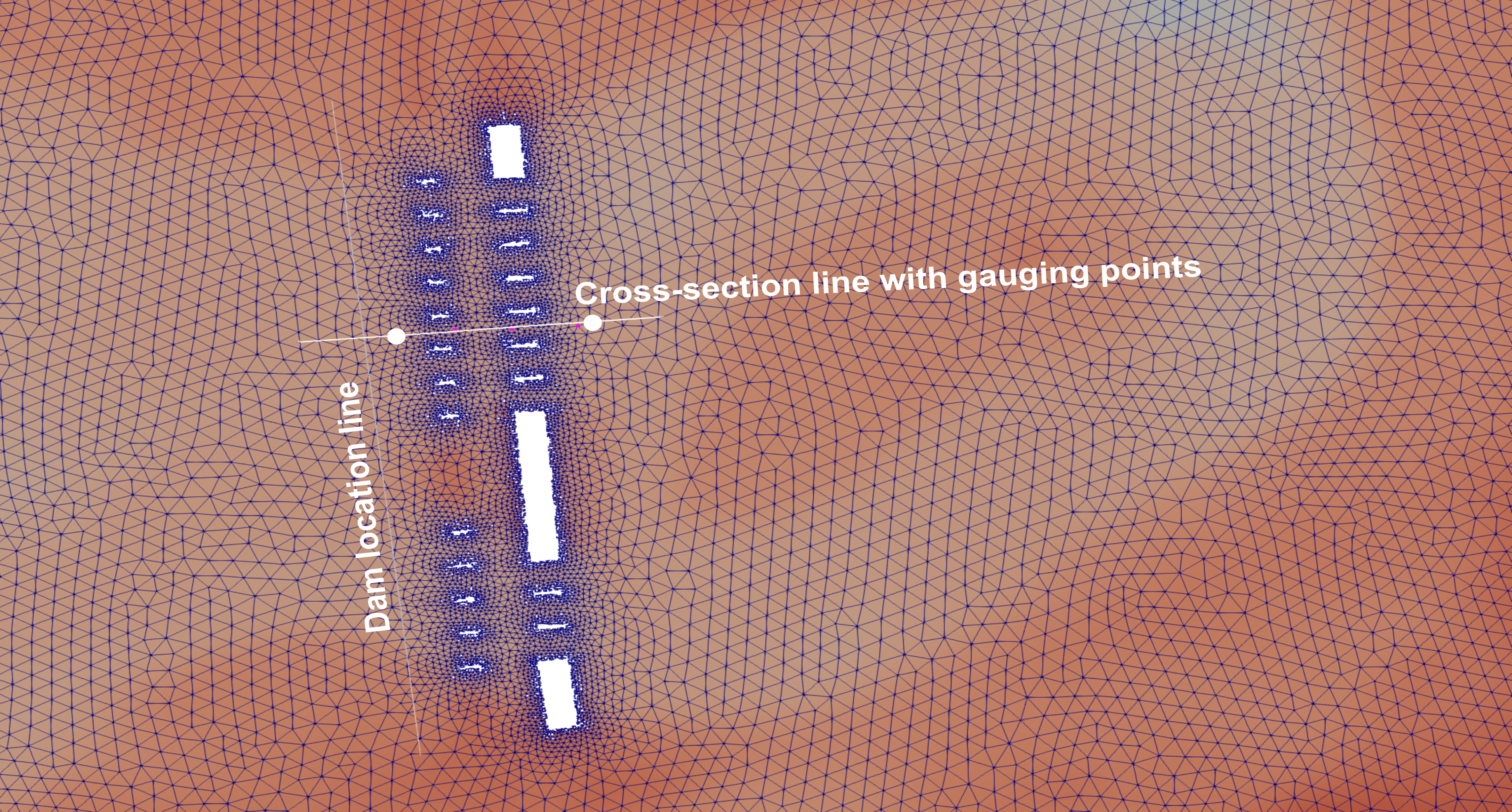}
         \caption{Zoom of the studied sub-domain}
         \label{fig:M_Iles_bathy_zoom}
     \end{subfigure} 
   \caption{Sketch of the reach of the Mille Iles River. (a): Bathymetry of the reach of the Mille Iles River; (b): Zoom of the sub-domain of study; the horizontal line ($x_{line}$) indicates the cross-section over which results are presented, whereas the vertical line delineates the position of the dam from which the breaching process was initiated. The two points (from point 1 near the dam location to point 2 downstream) represent the gauging positions where results are depicted as a function of time (for interpretation of the references to color in this figure, please refer to the web version of this article).}
   \label{fig:M_Iles_domain}
\end{figure}

In addition to the two benchmark test cases presented above for steady and time-dependent problems, the aforementioned modeling approaches were applied to uncertainty propagation analysis regarding a hypothetical dam-break with real terrain data within a reach of the Mille Iles river, which takes its source from the Lake of Two Mountains and discharging downstream into the Des Prairies river near its confluence with the Saint-Laurent river (province of Qu\'{e}bec, Canada). The data related to the bathymetry and the Manning roughness coefficient of the studied reach were provided by the Communeaut\'{e} M\'{e}tropolitaine de Montr\'{e}al (CMM). The hydraulic behavior of the dam-break flow is described by a deterministic in-house finite volume CuteFlow solver for shallow water equations in its multi-GPU version \citep{zokagoa2012pod,delmas2020multi}. The discretized free surface equations were solved through an unstructured mesh composed of $243\,161$ nodes connected to $481\,930$ triangular elements, representing the whole computational domain, half of which is shown in Fig.\ref{fig:M_Iles_reach_bathy}. A part of the computational domain, delineated by a sub-domain surrounding the location of the dam as shown in Fig.\ref{fig:M_Iles_bathy_zoom}, composed of $N_{e}=10\,200$ nodes and $16\,763$ elements was considered as a domain of study from which snapshots of the quantities of interest were extracted to construct the reduced-order models.\\

This test case describes a fictitious sudden rupture of a dam located within the domain of the study area, as shown by the line indicating its position in Fig.\ref{fig:M_Iles_bathy_zoom}. The initial conditions are defined by unequal water levels on both sides of the dam. The downstream part is considered as dry ($\eta_{ds}=b$, where $b$ denotes the bathymetry of the domain), while the initial upstream water surface elevation ($\eta_{up}$) is considered a uniformly-sampled uncertain input parameter with a given sample size ($N_{s}$) within its plausible variability range ($\eta_{up}\in\mathcal{U}\left[29,\;32\right]\,m$). A single value is selected from the generated random space to run the CuteFlow deterministic solver for a simulation time of $50\,s$. The obtained high-fidelity solutions of quantities of interest over the $N_{e}$ computational nodes for each time-step, defined as a time-parameter dependent snapshot, were recorded to construct the snapshot matrix. The number of time points considered is $N_{t}=100$, which is much lower than the value from the adaptive time-steps used in the deterministic numerical solver.\\

The results of the dam-break test case are presented mainly in terms of the mean and standard deviation profiles of the water surface level and the confidence interval distribution contours of the flooding lines accompanying the propagation of the shock wave over the whole computational domain. The solutions obtained from the LHS method with $N_{s}=2000$ are considered as reference solutions with which to assess the accuracy of the proposed approaches. The comparisons of the obtained statistical moments are performed over a cross-section line in the domain of study at different simulation times, as well as at different gauging points as a function of time, as shown in Fig.\ref{fig:M_Iles_bathy_zoom}. The construction of the reduced-order models (POD-BSBEM and POD-ANN) was performed by a series of computational tests to select the most dominate modes constituting the POD basis. Thus, a dual-POD was performed for both reduced models, with $\epsilon_{t}=10^{-7}$ and $\epsilon_{s}=10^{-8}$ for POD-BSBEM, and $\epsilon_{t}=\epsilon_{s}=10^{-8}$ for POD-ANN, producing $L=824$ modes (with $N_{s}=60$) and $L=3387$ modes (with $N_{s}=300$), respectively. The Artificial Neural Network used in the POD-ANN model is composed of feedforward neural networks which each have three hidden layers with $50$ neurons in each ($H_{1,2,3}=50$). As for the former time-dependent test case, the constructed networks were trained using the scaled conjugate gradient backpropagation algorithm within the MATLAB neural network training toolbox (\textit{nntraintool}) with its available default options. The data set used in the learning process was obtained from the $N_{s}=300$ high fidelity solutions of the deterministic model that were used for the snapshot matrix.  These were divided into training, validation and testing subsets using the MATLAB \textit{dividerand} function with its default parameters.\\

\begin{figure}[ht!]
  \centering
    \begin{subfigure}[b]{0.49\textwidth}
      \centering
        \includegraphics[width=\textwidth]{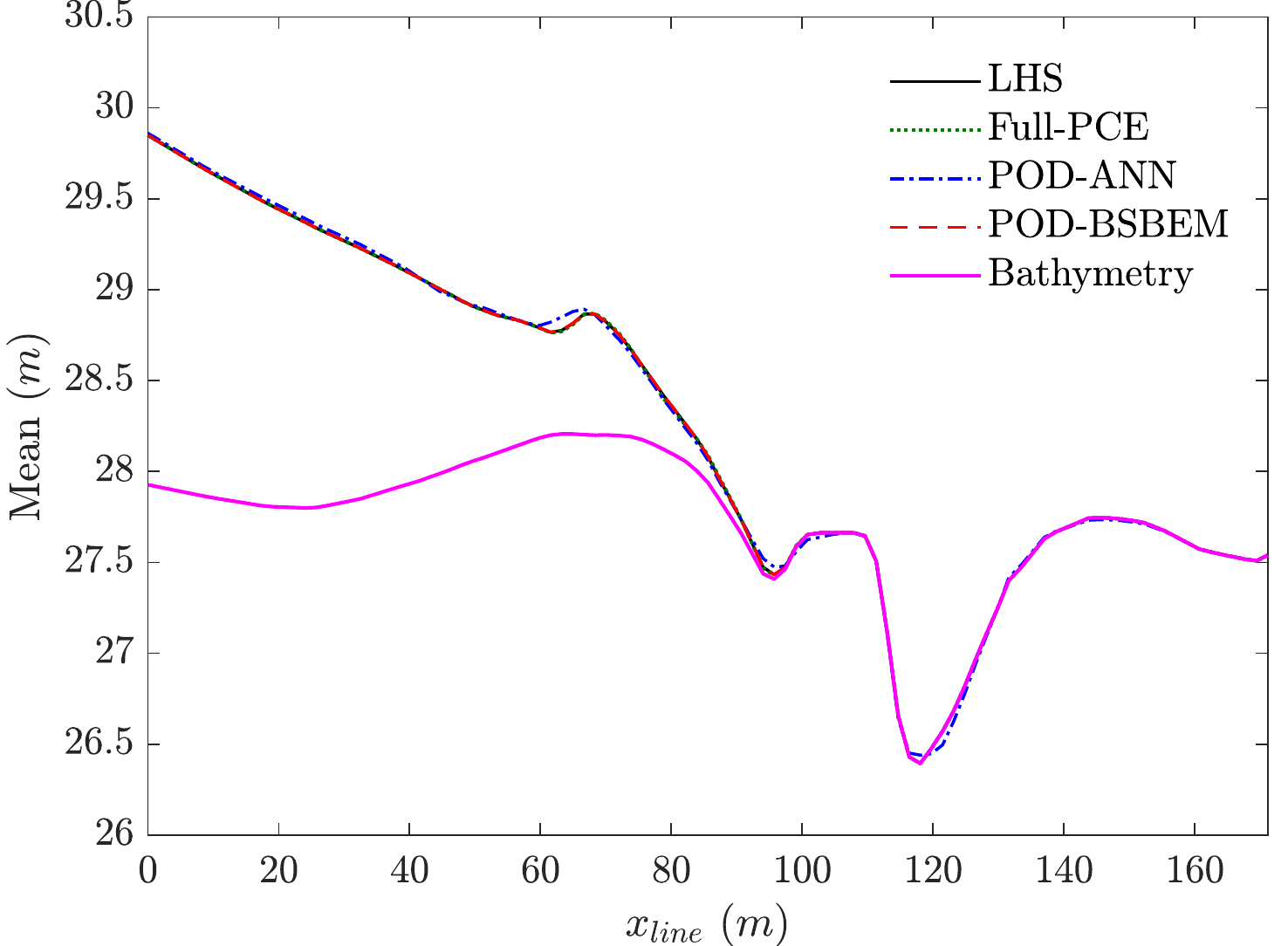}
         \caption{$t\approx\,10\;s$}
         \label{fig:M_Iles_mean_line_It_21}
    \end{subfigure} 
    \hfill 
    \begin{subfigure}[b]{0.49\textwidth}
      \centering
        \includegraphics[width=\textwidth]{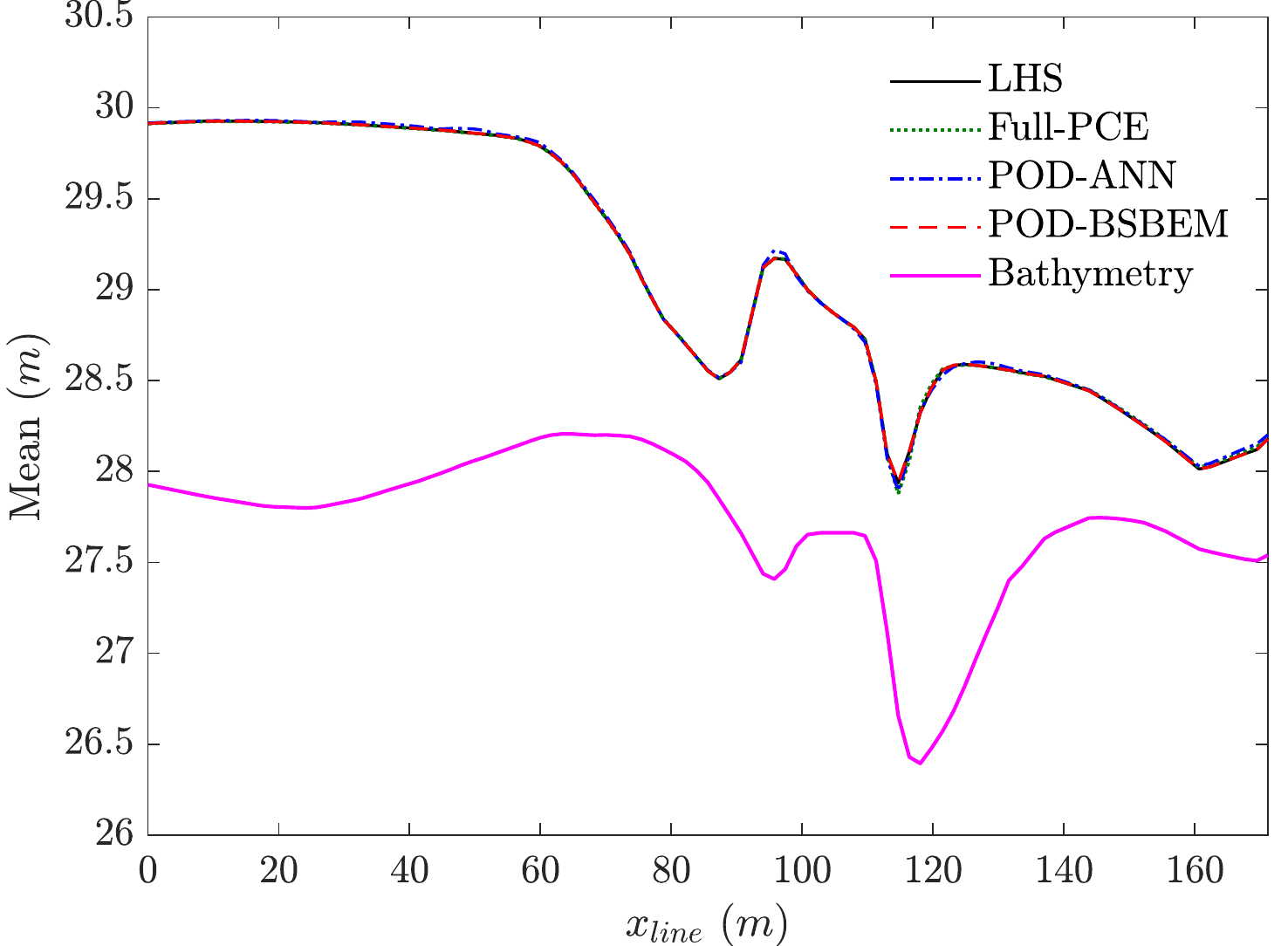}
         \caption{$t\approx\,45\;s$}
         \label{fig:M_Iles_mean_line_It_91}
     \end{subfigure}  
   \caption{The mean profiles of the water levels over the cross-section line at different simulation times, obtained with the Full-PCE ($p=6,\;n_{p}=2$), POD-ANN ($H_{1,2,3}=50,\;N_{s}=300$) and POD-BSBEM ($p=2,\;nx=20$), compared to those of the LHS reference solution ($N_{s}=2\,000 $). (a): $t\approx\,10\;s$, (b): $t\approx\,45\;s$}
   \label{fig:M_Iles_mean_line}
\end{figure}

\begin{figure}[ht!]
  \centering
    \begin{subfigure}[b]{0.49\textwidth}
      \centering
        \includegraphics[width=\textwidth]{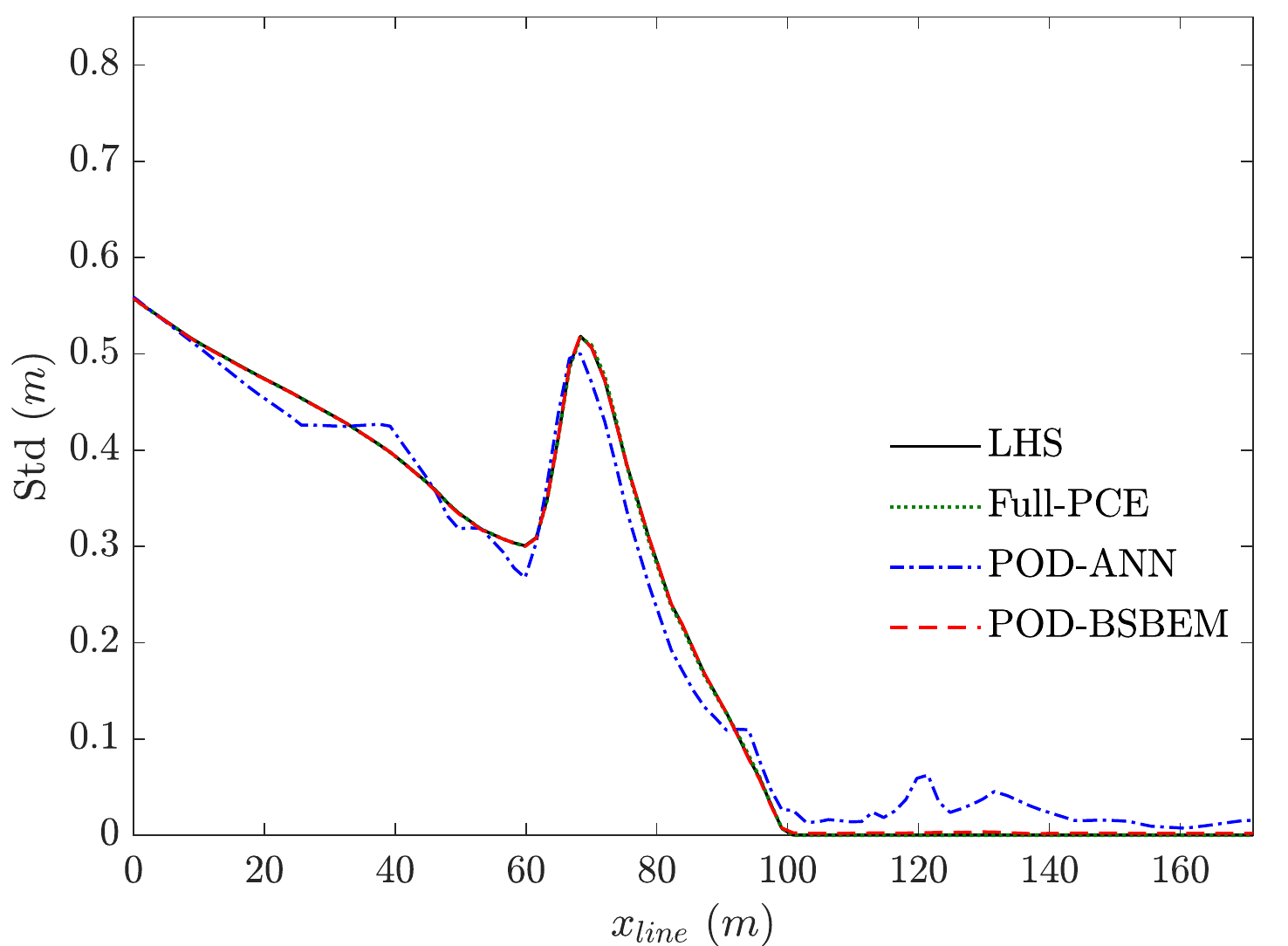}
         \caption{$t\approx\,10\;s$}
         \label{fig:M_Iles_std_line_It_21}
    \end{subfigure}  
  \hfill
    \begin{subfigure}[b]{0.49\textwidth}
      \centering
        \includegraphics[width=\textwidth]{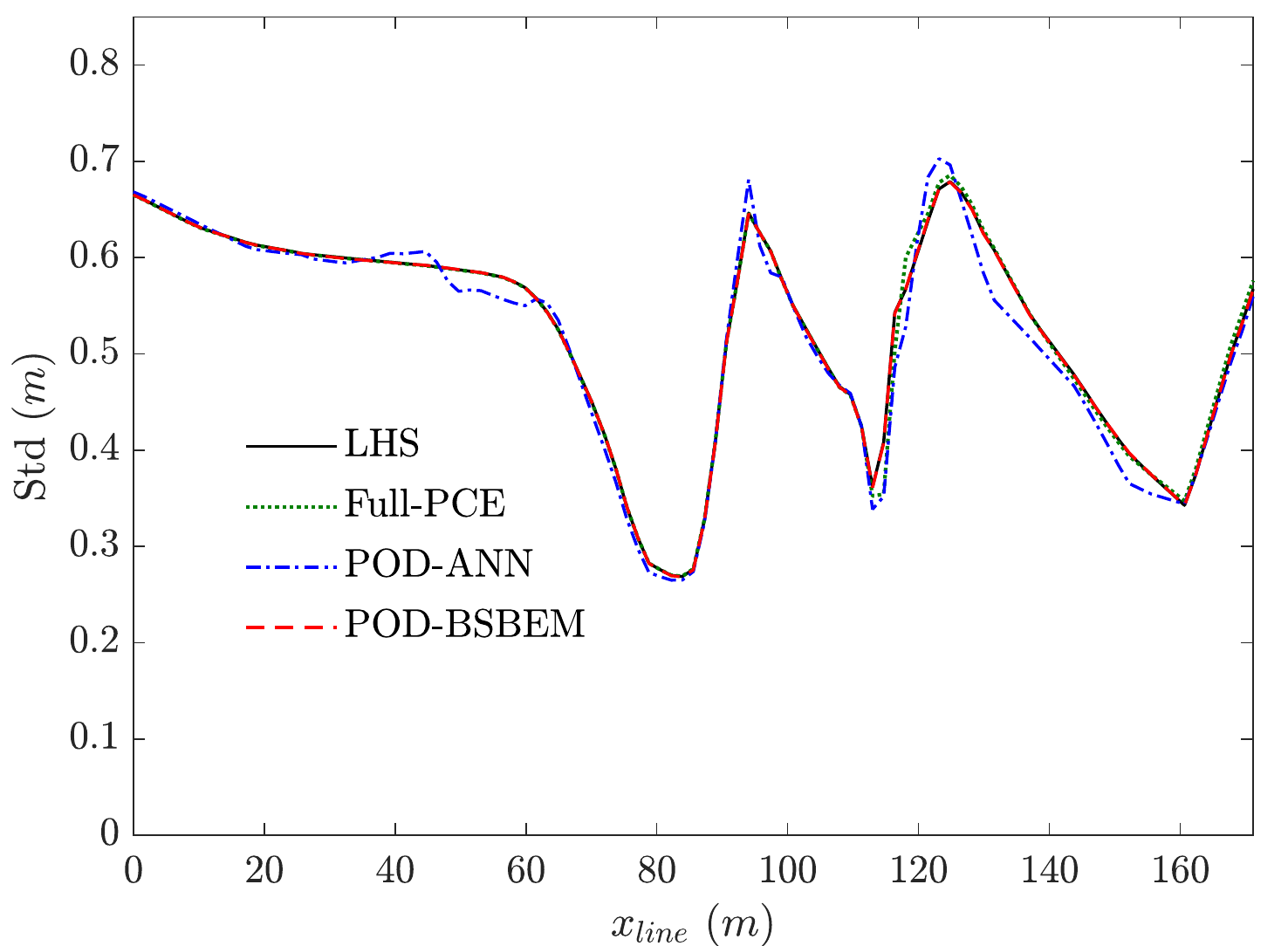}
         \caption{$t\approx\,45\;s$}
         \label{fig:M_Iles_std_line_It_91}
     \end{subfigure}  
   \caption{The standard deviation profiles of the water levels over the cross-section line at different simulation times, obtained with the Full-PCE ($p=6,\;n_{p}=2$), POD-ANN ($H_{1,2,3}=50,\;N_{s}=300$) and POD-BSBEM ($p=2,\;nx=20$), compared to those of the LHS reference solution ($N_{s}=2\,000$). (a): $t\approx\,10\;s$, (b): $t\approx\,45\;s$}
   \label{fig:M_Iles_std_line}
\end{figure}

Fig.\ref{fig:M_Iles_mean_line} shows the mean profiles of the water surface level over the cross-section line $x_{line}$ obtained with the Full-PCE ($p=6,\;n_{p}=2$), POD-ANN ($H_{1,2,3}=50,\;N_{s}=300$) and POD-BSBEM ($p=2,\;nx=20$) at different simulation times, compared with those from the LHS reference solution ($N_{s}=2\,000$ realizations). The bathymetry of the terrain is also represented in the figure in order to show the propagation of the flooding wave over the selected cross-section as the simulation time evolves. In general, the mean distribution profiles of all methods appear to be in good concordance with those from the reference LHS solution at the represented simulation times. However, the comparison of the standard deviation profiles, depicted in Fig.\ref{fig:M_Iles_std_line}, shows clear deviations between the POD-ANN and the LHS results, unlike those obtained by the POD-BSBEM method where an excellent superposition is observed with the LHS reference results for all simulation times. In addition to the relatively significant learning time required by the POD-ANN technique in its offline stage, a total of $30\,000$ ($N_{s}=300\times N_{t}=100$) high fidelity solutions were collected for computing the POD basis, versus  $6\,000$ ($60\times100$) for the proposed POD-BSBEM, which represents a reduction in the number of deterministic solver calls and thus in the computational effort required for the construction of the reduced basis space. Another point of comparison concerns the reduced-order aspect that characterizes the POD-BSBEM compared to the Full-order polynomial chaos expansion (Full-PCE). Indeed, while both techniques show a good approximation of the statistical moments, the Full-PCE is used to build a stochastic expansion that approximates the output response in each node over the whole computational space-time domain ($N_{e}\times N_{t}$), and the proposed POD-BSBEM only needs the construction of $\sum_{\ell}^{L}\mathcal{K}_{\ell}$ expansions over the $L$ modes of the reduced basis. \\

\begin{figure}[ht!]
  \centering
    \begin{subfigure}[b]{0.49\textwidth}
      \centering
        \includegraphics[width=\textwidth]{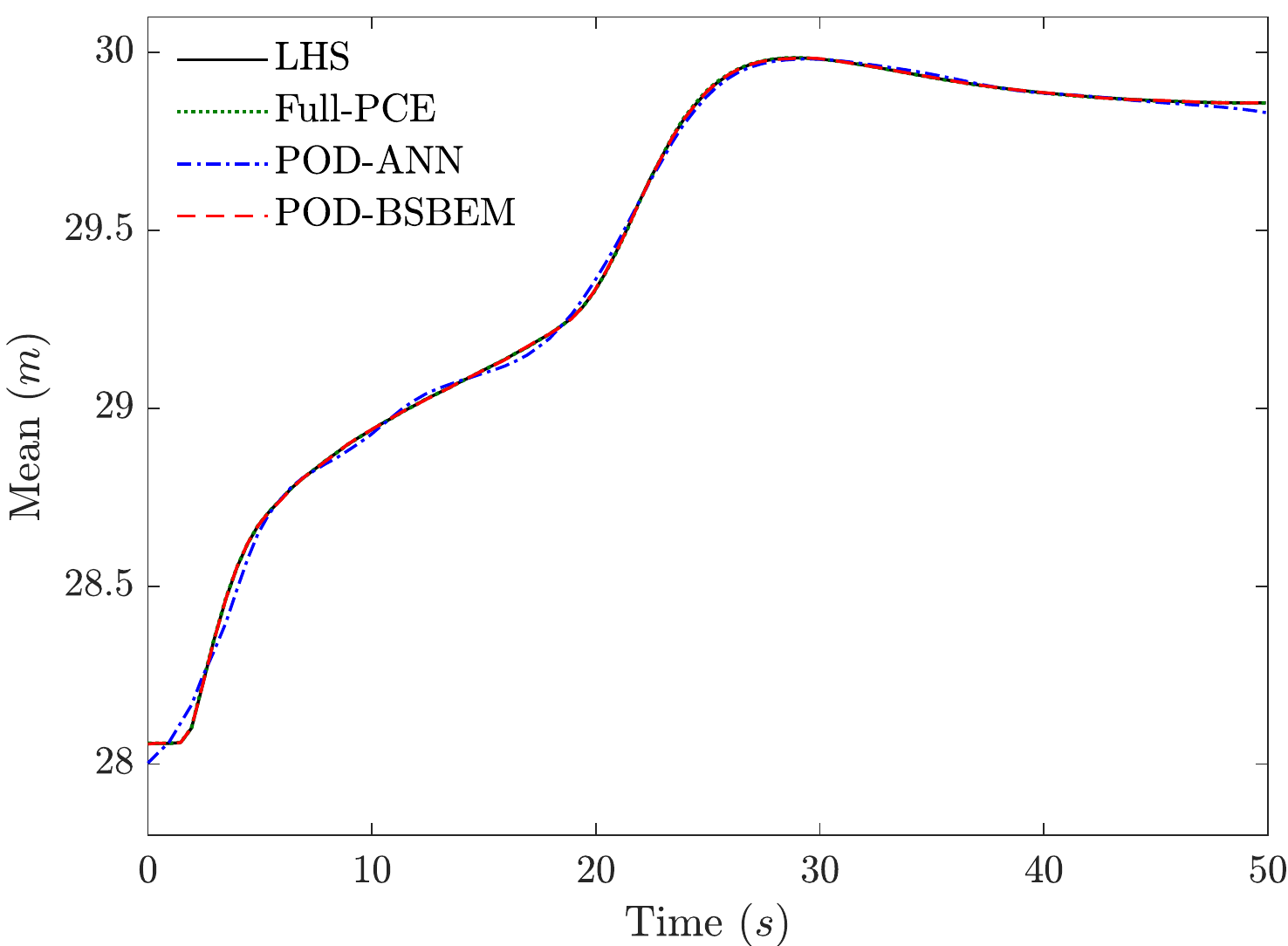}
         \caption{Point 1}
         \label{fig:M_Iles_mean_time_p1}
    \end{subfigure}  
     \hfill
    \begin{subfigure}[b]{0.49\textwidth}
      \centering
        \includegraphics[width=\textwidth]{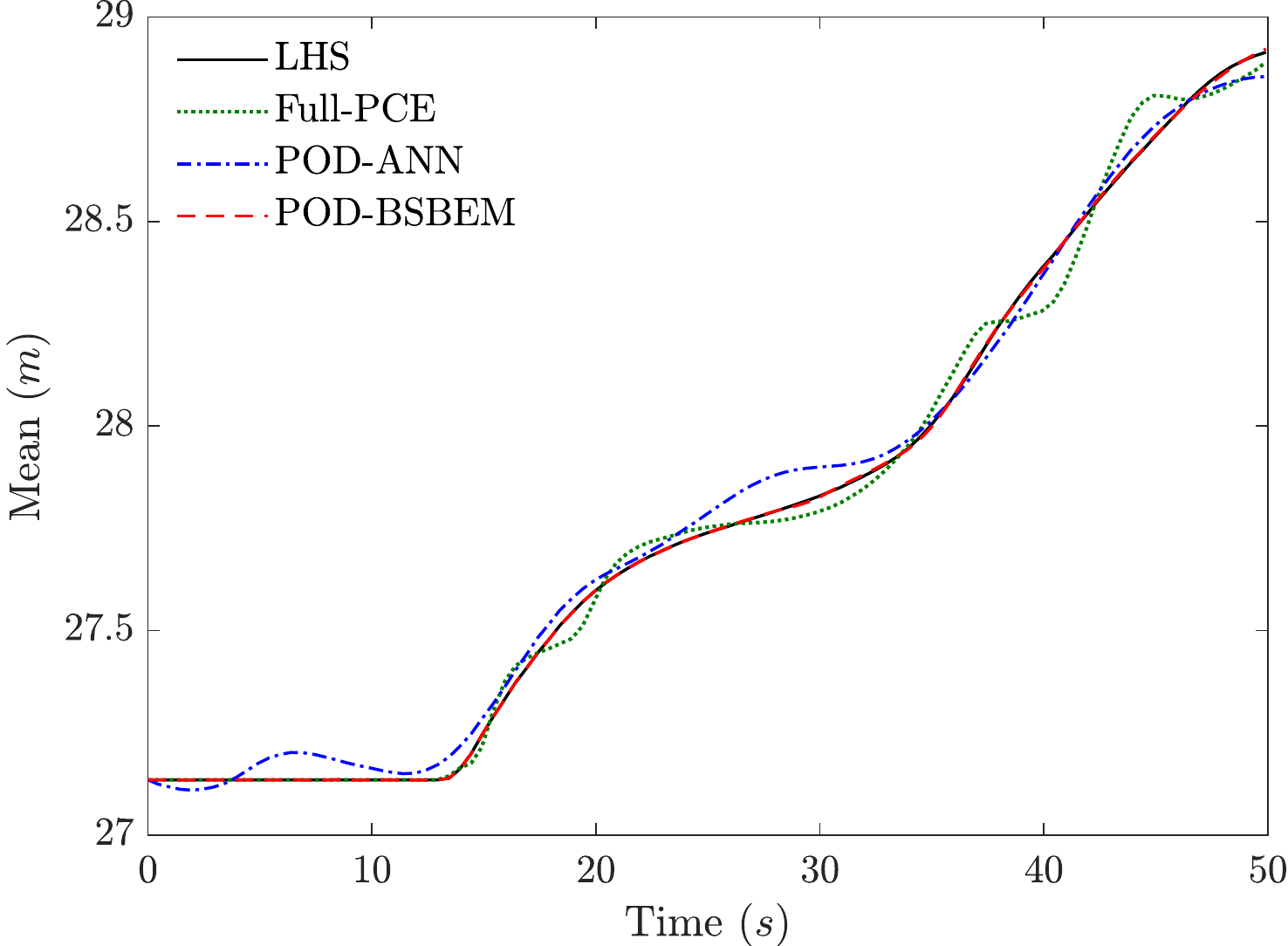}
         \caption{Point 2}
         \label{fig:M_Iles_mean_time_p4}
     \end{subfigure}  
   \caption{The mean water level profiles  at gauging points as a function of time, obtained with the Full-PCE ($p=6,\;n_{p}=2$), POD-ANN ($H_{1,2,3}=50,\;N_{s}=300$) and the POD-BSBEM ($p=2,\;nx=20$), compared to those from the LHS reference solution ($N_{s}=2\,000$). (a): Point 1, (b): Point 2.}
   \label{fig:M_Iles_mean_time}
\end{figure}

\begin{figure}[ht!]
  \centering
    \begin{subfigure}[b]{0.49\textwidth}
      \centering
        \includegraphics[width=\textwidth]{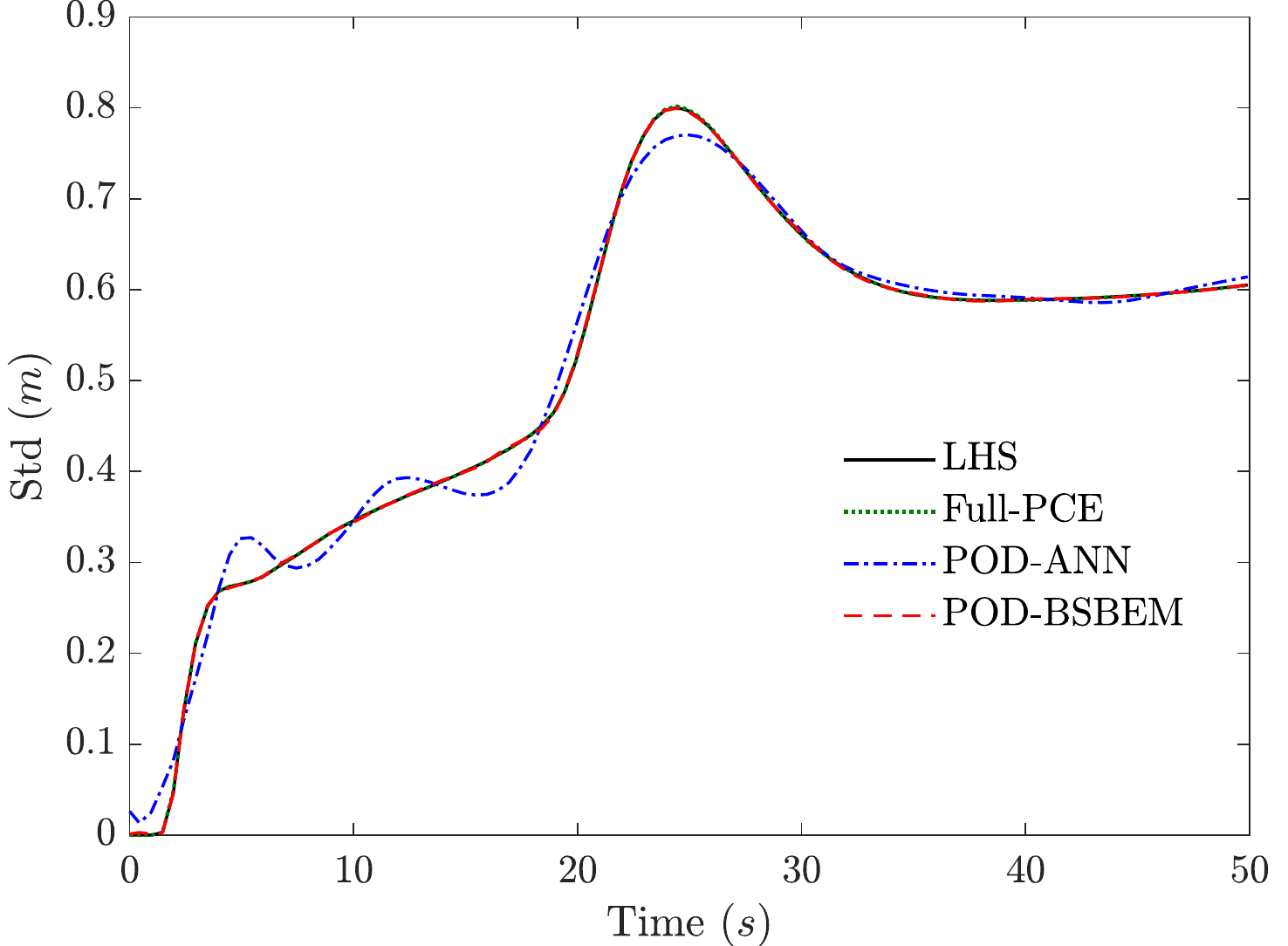}
         \caption{Point 1}
         \label{fig:M_Iles_std_time_p1}
    \end{subfigure}  
  \hfill
    \begin{subfigure}[b]{0.49\textwidth}
      \centering
        \includegraphics[width=\textwidth]{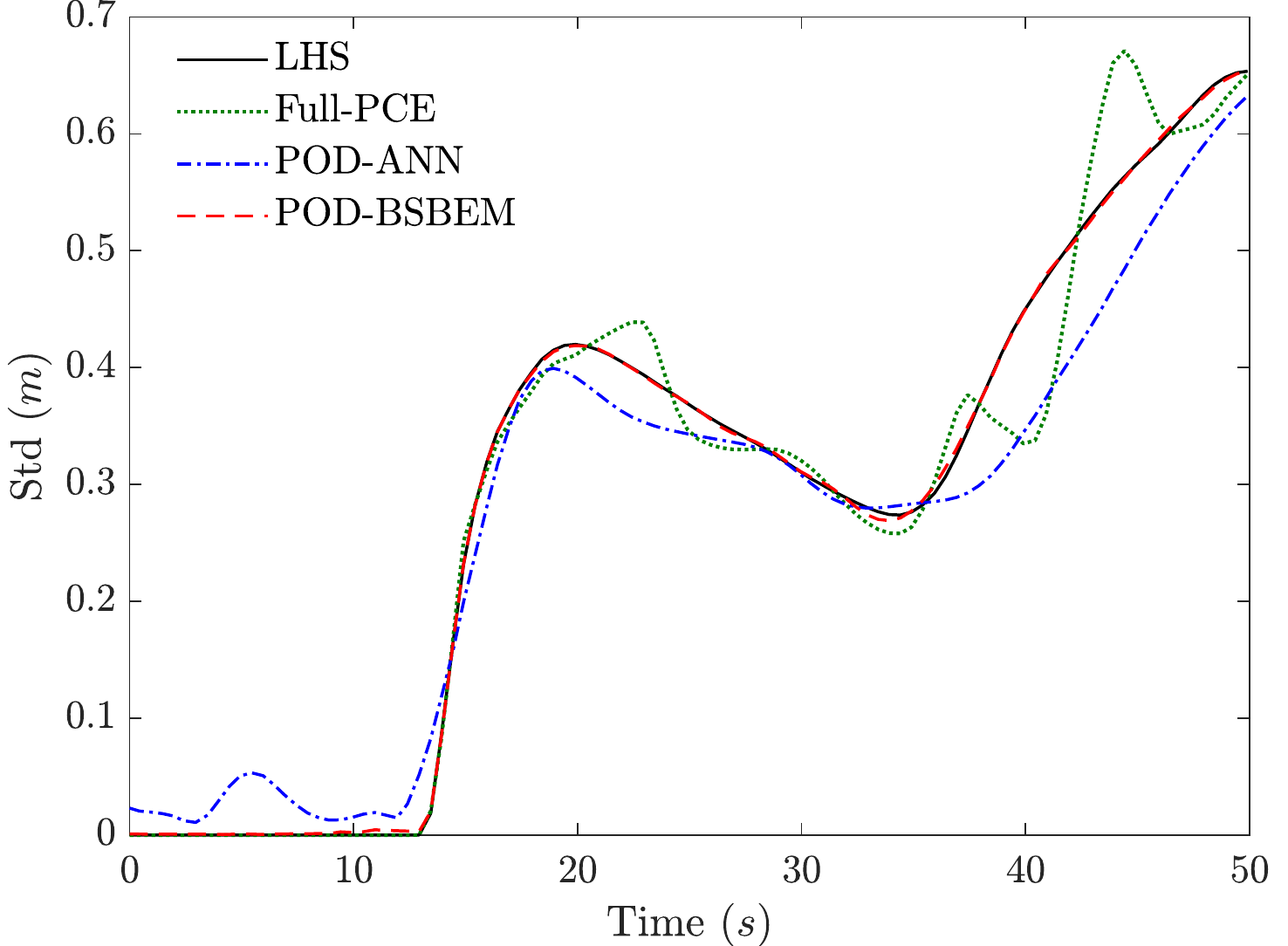}
         \caption{Point 2}
         \label{fig:M_Iles_std_time_p4}
     \end{subfigure}  
   \caption{The standard deviation profiles of the water levels at gauging points as a function of time, obtained with Full-PCE ($p=6,\;n_{p}=2$), POD-ANN ($H_{1,2,3}=50,\;N_{s}=300$) and POD-BSBEM ($p=2,\;nx=20$), compared to those from the LHS reference solution ($N_{s}=2\,000$). (a): Point 1, (b): Point 2.}
   \label{fig:M_Iles_std_time}
\end{figure}

In Fig.\ref{fig:M_Iles_mean_time}, the temporal evolution of the mean profiles of the water level is depicted at two gauging points. These points are selected from the sub-domain of study near the cross-section line to show the abilities of the proposed POD-BSBEM to accurately approximate the statistical moments of the output responses with hyperbolic behavior as the simulation time evolves. The results are compared to those of the POD-ANN and Full-PCE techniques, as well as to the results of the reference LHS solutions. These plots reveal that both the POD-ANN and the Full-PCE techniques present some discrepancies in their profiles; a slight oscillatory behavior can be observed, particularly at gauging point $P2$ as shown in Fig.\ref{fig:M_Iles_mean_time_p4}. In contrast, the POD-BSBEM gives particularly good predictions in the mean-field reconstruction of the water level, where an excellent agreement can be observed between the proposed POD-BSBEM and LHS profiles at all gauging points over the whole temporal domain. This trend is confirmed by Fig.\ref{fig:M_Iles_std_time} where the temporal evolution of the standard deviation of the water level is represented, and which seems to be much more sensitive than the mean profile. Indeed, it is clear that both the POD-ANN and the Full-PCE techniques present oscillations, and that their predictions of the standard deviation are less satisfactory, revealing non-physical behavior at some gauging points Fig.\ref{fig:M_Iles_std_time_p4}. The PCE basis functions are known to be less accurate than other approaches at estimating the strong hyperbolic behavior that characterizes the propagation of a flooding wave. On the other hand, an excellent agreement between the POD-BSBEM and LHS results can be observed for all gauging points. This gives a clear idea about the abilities of the POD-BSBEM to capture the temporal dynamics of the output quantities of interest and also indicates the potential that the piecewise B-splines basis functions offer to alleviate oscillatory behavior due to the mathematical features they procure and the multi-elements' aspect of the proposed approach.\\

\begin{figure}[ht!]
  \centering
    \begin{subfigure}[b]{0.49\textwidth}
      \centering
        \includegraphics[width=\textwidth]{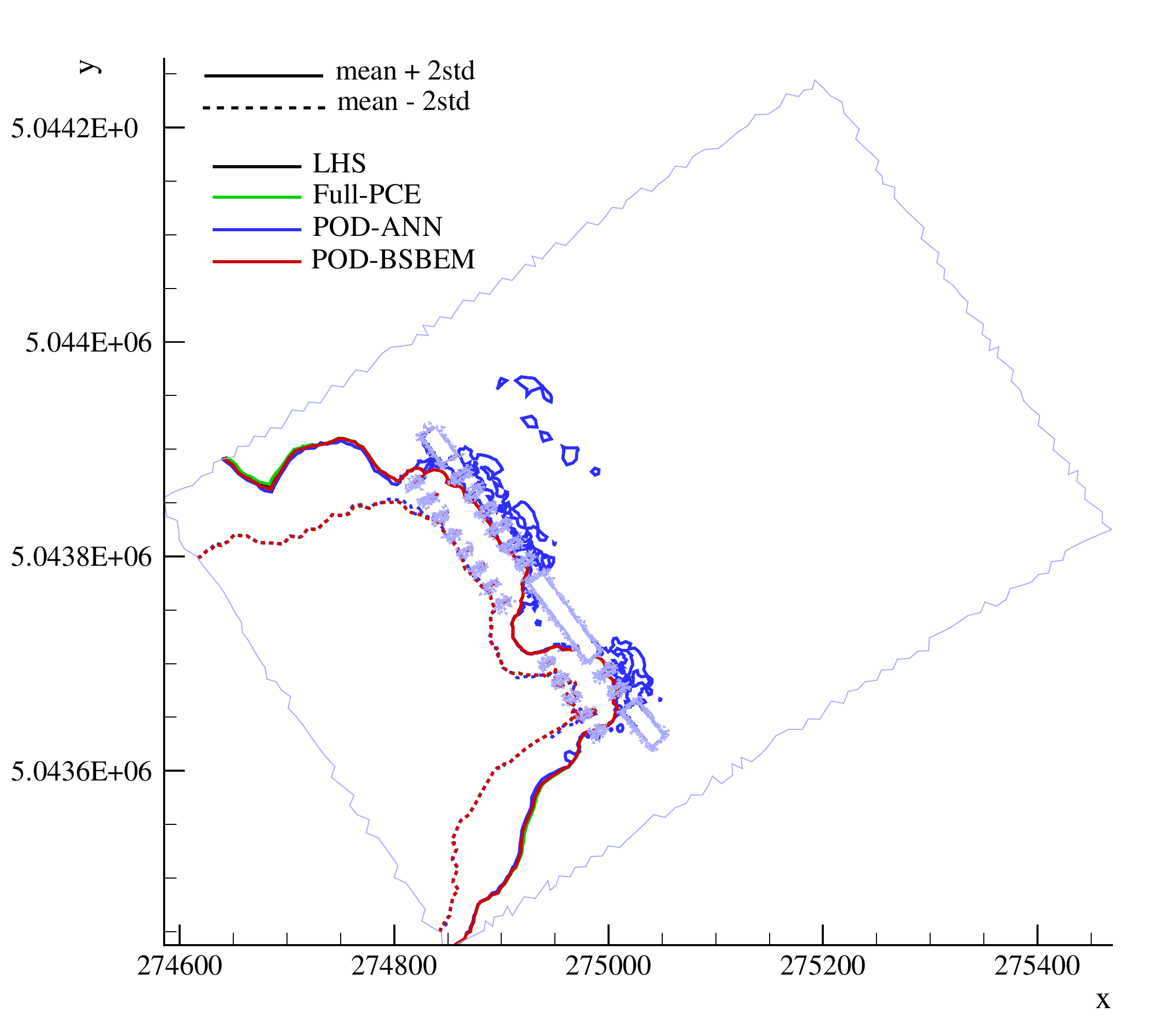}
         \caption{$t\approx\,10\;s$}
         \label{fig:M_Iles_inundation_21}
    \end{subfigure}  
  \hfill
    \begin{subfigure}[b]{0.49\textwidth}
      \centering
        \includegraphics[width=\textwidth]{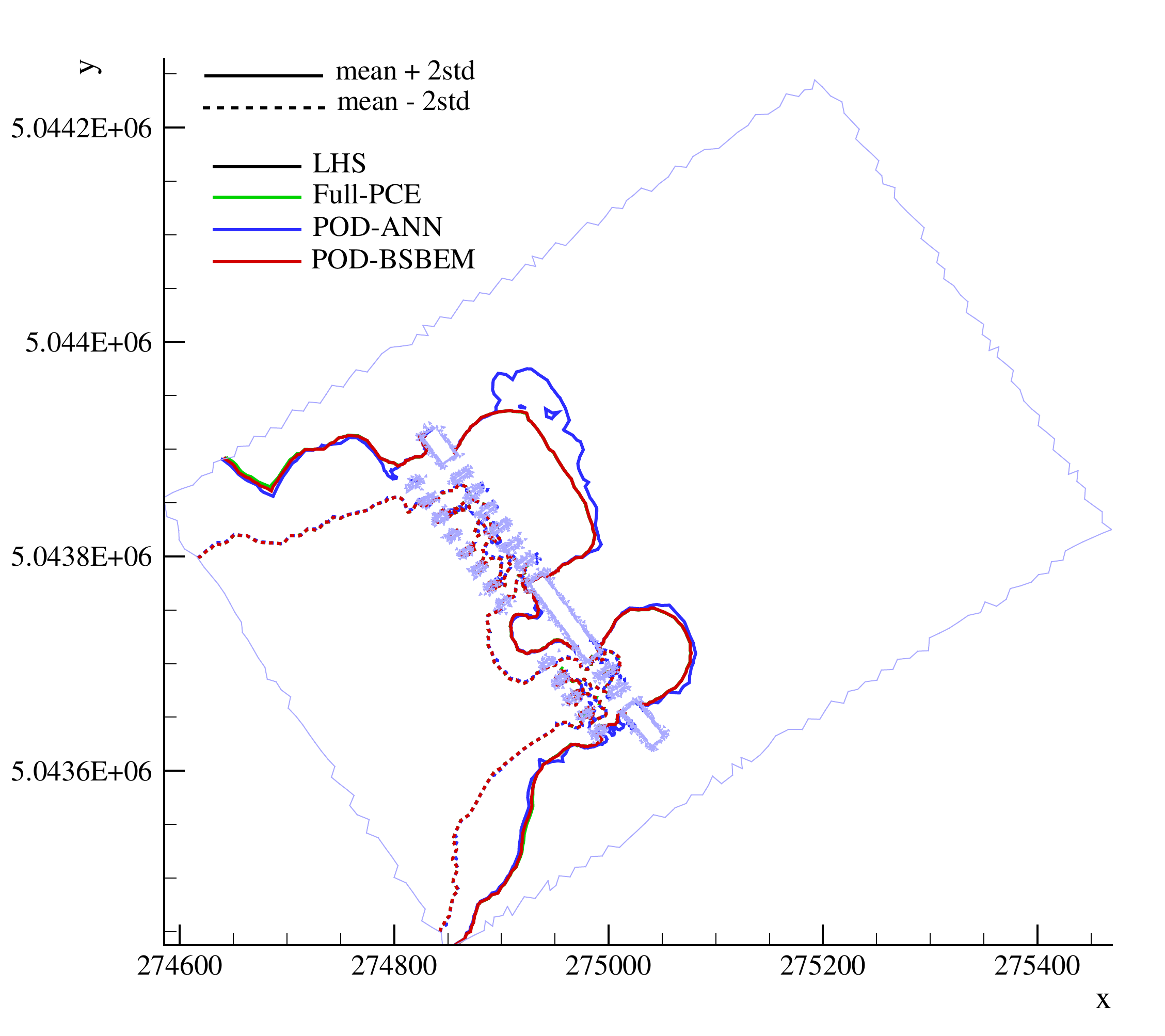}
         \caption{$t\approx\,25\;s$}
         \label{fig:M_Iles_inundation_51}
     \end{subfigure} 
     \hfill
    \begin{subfigure}[b]{0.49\textwidth}
      \centering
        \includegraphics[width=\textwidth]{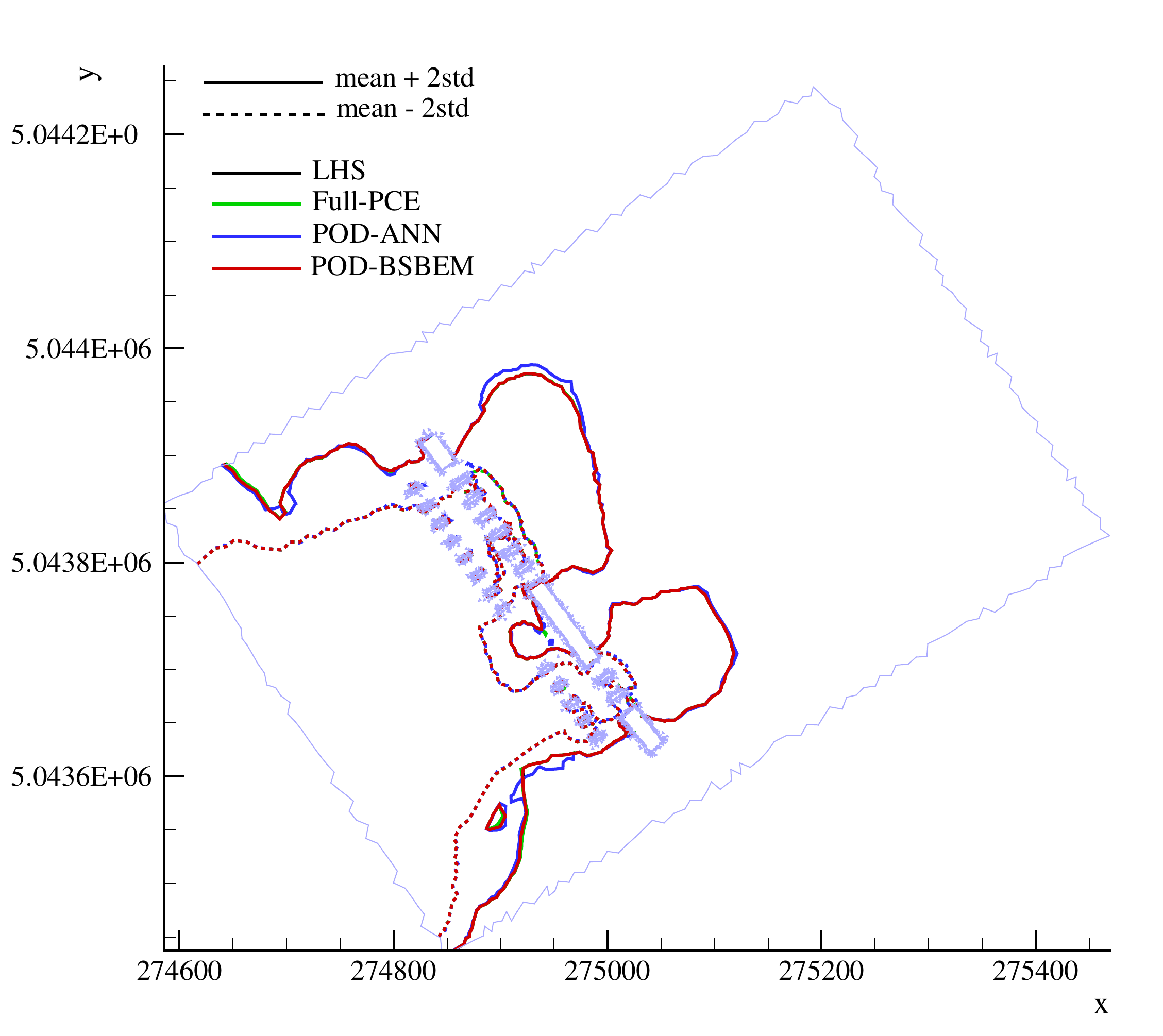}
         \caption{$t\approx\,37\;s$}
         \label{fig:M_Iles_inundation_75}
     \end{subfigure} 
     \hfill
    \begin{subfigure}[b]{0.49\textwidth}
      \centering
        \includegraphics[width=\textwidth]{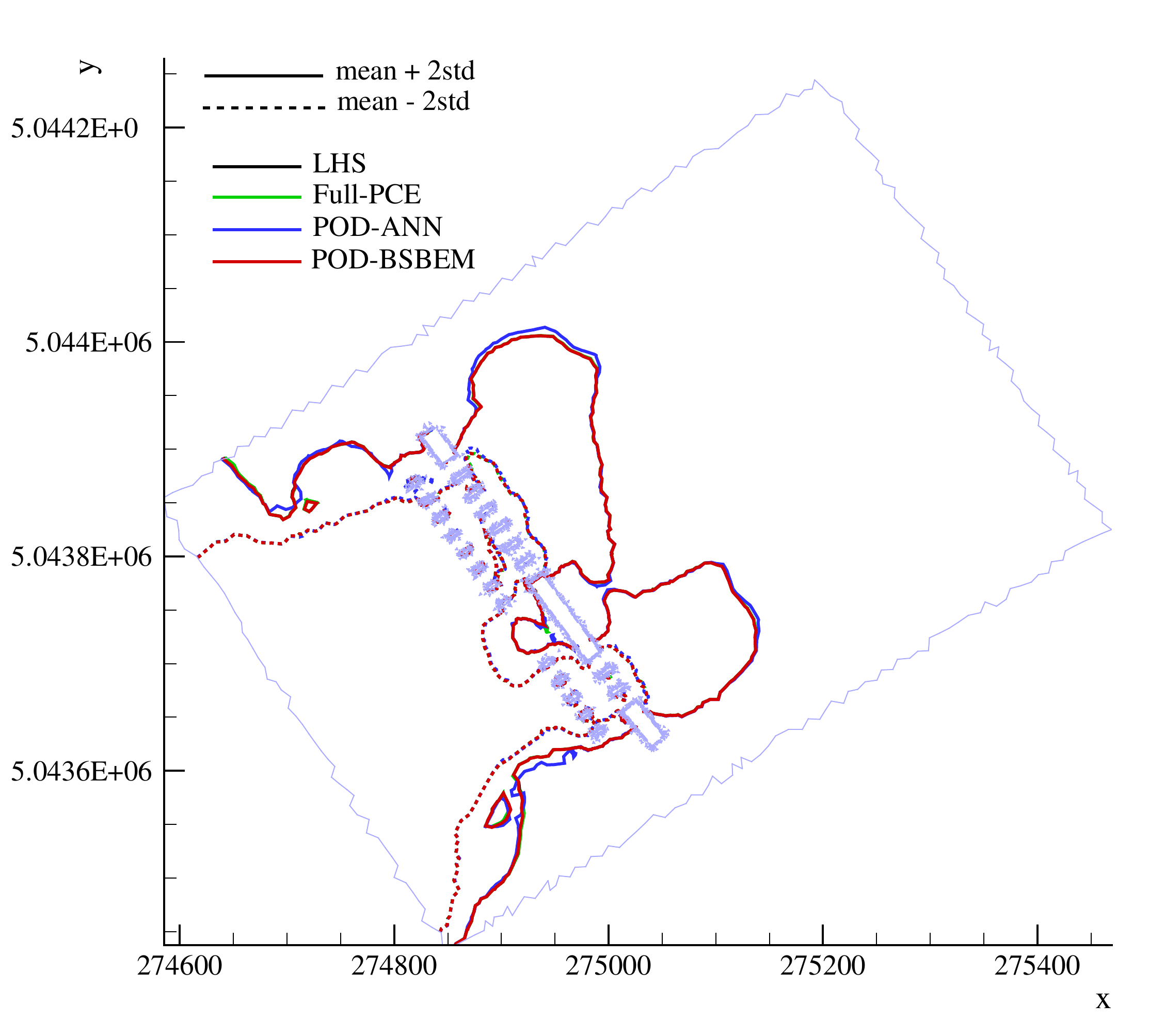}
         \caption{$t\approx\,45\;s$}
         \label{fig:M_Iles_inundation_91}
     \end{subfigure}  
   \caption{The $95\,\%$ confidence interval contours of flooding lines ($h=0.05\,m$) at different simulation times, obtained with the Full-PCE ($p=6,\;n_{p}=2$), POD-ANN ($H_{1,2,3}=50,\;N_{s}=300$) and the POD-BSBEM ($p=2,\;nx=20$), compared to those from the LHS reference solution ($N_{s}=2\,000$). The dashed and solid lines delineate the lower and upper bounds, respectively, of the $95\,\%$ CI of the flooding line. (a): $t\approx\,10\;s$, (b): $t\approx\,25\;s$, (c): $t\approx\,37\;s$, (d): $t\approx\,45\;s$. (For interpretation of the references to color in this figure, please refer to the web version of this article ).}
   \label{fig:M_Iles_Inundation_CI}
\end{figure}

Fig.\ref{fig:M_Iles_Inundation_CI} illustrates the $95\,\%$ confidence interval ($95\,\%$ CI)-bound contours of the water depth, representing the probabilistic inundation map obtained by the POD-BSBEM and the POD-ANN, compared to the contours of the  Full-PCE and the LHS reference solutions at different simulation times. The CI is quantified in each node of the computational domain by the mean, plus and minus two standard deviations of the water depth. The flooding line is selected by the threshold water depth of $5\;cm$ defining the borderline between the wet and the dry areas. Thus, the projection of the lower and upper bounds of the CI (represented by the dashed and solid lines, respectively) corresponding to the aforementioned threshold of $5\;cm$ onto the topography of the studied sub-domain delineate the variability range surrounding the flooding line. These plots show the temporal evolution of the variability range surrounding the flooding line, which becomes widespread as the simulation time evolves. This evolution traces the propagation of the uncertainty stemming from the initial input water level where the break was initiated. The POD-BSBEM clearly yields accurate predictions of the CI contours of the flooding line for all simulation times. Indeed, the superimposition of the flooding line contours illustrates an excellent agreement with the results from the Full-PCE and the reference LHS solutions (with $2\,000$ realizations). However, obvious deviations can be seen in the predictions from the POD-ANN technique, particularly at times $t=10$ and $25\;s$, as shown in Figs.\ref{fig:M_Iles_inundation_21} and \ref{fig:M_Iles_inundation_51}, where the obtained variability range contours do not match with the contour from the reference LHS solution.\\

To illustrate the computational efficiency of the proposed method, its time cost is compared with the costs of other techniques, as presented in Table \ref{tab:Tab_comp_cpu_M_iles}. The comparison mainly concerns the reduced-order models in their offline and online phases. The high-fidelity solutions of all these methods were obtained separately by running the numerical solver CuteFlow on a $V100$ GPU whose average computation time in the half-domain of study is almost $30\;s$. Therefore, this aspect is not included in the offline phase of the POD-ANN and POD-BSBEM approaches. Instead, it is considered separately by comparing the number of numerical solver calls ($N_{s}$) of each technique, as reported in the third row of Table \ref{tab:Tab_comp_cpu_M_iles}. The POD-ANN implementation was performed on a V100 GPU with $32$ GB of memory, with the MATLAB toolbox utilized to train the Artificial Neural Network with the aforementioned architecture, while the POD-BSBEM and Full-PCE were implemented serially and their computations were performed on a personal computer: an Intel (R) Xeon (R) CPU E3-125 v6 @3.30 GHz with 16 GB of memory. The data in Table \ref{tab:Tab_comp_cpu_M_iles} indicates that the POD-BSBEM presents a much lower time cost for both offline and online stages than the POD-ANN. The number of high-fidelity solutions required to build the surrogate reduced-order models also shows a significant difference between the two techniques: $N_{s}=300$ for the POD-ANN and only $N_{s}=60$ for the POD-BSBEM. 

\begin{table}[h!]
\caption{Comparison of the time cost required by the POD-ANN ($H_{1,2,3}=50,\;\epsilon_{t}=\epsilon_{s}=10^{-8}$), POD-BSBEM ($p=2,\;nx=20,\;\epsilon_{t}=10^{-7},\;\epsilon_{s}=10^{-8}$), Full-PCE ($p=6,\;n_{p}=2$) and the LHS ($N_{s}=2\,000$). The unit of time cost is the second ($s$).}
\centering
\begin{tabular*}{0.88\textwidth}{@{\extracolsep{\fill}}l r  rrr }
\hline        
           &POD-ANN & POD-BSBEM& Full-PCE& LHS\\
\hline        
Offline ($s$)& $5\,839.93$& $198.64$& -& -\\
Online ($s$)& $1\,612.04$& $382.21$& -& -\\
Number of solver calls ($N_{s}$)& $300$& $60$& $14$& $2\,000$\\
Time cost per snapshot ($s$)& $1.612$& $0.350$& $0.144$& $30$\\
\hline 
\end{tabular*}
\label{tab:Tab_comp_cpu_M_iles}
\end{table}
The time cost required for each snapshot reported in the last row of Table \ref{tab:Tab_comp_cpu_M_iles} shows that the use of the constructed surrogate models significantly reduces the computational cost compared to the high-fidelity numerical solver. The POD-BSBEM is more efficient than the POD-ANN, reducing the computation effort while maintaining good accuracy in the prediction of the output response. The results also reveal that the Full-PCE requires a lower number of high-fidelity solutions and therefore less computation cost per snapshot, despite its less accurate predictions with oscillatory behavior. This lower requirement is due to the univariate basis functions of polynomial chaos (one input random variable), which involves a minimal number of expansion terms. Indeed, the total number of expansion terms (coefficients) that have to be computed over the $N_{e}\times N_{t}=10200\times 100$ space-time nodes is only on the order of $M=p+1=7$. However, this number grows significantly once the number of random variables increases, expressed as of $M=\frac{(p+m)!}{p!m!}$, and therefore the time required to evaluate each expansion increases as well (as shown in Table \ref{tab:Tab_comp_cpu_Ack} for three random variables), particularly for very large computational domains with millions of meshing nodes. In contrast, the POD-BSBEM approach only constructs expansions over the $L$-modes of the reduced basis, which offers a significant gain in the computational cost, in addition to procuring accurate predictions for time-dependent problems with strong hyperbolic behavior.


\section{Conclusion}\label{conc}
This paper proposes a non-intrusive reduced-order model-based B-splines B\'{e}zier elements method (POD-BSBEM) for uncertainty propagation analysis for stochastic steady-state and time-dependent problems. The method combines the advantages of both the model reduction technique based on proper orthogonal decomposition (POD) and a stochastic approach, the B-splines B\'{e}zier elements method (BSBEM), which is based on splitting the random parameter domain into subspaces called B\'{e}zier elements. The construction of the reduced basis is performed from a collection of high-fidelity solutions through a two-level POD procedure by compressing the time-trajectory snapshots in the first POD level and then extracting the final reduced basis via a second POD level. Each associated projection coefficient dataset, which is randomly parametrized and time-dependent, is expressed as a regression constituted as a series of time-dependent reduced modes, generated via a third POD level, and their associated unknown stochastic parameter-dependent coefficients. These coefficients are then approximated as a local stochastic expansion using the local B-splines basis functions defined within the local B\'{e}zier element constituting a part of the random parametric domain.\\

Two benchmark test cases, representing the steady-state Ackley function, with three random input  parameters , and a univariate time-dependent viscous Burgers' equation, were used to investigate the performance and accuracy of the proposed method with respect to the non-intrusive reduced and full-order techniques POD-ANN and Full-PCE, respectively. The numerical results presented in this paper demonstrate the abilities of POD-BSBEM to predict the statistics and the estimated probability density function of the output quantities of interest, particularly for the time-dependent viscous Burgers' equation where a smooth profile of the standard deviation is obtained, unlike the POD-ANN approach which shows an oscillatory behavior. The results also highlighted the computational efficiency of the proposed approach, as they reveal a significant reduction in the computational effort required for both offline and online stages compared to the other techniques.\\

The proposed approach was then applied to investigate uncertainty propagation through a numerical model describing the highly non-linear flood wave stemming from a hypothetical dam-break with real terrain data of a reach of the Mile-Iles river located in the province of Qu\'{e}bec. The results reveal that predictions from POD-BSBEM provide smooth profiles of the standard deviation as a function of time at the gauging points, where excellent agreements with the reference LHS profiles were observed, unlike the predictions from POD-ANN and Full-PCE that present a strong oscillatory behavior. Another meaningful finding is the ability of POD-BSBEM to provide accurate probabilistic inundation maps that consider the uncertainty stemming from the inputs for better flood hazard predictions while considerably reducing the computation time. Indeed, the proposed technique presents a much lower computational effort for both offline and online stages, with a significant speed-up ratio compared to the POD-ANN method. Moreover, the introduced method only requires the construction of expansions over a few limited numbers of modes of the reduced basis, unlike the full-order Full-PCE model that needs to construct an expansion over each node of the computational domain, and must do so for each time step, a cumbersome computational effort particularly for high-dimensional stochastic problems, i.e., with a high number of input random parameters. Thus, the proposed non-intrusive reduced-order model-based B-splines B\'{e}zier element method presents a promising tool for uncertainty propagation analysis for stochastic time-dependent problems with a strong hyperbolic behavior.

\section*{Acknowledgments}
This research was supported by the Natural Sciences and Engineering Research Council of Canada and Hydro Qu\'{e}bec. Their financial support is gratefully acknowledged.
 
\clearpage
\bibliography{mybibfile}
\bibliographystyle{elsarticle-num}
\end{document}